\documentclass[11pt]{article}

\usepackage[margin=1in]{geometry}  
\usepackage{setspace}              
\usepackage{amsmath,amssymb}       
\usepackage{graphicx}              
\usepackage{multirow}%
\usepackage{amsmath,amssymb,amsfonts}%
\usepackage{amsthm}%
\usepackage{mathrsfs}%
\usepackage[title]{appendix}%
\usepackage{xcolor}%
\usepackage{textcomp}%
\usepackage{ragged2e}
\usepackage{manyfoot}%
\usepackage{booktabs}%
\usepackage{arydshln}
\usepackage{algorithm}
\usepackage{algpseudocode} 
\algrenewcommand\algorithmicrequire{\textbf{Input:}}
\algrenewcommand\algorithmicensure{\textbf{Output:}}
\usepackage{listings}%
\usepackage{enumitem}
\usepackage{makecell}
\usepackage{subcaption}

\usepackage{tabularx}
\usepackage{tablefootnote}
\usepackage{endnotes}

\usepackage{threeparttable}
\usepackage{url}
\usepackage[utf8]{inputenc}

\newtheorem{theorem}{Theorem}%
\newtheorem{proposition}[theorem]{Proposition}%
\newtheorem{lemma}[theorem]{Lemma}%
\newtheorem{example}{Example}%
\newtheorem{remark}{Remark}%
\newtheorem{problem}{Problem}

\newtheorem{definition}{Definition}%
\newtheorem{assumption}{Assumption}%

\numberwithin{equation}{section}

\begin{document}

\title{Robust Out-of-Distribution Stochastic Optimization}
\author{
    \large Xianyu Li \\
    \footnotesize Department of Automation, Tsinghua University \\
    \footnotesize li-xy22@mails.tsinghua.edu.cn \\[2ex]
    \large Huan Xu \\
    \footnotesize Antai College of Economics \& Management, Shanghai Jiao Tong University \\
    \footnotesize xuhuan\_antai@sjtu.edu.cn \\[2ex]
    \large Xiaolin Huang \\
    \footnotesize Institute of Image Processing and Pattern Recognition, Shanghai Jiao Tong University \\
    \footnotesize xiaolinhuang@sjtu.edu.cn\\[2ex]
    \large Chao Shang\\
    \footnotesize Department of Automation, Tsinghua University \\
    \footnotesize c-shang@tsinghua.edu.cn}

\date{}
\maketitle

\begin{abstract}
    Data-driven decision-making under uncertainty typically presumes the attainability of historical data from an unknown target probability distribution. However, one may have no access to any data from the target distribution prior to decision-making. To address this challenge, we propose \textit{robust out-of-distribution stochastic optimization} (RooD-SO), a novel data-driven framework that effectively utilizes relevant data distributions for robust decision-making under unseen distributions. A central premise of our RooD-SO framework is that all data distributions are randomly generated from a meta-distribution over distributions. To describe uncertainty in distribution generation, we propose to learn a data-driven uncertainty set in a reproducing kernel Hilbert space (RKHS) from relevant data distributions, with adjustable conservatism. We then incorporate this set into a min-max stochastic program to make robust decisions. Notably, under randomness in distribution generation, we establish rigorous \emph{out-of-distribution generalization guarantees} for the uncertainty set as well as the solution. To effectively resolve the semi-infinite problem, we propose a safe constraint-tightening approximation and design an acceleration strategy based on Nystr\"om approximation. Extensive numerical experiments on multi-item newsvendor and portfolio optimization demonstrate the superior out-of-distribution performance of our decision-making framework under unseen data distributions, even when only a small or moderate number of relevant sources are available.
\end{abstract}

\section{Introduction}\label{sec:intro}
Many operational decisions recur across related but distinct uncertain environments. For a single uncertain environment, the decision problem is commonly formulated as a stochastic program (SP)
\begin{equation} \label{eq:SP}
\min_{\mathbf{x}\in\mathcal{X}}~\mathbb{E}_{\boldsymbol{\xi}\sim\mathbb{P}}\{f(\mathbf{x},\boldsymbol{\xi})\},
\end{equation}
where \(\mathbf{x}\in\mathcal{X}\subseteq\mathbb{R}^d\) denotes the decision variable, \(\boldsymbol{\xi}\in\Xi\) represents the uncertainty, \(f:\mathcal{X}\times\Xi\rightarrow\mathbb{R}\) is the cost function, and \(\mathbb{P}\in\mathcal{P}(\Xi)\) denotes the distribution governing the uncertainty. Historical observations from this environment are typically assumed to be accessible so that sample average approximation (SAA) can be applied to construct a data-driven approximation to problem~\eqref{eq:SP}. However, when decisions concern a new environment, such observations may be unavailable, and decision-makers have to instead rely on data from related environments with potentially different underlying distributions.

\begin{figure}[!ht]
    \centering
    \includegraphics[width=\linewidth]{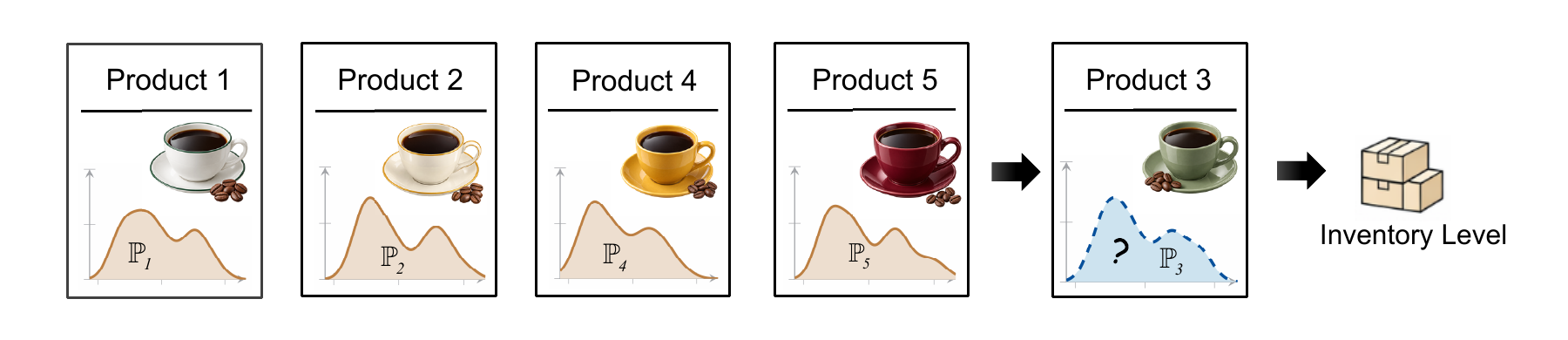}
    \caption{Illustrative example: inventory planning for a new coffee product.}
    \label{fig:top}

    \par\smallskip
    \begin{minipage}{\linewidth}
    \small
    \justifying
    \noindent
    \emph{Notes.} This illustration is based on real sales data for regular-size coffee from the \emph{Coffee Shop Sales} dataset used in our numerical experiments (Section~\ref{sec:experiment_coffee_sales}). The estimated distributions exhibit product-specific differences while sharing similar overall shapes. Product 3 is held out as the target product, and Products 1, 2, 4, and 5 serve as related source products. The target-product curve is shown for illustration; its sales data are withheld during decision-making.
    \end{minipage}
\end{figure}

Such situations are not rare in practice~\cite{lin2024temporally}. A simple example is inventory planning illustrated in Figure~\ref{fig:top}, in which a coffee retailer routinely introduces new products to refresh its assortment and encourage store visits~\cite{bernstein2017dynamic}. Before each launch, the retailer seeks an initial inventory level that balances the costs of overstocking and understocking, despite having no sales history for the new product. Fortunately, data are available from existing relevant coffee products, but each has its own demand distribution. As shown in Figure~\ref{fig:top}, the distributions estimated from real sales data share broad features while exhibiting product-specific differences. This empirical pattern suggests both commonality and heterogeneity among the underlying demand distributions. It motivates seeking an inventory decision that is robust across related products and performs reliably for a new product. Importantly, this issue extends beyond commercial domains and also frequently arises in public health~\cite{ferguson2020report} and medical services~\cite{villar2015multi, awasthi2022distributionally, liu2024multi}. For instance, for a newly established hospital, local health records are not available, and one has to rely on data from other hospitals with similar geographic and demographic profiles to inform resource planning. More formally, this can be posed as the following out-of-distribution decision-making problem:
\begin{problem}
    For an unseen target distribution \(\mathbb{P}\), how can we leverage other relevant datasets \(\mathcal{D}_1,\cdots,\mathcal{D}_M\), each including $n_i$ samples drawn from a distribution \(\mathbb{P}_i\) that is not identical with but similar to \(\mathbb{P}\), to solve problem~\eqref{eq:SP} and inform robust decision-making under \(\mathbb{P}\)?
    \label{Q}
\end{problem}

In the absence of precise knowledge of probability distributions, distributionally robust optimization (DRO) offers a promising data-centric decision-making paradigm. It employs a set of distributions \(\mathcal{U}\) to describe the uncertainty of the target distribution, and optimizes the objective of~\eqref{eq:SP} in the worst case:
\begin{equation} \label{eq:DRO}
\min_{\mathbf{x}\in\mathcal{X}}\;\sup_{\mathbb{P}\in\mathcal{U}}\;\mathbb{E}_{\boldsymbol{\xi}\sim\mathbb{P}}\{f(\mathbf{x},\boldsymbol{\xi})\}.
\end{equation}
Given a dataset drawn from the target distribution, \(\mathcal{U}\) shall be designed to include distributions that plausibly resemble the empirical distribution \(\hat{\mathbb{P}}\). A prevalent choice is to define \(\mathcal{U}\) as a ball centered at the empirical distribution \(\hat{\mathbb{P}}\) with respect to some probability metric, such as divergence~\cite{ben2013robust, duchi2021learning}, Wasserstein distance~\cite{mohajerin2018data}, or maximum mean discrepancy (MMD)~\cite{staib2019distributionally}. However, this is not applicable under the setup of Problem~\ref{Q} where multi-source data $\{ \mathcal{D}_i \}$ do not provide a nominal distribution and the target distribution $\mathbb{P}$ is completely unobserved. This has stimulated recent extensions of DRO models aggregating data from relevant heterogeneous distributions, which can be classified into three categories in terms of uncertainty representation.

\begin{figure}[!ht]
    \begin{subfigure}[b]{0.296727\linewidth}
        \centering
        \includegraphics[width=\linewidth]{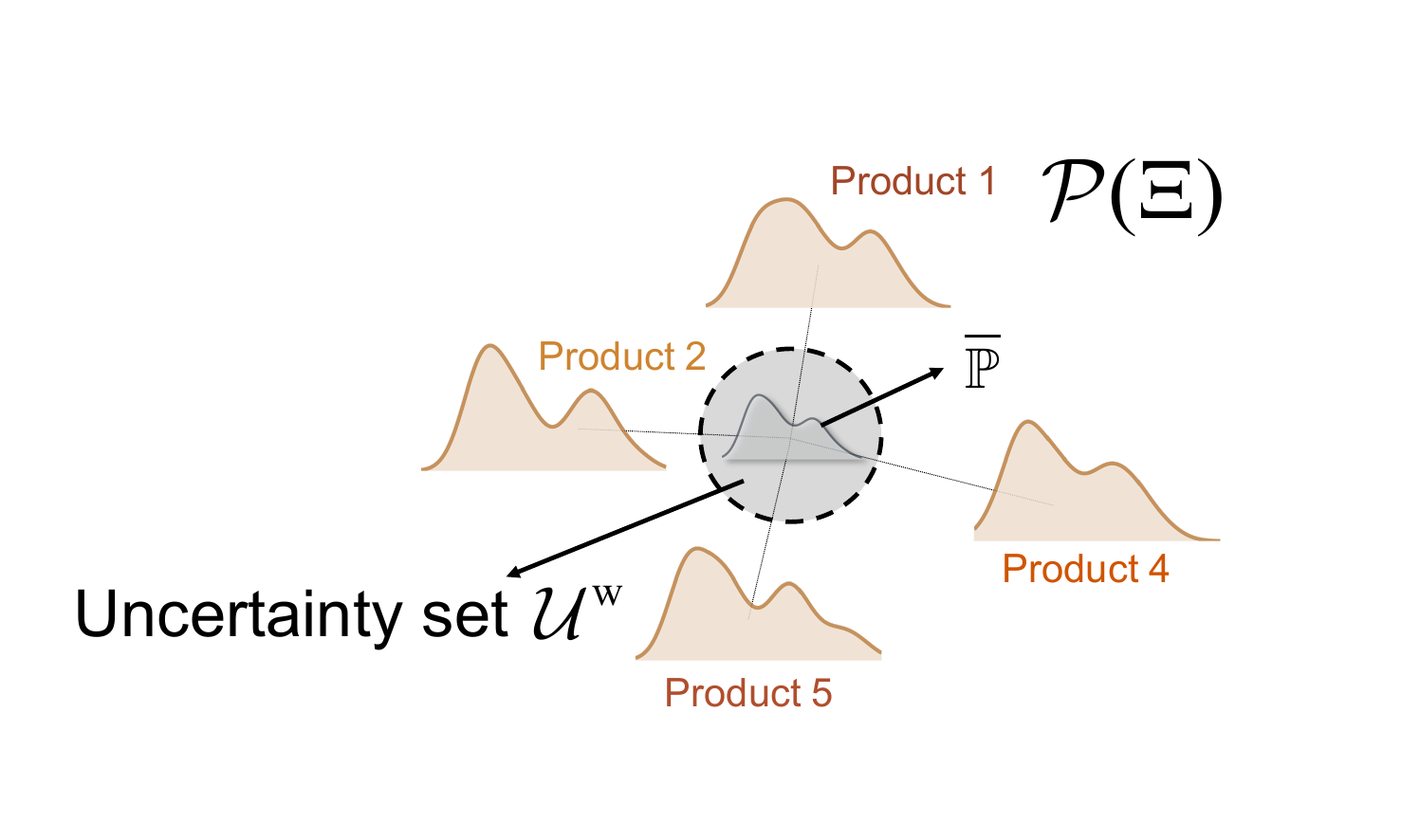}
        \caption{Distribution averaging.}
        \label{fig:b0}
      \end{subfigure}
      \hspace{0.02\linewidth}
  \begin{subfigure}[b]{0.356072\linewidth}
    \centering
    \includegraphics[width=\linewidth]{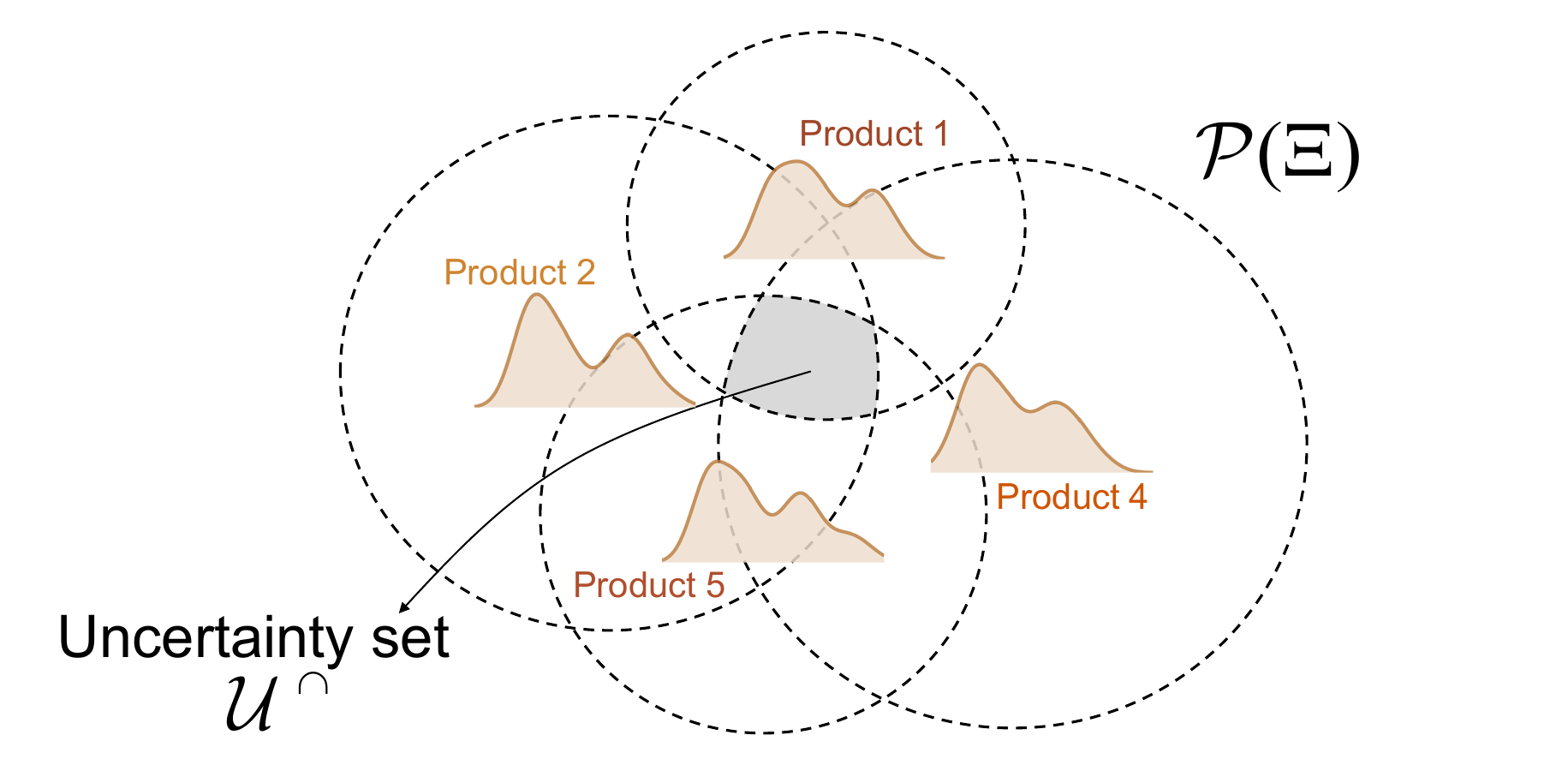}
    \caption{Intersection of Wasserstein balls.}
    \label{fig:b2}
  \end{subfigure}
    \begin{subfigure}[b]{0.296727\linewidth}
    \centering
    \includegraphics[width=\linewidth]{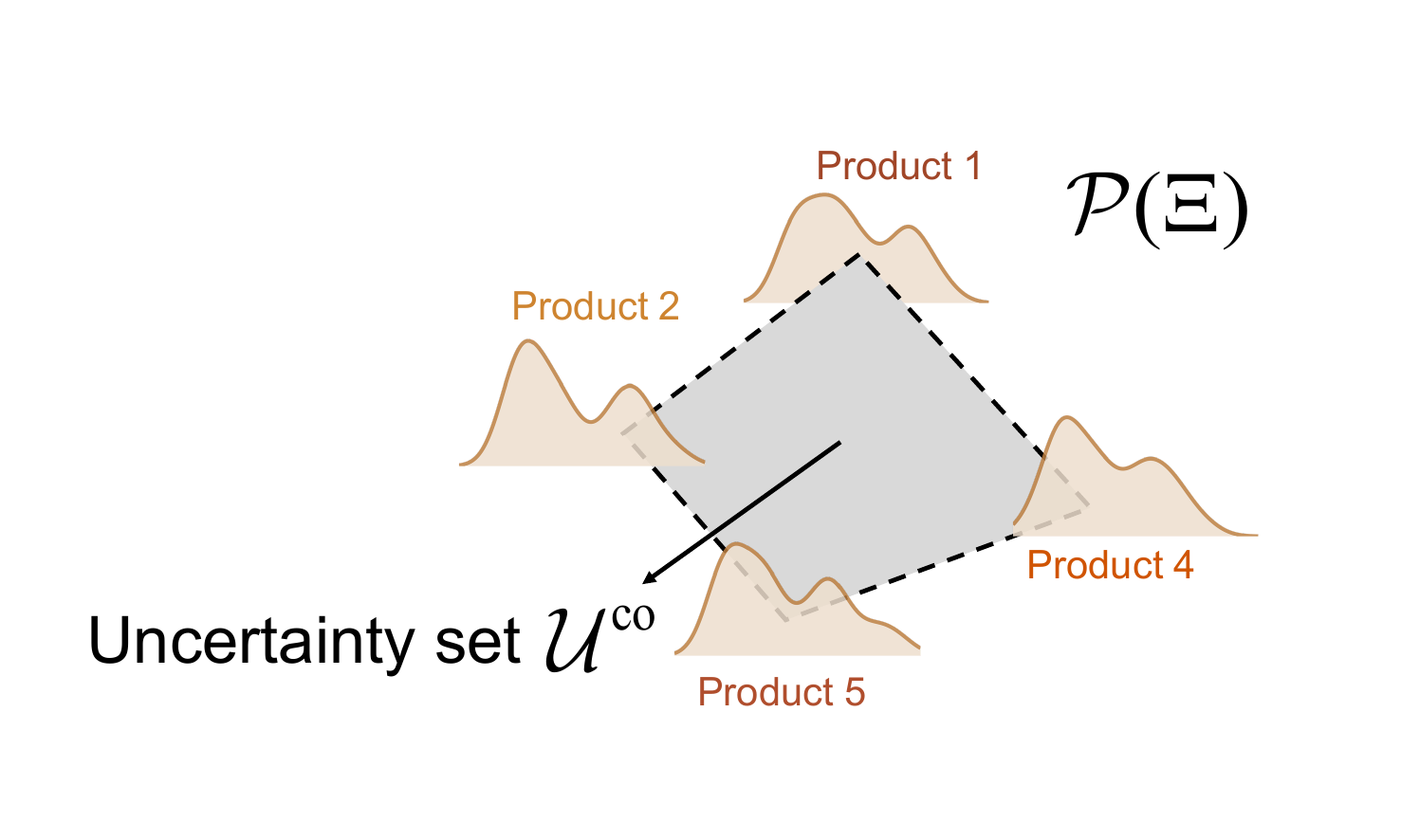}
    \caption{Convex hull of distributions.}
    \label{fig:b1}
  \end{subfigure}
  
  \caption{Schematic of existing methods for aggregating multiple data distributions.}
\end{figure}

\textbf{Distribution averaging via barycenter:} To aggregate multi-source datasets, a simple idea is to average all heterogeneous distributions into a mean component \(\bar{\mathbb{P}}\), as shown in Figure~\ref{fig:b0}. In practice, one can first pool samples from all sources into a single aggregate dataset \(\mathcal{D}^{\cup}:=\mathcal{D}_1\cup\cdots\cup\mathcal{D}_M\), and then use SAA to solve~\eqref{eq:SP} approximately~\cite{gupta2022data, dursun2022data}. This amounts to taking a weighted average of empirical distributions \(\bar{\mathbb{P}}=\sum_{i=1}^M \frac{n_i}{\sum_{i=1}^Mn_i}\hat{\mathbb{P}}_i\), interpretable as a barycenter in the Euclidean space. Rooted in the optimal transport theory, the Wasserstein barycenter (WB) plays a crucial role in machine learning and statistics thanks to its ability to capture structural and geometric features~\cite{srivastava2018scalable, schmitz2018wasserstein, montesuma2021wasserstein, rychener2024wasserstein}. To immunize against the discrepancy between the barycenter \(\bar{\mathbb{P}}\) and the true target distribution, one can follow the rationale of DRO by constructing a Wasserstein ambiguity set \(\mathcal{U}^{\rm w}:=\mathbb{B}_{\epsilon}(\bar{\mathbb{P}})\), i.e., a ball of radius \(\epsilon\) centered at the barycenter \(\bar{\mathbb{P}}\), and solving \eqref{eq:DRO}~\cite{montesuma2021wasserstein}. Crucially, the effectiveness of this framework hinges on the proper definition of the barycenter and a judicious selection of the Wasserstein radius \(\epsilon\).

\textbf{Intersection of Wasserstein balls:} An alternative idea is to assume the unseen target distribution \(\mathbb{P}\) to be concurrently ``close'' to all source distributions and place a Wasserstein ball around each of them. This gives rise to an intersection-based uncertainty set, as illustrated in Figure~\ref{fig:b2}:
\begin{equation} \label{eq:intersection_ambiguity_set}
    \mathcal{U}^{\cap} = \bigcap_{i=1}^M \mathbb{B}_{\epsilon_i}(\hat{\mathbb{P}}_i),
\end{equation}
where the ball radius \(\epsilon_i\) characterizes the possible discrepancy between \(\mathbb{P}_i\) and \(\mathbb{P}\), allowing us to encode prior knowledge about the proximity to each source distribution. The simplest case of \(M=2\) has already been studied in~\cite{awasthi2022distributionally, pmlr-v139-taskesen21a, selvi2024distributionally, wang2024contextual}, where an auxiliary dataset provides an additional source of information about \(\mathbb{P}\) and helps to exclude unrealistic distributions. More recently, \cite{rychener2024wasserstein} presented a more general framework considering \(M\geq 2\) source distributions and derived the convex problem reformulations. To suitably choose \(\epsilon_i\), prior knowledge about the underlying bias is required, which is however non-trivial to acquire in practice. In addition, the computational complexity grows exponentially with $M$, rendering the solution prohibitive even under moderate sample sizes.

\textbf{Convex hull of distributions:} By adopting a set-membership perspective, we can use a mixture of heterogeneous source distributions to capture the uncertainty of the target distribution:
\begin{equation} \label{eq:Group_DRO_ambiguity_set}
    \mathcal{U}^{\rm co}=\left \{\mathbb{P}\in\mathcal{P}(\Xi)\,\middle|\,\boldsymbol{\alpha}\in\Delta^M,~ \mathbb{P}=\sum_{i=1}^M\alpha_i\hat{\mathbb{P}}_i\right \}.
\end{equation}
As depicted in Figure~\ref{fig:b1}, the uncertainty set \(\mathcal{U}^{\rm co}\), including all possible mixtures, is essentially a convex hull in the space $\mathcal{P}(\Xi)$. The usage of \(\mathcal{U}^{\rm co}\) is prevalent in machine learning problems that admit group DRO formulations~\cite{Sagawa2020Distributionally, carmon2022distributionally, zhang2023stochastic, pmlr-v235-yu24a, guo2026statisticalanalysisconditionalgroup}, e.g., to address fairness concerns~\cite{hashimoto2018fairness, lahoti2020fairness} or to pursue group robustness in federated learning~\cite{mohri2019agnostic, deng2020distributionally}. However, the inflexibility in regulating conservatism of \(\mathcal{U}^{\rm co}\) can lead to overly conservative solutions to the min-max problem~\eqref{eq:DRO}. 

\subsection*{Open Challenge and Our Proposal}

Conceptually, the three DRO extensions adopt different philosophies for modelling $\mathbb{P}$. Both computing a Barycenter and intersecting Wasserstein balls can be viewed as different means of identifying $\mathbb{P}$ via averaging multiple biased distributions, where the former assigns sample-size-dependent weights to eliminate biases of multiple sources and the latter allows for adjustable confidences through user-specified ball radii. Alternatively, the convex hull paradigm takes a deterministic set-membership viewpoint to capture the uncertainty of $\mathbb{P}$. 

\begin{figure}[!ht]
    \centering
    \includegraphics[width=\linewidth]{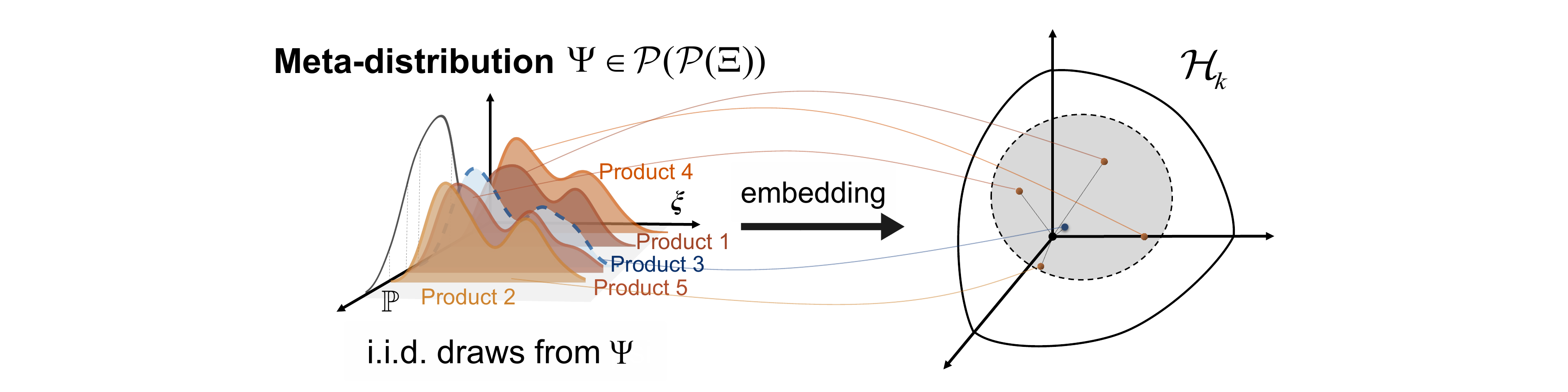}
    \caption{Our proposal: Meta-distributional modeling and embeddings in RKHS.}
    \label{fig:b3}
\end{figure}

Generally speaking, these current endeavors did not put the randomness of the target distribution and source distributions on a stochastic footing. We seek to address this open challenge from a novel stochastic distribution-generating perspective, which meaningfully reflects the uncertainty the decision
will finally face in Problem 1, and design a learning-based decision-making framework called \textit{robust out-of-distribution stochastic optimization} (RooD-SO). Mathematically, we interpret both the unobserved target distribution \(\mathbb{P}\) and observed source distributions \(\{\mathbb{P}_i\}\) as independent realizations of an unknown \textit{meta-distribution} \(\Psi\in\mathcal{P}(\mathcal{P}(\Xi))\), i.e., a distribution over distributions, where \(\mathcal{P}(\mathcal{P}(\Xi))\) is the set of all probability measures defined on the sample space \(\mathcal{P}(\Xi)\). As illustrated in Figure~\ref{fig:b3}, the pair \((\mathbb{P}, \boldsymbol{\xi})\) features a two-layer randomness structure, where one can first draw a distribution \(\mathbb{P}\) from \(\Psi\) and then sample a random realization \(\boldsymbol{\xi}\in\Xi\) from \(\mathbb{P}\). This new hierarchy of generation inspires us to convert the construction of $\mathcal{U}$ into a problem of estimating the support of \(\Psi\) from \(M\) source distributions \(\{\mathbb{P}_i\}_{i=1}^M\). Specifically, with \(\{\mathbb{P}_i\}_{i=1}^M\) projected into a reproducing kernel Hilbert space (RKHS), similarities between distributions can be evaluated through inner products in the RKHS. The support set induced by one-class clustering in RKHS yields a data-driven uncertainty set \(\mathcal{U}(\mathbb{P}_1,\cdots,\mathbb{P}_M)\) that learns possible realizations of distributions generated by \(\Psi\). In doing so, it captures not only the structural similarities among \(\mathbb{P}_1,\cdots,\mathbb{P}_M\) but also the underlying distributional uncertainty, with probabilistic coverage guarantees for other unseen distributions drawn from $\Psi$. This turns out to secure robust decision-making with stable out-of-distribution performance when resolving the min-max problem~\eqref{eq:DRO}. Our technical contributions are summarized as follows.

\begin{itemize}
    \item We embed source distributions \(\{\mathbb{P}_i\}_{i=1}^M\) into an RKHS using their kernel mean embeddings and utilize support measure clustering (SMC)~\cite{muandet2013one} to learn a meaningful support estimate of \(\Psi\) as a ball in the RKHS. The conservatism of the resulting uncertainty set can be flexibly regulated by two tuning parameters. Since the ground truth of all source distributions is unknown and their empirical substitutes \(\{\hat{\mathbb{P}}_i\}_{i=1}^M\) have to be used, we establish rigorous finite-sample bounds to characterize the approximation errors.
    \item We establish novel out-of-distribution probabilistic guarantees under the two-layer randomness arising from finitely many source distributions and finite observations within each source. In particular, we prove that the learned uncertainty set covers an unseen target distribution with high probability and that the resulting optimal value of the downstream min-max problem provides a probabilistic performance certificate for the decision of RooD-SO. We further establish the asymptotic validity of these guarantees under certain sampling conditions.
    \item We present a novel constraint-tightening formulation to approximately solve the semi-infinite dual problem of RooD-SO, with feasibility preserved and consistency established. We also devise an acceleration solution strategy based on Nystr\"om approximation, which yields a bounded loss of optimality and preserves the desired probabilistic performance guarantees.
    \item To evaluate the efficacy of our RooD-SO, we make comprehensive empirical comparisons against known approaches on both multi-item newsvendor and portfolio optimization problems. Numerical experiments on synthetic and real data showcase that our framework exhibits desirable out-of-distribution generalization performance and secures significant improvements over other schemes in unseen uncertain environments with the aid of a small to moderate number of sources.
\end{itemize}

\paragraph{Structure.} The layout of this paper is structured as follows. Section~\ref{sec:method} elaborates on the proposed framework, covering preliminaries, the problem setup and assumptions, uncertainty set construction, and the resultant optimization problem. We also provide a comprehensive theoretical analysis that accommodates the two-layer randomness. Subsequently, Section~\ref{sec:solution} devises an efficient solution methodology for our proposed framework. Section~\ref{sec:experiments} reports numerical results evaluated using synthetic data, real-world coffee shop sales data, and ETF data. Finally, Section~\ref{sec:conclusion} closes the paper with concluding remarks. All technical proofs are relegated to Appendix~\ref{Appendix:proofs}. 

\paragraph{Notation.} Let \((\Xi,\mathcal{A})\) be a measurable space, where \(\Xi\) is a non-empty observation space and \(\mathcal{A}\) is the associated \(\sigma\)-algebra. Let \(\mathbb{P}\) denote a probability measure on \((\Xi,\mathcal{A})\), and \(\mathcal{P}(\Xi)\) be the set of all such probability measures. For any distribution \(\mathbb{P}\), we write \(\mathrm{supp}(\mathbb{P})\) for its support. We write \(\boldsymbol{\xi}_1,\cdots,\boldsymbol{\xi}_n\overset{\rm i.i.d.}{\sim} \mathbb{P}\) to indicate that \(\boldsymbol{\xi}_1,\cdots,\boldsymbol{\xi}_n\) are independent and identically distributed (i.i.d.) according to \(\mathbb{P}\). We write \(\delta_{\boldsymbol{\xi}}\) for the Dirac measure at \(\boldsymbol{\xi}\in\Xi\). We use $\otimes$ to denote the product of probability measures and $\bigotimes$ for its indexed form. In particular, $\mathbb{P}^{m}:=\bigotimes_{k=1}^{m}\mathbb{P}$ denotes the $m$-fold product measure of $\mathbb{P}$ on $\Xi^m$. A function \(f:\mathcal{X}\rightarrow[-\infty,\infty]\) is upper semicontinuous on \(\mathcal{X}\) if, for every \(x_0\in\mathcal{X}\), \(\limsup_{x\rightarrow x_0} f(x)\leq f(x_0)\). Likewise, \(f\) is lower semicontinuous if, for every \(x_0\in\mathcal{X}\), \(\liminf_{x\rightarrow x_0} f(x)\geq f(x_0)\). We use $\operatorname{ri}(\mathcal{C})$ to denote the interior of a convex set $\mathcal{C}$ relative to its affine hull, with the topology inherited from the ambient space. We denote by \(\boldsymbol{1}_J\in\mathbb{R}^J\) the vector of length \(J\) whose entries are all equal to one. We use \(x\lesssim n\) to denote \(x\leq Dn\) for some universal constant \(D>0\). For any event \(A\), \({\rm Pr}_\mathbb{P}\{A\}\) denotes the probability of \(A\) under the probability measure \(\mathbb{P}\). We denote the set of positive integers by \(\mathbb{Z}_{>0}\). The notation \(\mathcal{W}_p(\mathbb{P}_1,\mathbb{P}_2)\) stands for the \(p\)-Wasserstein distance between \(\mathbb{P}_1\) and \(\mathbb{P}_2\). \({\rm Unif}([a,b])\) denotes the continuous uniform distribution on the closed interval \([a,b]\) with \(a<b\). For any symmetric matrix \(\mathbf{K}\in\mathbb{R}^{M\times M}\), we use \(\lambda_{\max}(\mathbf{K})\) and \(\lambda_{\min}(\mathbf{K})\) to denote its largest and smallest eigenvalues, respectively. The notation \({\rm diag} (\mathbf{K})\in\mathbb{R}^M\) denotes the vector consisting of the diagonal elements of \(\mathbf{K}\). We use $\|\cdot\|_2$ for the Euclidean norm of
vectors and the induced spectral norm of matrices. Furthermore, the spectral norm of a positive semidefinite matrix $\mathbf{K}$ is defined as \(\|\mathbf{K}\|_2:=\sup_{\boldsymbol{\xi}\neq 0}\|\mathbf{K}\boldsymbol{\xi}\|_2/\|\boldsymbol{\xi}\|_2=\lambda_{\max}(\mathbf{K})\). The superscript \(\dagger\) denotes the Moore-Penrose pseudoinverse of a matrix. For a non-empty closed convex set \(\mathcal{M}\subset\mathcal{H}_k\), its support function is defined as \(h^*_{\mathcal{M}}(g):=\sup_{\mu\in \mathcal{M}}\langle g,\mu \rangle_{\mathcal{H}_k}\). For any positive integer \(N\), we use \([N]\) to denote the index set \(\{1,\cdots,N\}\). For a set of points $\mathcal{D} = \{\xi_1,\dots,\xi_{n}\}$, its cardinality is denoted by $|\mathcal{D}| = n$. We denote by $\Delta^M := \{\boldsymbol{\alpha}\in\mathbb{R}^M : \boldsymbol{\alpha}\ge 0,\ \mathbf{1}^{\top}\boldsymbol{\alpha}=1\}$ the standard simplex in $\mathbb{R}^M$. For a point \(\mathbf x\in\mathbb R^m\) and a nonempty set \(\mathcal A\subseteq\mathbb R^m\), we write \(\operatorname{dist}(\mathbf x,\mathcal A):=\inf_{\mathbf y\in \mathcal A}\|\mathbf x-\mathbf y\|_2\).

\section{Proposed methodology}
\label{sec:method}

In this section, we formally present our proposed RooD-SO framework. To set the stage, we first discuss a new meta-distribution perspective and then formalize our problem setup. Given \(M\) different source distributions \(\{\mathbb{P}_i\}_{i=1}^M\), we assume that they are independent realizations from an unknown meta-distribution \(\Psi \in\mathcal{P}(\mathcal{P}(\Xi))\), and so is the target distribution \(\mathbb{P}\). That is,
\begin{equation}
    \mathbb{P}_1, \dots, \mathbb{P}_M, \mathbb{P} \overset{\rm i.i.d.}{\sim} \Psi.
\end{equation}
We seek to estimate \(\mathrm{supp}(\Psi)\) in the probability space $\mathcal{P}(\Xi)$ by capturing structural similarities among \(\{\mathbb{P}_i\}_{i=1}^M\). The learned support yields an uncertainty set \(\mathcal{U}(\mathbb{P}_1,\cdots,\mathbb{P}_M)\) enclosing adverse but plausible realizations of the target distribution $\mathbb{P}$ on which we want our decision to be well-performing. This is achieved by formulating and solving the following min-max problem
\begin{equation} \label{eq:main_data_driven}
    \min_{\mathbf{x} \in \mathcal{X}} \; \sup_{\mathbb{P} \in \mathcal{U}(\mathbb{P}_1,\ldots,\mathbb{P}_M)} 
   \; \mathbb{E}_{\boldsymbol{\xi} \sim \mathbb{P}} \left\{ f(\mathbf{x},\boldsymbol{\xi}) \right\}.
\end{equation} 
If \(\mathcal{U}(\mathbb{P}_1,\cdots,\mathbb{P}_M)\) can adequately enclose the true target distribution $\mathbb{P}$, we can make a data-driven decision \(\mathbf{x}\in\mathcal{X} \) with robust out-of-distribution performance by hedging against the least-favorable plausible distribution in \(\mathcal{U}(\mathbb{P}_1,\cdots,\mathbb{P}_M)\). Here we refer to \(\mathcal{U}(\mathbb{P}_1,\cdots,\mathbb{P}_M)\) as an \textit{uncertainty set} to emphasize that it captures distributional uncertainty induced by the meta-distribution \(\Psi\), in contrast to the \textit{ambiguity set} enclosing an existent but unknown data-generating distribution in conventional DRO. Throughout this article we make the following assumption, which ensures the structural properties necessary for the existence and duality results to be established in the sequel.
\begin{assumption}\label{ass:loss_function}
    \(\Xi\) is a compact metric space and \(\mathcal{X}\subseteq \mathbb{R}^{d}\) is a closed convex set. The cost function \(f(\mathbf{x},\cdot):\Xi\rightarrow\mathbb{R}\) is continuous for each \(\mathbf{x}\in\mathcal{X}\) and \(f(\cdot,\boldsymbol{\xi}):\mathcal{X}\rightarrow\mathbb{R}\) is convex for every \(\boldsymbol{\xi}\in\Xi\). 
\end{assumption}
Assumption~\ref{ass:loss_function} is mild and holds for a broad class of cost functions of practical interest, including linear and piecewise-linear convex losses, as well as quadratic losses commonly arising in risk and regression models. 

As discussed, constructing the uncertainty set $\mathcal{U}(\mathbb{P}_1,\cdots, \mathbb{P}_M)$ in~\eqref{eq:main_data_driven} substantially entails estimating $\mathrm{supp}(\Psi)$ over the probability space \(\mathcal{P}(\Xi)\), which is non-trivial since \(\mathcal{P}(\Xi)\) is an infinite-dimensional space of measures subject to nonnegativity and unit-mass constraints. To combat this, we employ the kernel mean embedding~\cite{gretton2006kernel} to map distributions in \(\mathcal{P}(\Xi)\) into an abstract RKHS. This projection enables us to utilize the Hilbert space geometry and full arsenal of kernel methods~\cite{Hofmann2008kernelmethods} to operate efficiently on distributions. In the remainder of this section, we shall first review some basic knowledge of RKHS and kernel mean embedding (Section~\ref{sec:preliminaries}). Then, we investigate an idealized setup with full knowledge of \(\{\mathbb{P}_i\}_{i=1}^M\) to construct the uncertainty set and formulate the associated min-max problem (Section~\ref{sec:SMC_ambiguity_set}), and then turn to the finite-sample case and establish concentration inequalities for the induced errors (Section~\ref{sec:data_driven_setup}). Finally, we establish rigorous out-of-distribution guarantees for the uncertainty set and the solution under two-layer randomness (Section~\ref{sec:theoretical_analysis}).


\subsection{Distribution embeddings in RKHS}
\label{sec:preliminaries}
Next, we provide an introductory guide to RKHS and kernel mean embedding. A detailed exposition can be found in~\cite{aronszajn1950theory, muandet2017kernel}. Given a symmetric and positive definite kernel \(k:\Xi\times \Xi\rightarrow\mathbb{R}\), there exists a Hilbert space \(\mathcal{H}_k\) and a feature map \(\phi:\Xi\rightarrow\mathcal{H}_k\) such that \(k(\boldsymbol{\xi},\boldsymbol{\xi}')=\langle\phi(\boldsymbol{\xi}),\phi(\boldsymbol{\xi}')\rangle_{\mathcal{H}_k}\) for all \(\boldsymbol{\xi},\boldsymbol{\xi}'\in\Xi\).\endnote{For readers less familiar with kernel methods, a kernel \(k(\boldsymbol{\xi},\boldsymbol{\xi}')\) may be interpreted as a similarity score between two inputs \(\boldsymbol{\xi}\) and \(\boldsymbol{\xi}'\). A key property is that such a kernel defines an associated RKHS, within which nonlinear functional relationships can often be handled using linear operations in an implicit feature space. This avoids the need to explicitly construct that feature space and provides a convenient way to model complex nonlinear structures.} This space \(\mathcal{H}_k\), an RKHS associated with \(k\), consists of functions \(f:\Xi\rightarrow\mathbb{R}\) wherein the evaluation functional is continuous, yielding the following reproducing property:
\begin{equation}
    f(\boldsymbol{\xi})=\langle f,\phi(\boldsymbol{\xi}) \rangle_{\mathcal{H}_k}, \quad \forall f\in\mathcal{H}_k,~\boldsymbol{\xi}\in\Xi.
\end{equation}
This property enables pointwise evaluation of functions in \(\mathcal{H}_k\) through inner products. Additionally, RKHSs are linear such that any linear combination of functions still lies in the space and the inner product distributes over addition in each argument.

Thanks to the reproducing property and the linearity of RKHSs, one can naturally extend pointwise evaluation to a distribution \(\mathbb{P}\in\mathcal{P}(\Xi)\) by taking expectations of \(\phi(\boldsymbol{\xi})\) under \(\mathbb{P}\). Under the (sufficient) integrability condition \(\mathbb{E}_{\boldsymbol{\xi}\sim\mathbb{P}}\{\sqrt{k(\boldsymbol{\xi},\boldsymbol{\xi})}\} < \infty\), the Bochner integral 
\(\mathbb{E}_{\boldsymbol{\xi} \sim \mathbb{P}}\{\phi(\boldsymbol{\xi})\}\) 
is well-defined in \(\mathcal{H}_k\)~\cite{smola2007hilbert}. Then, for any \( f \in \mathcal{H}_k \) and distribution \( \mathbb{P} \in \mathcal{P}(\Xi) \), we have
\begin{equation}
    \mathbb{E}_{\boldsymbol{\xi} \sim \mathbb{P}}\left\{f(\boldsymbol{\xi})\right\} = \mathbb{E}_{\boldsymbol{\xi} \sim \mathbb{P}}\left\{\langle f, \phi(\boldsymbol{\xi}) \rangle_{\mathcal{H}_k}\right\} = \left\langle f, \mathbb{E}_{\boldsymbol{\xi} \sim \mathbb{P}}\left\{\phi(\boldsymbol{\xi})\right\} \right\rangle_{\mathcal{H}_k},
\end{equation}
which enables us to define the kernel mean embedding of a probability distribution \(\mathbb{P}\). 
\begin{definition}[Kernel mean embedding \cite{smola2007hilbert, gretton2006kernel}] \label{defn:kernel_mean_embedding}
    Given a kernel function \(k:\Xi\times\Xi\rightarrow\mathbb{R}\) , the kernel mean embedding of a probability measure \(\mathbb{P}\in\mathcal{P}(\Xi)\) is 
    \begin{equation} \label{eq:kernel_mean_embedding}
        \begin{aligned}
            \mu_{\mathbb{P}}:=\int_{\Xi}\phi(\boldsymbol{\xi})\mathrm{d}\mathbb{P}(\boldsymbol{\xi})\in\mathcal{H}_k.
        \end{aligned}
    \end{equation}
\end{definition}
Intuitively, the kernel mean embedding maps any distribution to a single point \(\mu_{\mathbb{P}}\) in the RKHS, which can be explained as a generalized moment vector carrying rich statistical information of \(\mathbb{P}\)~\cite{zhu2021kernel}. By applying Fubini's theorem and the reproducing property, we obtain
\begin{equation} \label{eq:inner_product_distribution}
    \langle \mu_{\mathbb{P}_1},\mu_{\mathbb{P}_2} \rangle_{\mathcal{H}_k}=\int\int \langle \phi(\boldsymbol{\xi}),\phi(\boldsymbol{\zeta}) \rangle_{\mathcal{H}_k}\mathrm{d}\mathbb{P}_1(\boldsymbol{\xi})\mathrm{d}\mathbb{P}_2(\boldsymbol{\zeta}) =\int\int k(\boldsymbol{\xi},\boldsymbol{\zeta})\mathrm{d}\mathbb{P}_1(\boldsymbol{\xi})\mathrm{d}\mathbb{P}_2(\boldsymbol{\zeta}),
\end{equation}
which allows for distribution-level evaluations. In particular, the induced RKHS distance is given by $d_k(\mathbb P_1,\mathbb P_2):=\|\mu_{\mathbb{P}_1}-\mu_{\mathbb{P}_2}\|_{\mathcal{H}_k}=\sqrt{\langle\mu_{\mathbb{P}_1},\mu_{\mathbb{P}_1}\rangle_{\mathcal{H}_k}+\langle\mu_{\mathbb{P}_2},\mu_{\mathbb{P}_2}\rangle_{\mathcal{H}_k}-2\langle\mu_{\mathbb{P}_1},\mu_{\mathbb{P}_2}\rangle_{\mathcal{H}_k}}$, which is also known as MMD and naturally defines a metric for evaluating discrepancy between two distributions. Moreover, the embeddings shall encode all information about distributions, which is ensured by the so-called characteristic kernels.

\begin{theorem}[Characteristic kernels: injectivity and a sufficient condition~\cite{ sriperumbudur2010hilbert}]
\label{thm:injectivity}
Let \(k:\Xi\times\Xi\to\mathbb{R}\) be a positive semidefinite kernel on \(\Xi\) with RKHS \(\mathcal{H}_k\), and let
\(\mu_{\mathbb{P}}\in\mathcal{H}_k\) denote the kernel mean embedding for \(\mathbb{P}\in\mathcal{P}(\Xi)\).
\begin{enumerate}[label=(\roman*), leftmargin=*]
    \item Kernel \(k(\cdot,\cdot)\) is said to be characteristic if and only if the map \(\mathbb{P}\mapsto\mu_{\mathbb{P}}\) is injective; equivalently, for any \( \mathbb{P}, \mathbb{Q} \in \mathcal{P}(\Xi) \), \( \mu_{\mathbb{P}} = \mu_{\mathbb{Q}} \) implies \( \mathbb{P} = \mathbb{Q} \).
    \item If kernel \(k(\cdot,\cdot)\) is integrally strictly positive definite on \(\Xi\), then \(k(\cdot,\cdot)\) is characteristic to \(\mathcal{P}(\Xi)\).
\end{enumerate}
\end{theorem}



Many popular kernel functions, like the RBF (Radial Basis Function), Laplace, and Mat\'ern kernels, are integrally strictly positive definite and therefore characteristic~\cite{sriperumbudur2010hilbert}. For these kernels, Theorem~\ref{thm:injectivity} establishes the injectivity of the kernel mean embeddings, implying that distributions can be uniquely represented in RKHS. Consequently, \(\|\mu_{\mathbb{P}_1}-\mu_{\mathbb{P}_2}\|_{\mathcal{H}_k}\) becomes a valid metric on the space of distributions, which is known to be upper bounded by some probability metrics such as the total variation and the Wasserstein metrics~\cite{sriperumbudur2010hilbert}. Thus, distributions that are close in these conventional metrics manifest adjacency of their embeddings in \(\mathcal{H}_k\), thereby suggesting continuity of \(\mathbb{P}\mapsto\mu_{\mathbb{P}}\). These merits have stimulated a series of applications of kernel mean embedding, e.g., two-sample homogeneity tests \cite{gretton2012kernel, doran2014permutation} and group anomaly detection \cite{muandet2013one}. Back to our problem setting, we make the following assumption throughout this article.

\begin{assumption}\label{ass:kernel}
The kernel \(k:\Xi\times\Xi\rightarrow\mathbb{R}\) is continuous, symmetric, and integrally strictly positive definite. Moreover, it is uniformly bounded on the diagonal, i.e., \(\sup_{\boldsymbol{\xi}\in\Xi} k(\boldsymbol{\xi},\boldsymbol{\xi}) = C < \infty\).
\end{assumption}

Assumption~\ref{ass:kernel} imposes only mild regularity conditions on the kernel, which hold for a broad class of kernel functions including the RBF, Laplace, and Mat\'ern kernels. Consequently, given \(M\) source distributions \(\{\mathbb{P}_i\}_{i=1}^M\), we can leverage kernel mean embedding defined in~\eqref{eq:kernel_mean_embedding} to embed them into an RKHS \(\mathcal{H}_k\) without loss of information, yielding the embeddings \(\{\mu_{\mathbb{P}_i}\}_{i=1}^M\).

\subsection{Uncertainty set learning}
\label{sec:SMC_ambiguity_set}
In this subsection, we describe a new learning-based approach to uncertainty set construction from a finite number of heterogeneous data distributions. Drawing ideas from SMC~\cite{muandet2013one}, we characterize possible realizations of distributions by enclosing their kernel mean embeddings with a ball in RKHS. This yields an uncertainty set capturing structural similarity among distributions and encodes support information of the meta-distribution \(\Psi\). To ease understanding, we assume throughout this subsection that the ground truth of \(M\) source distributions \(\{\mathbb{P}_i\}_{i=1}^M\) is exactly known. Then the support of \(\Psi\) can be estimated as a ball in an RKHS by solving the following one-class clustering problem~\cite{muandet2013one}:
\begin{equation}
    \begin{aligned}
        \min_{\mu_c,R,\mathbf{s}}~&R^2+\frac{1}{M\nu}\sum_{i=1}^Ms_i&\\
        \text{s.t.}~~&\|\mu_{\mathbb{P}_i}-\mu_c\|_{\mathcal{H}_k}^2\leq R^2+s_i,&\quad  i\in[M],\\
        &s_i \geq 0, &\quad i\in[M],\\
    \end{aligned}
    \label{eq:soft_minimum_ball}
\end{equation}
where \(\mu_c\) and \(R\) denote the center and the radius of the ball, respectively. The nonnegative slacks \(\{s_i\}_{i=1}^M\) in~\eqref{eq:soft_minimum_ball} allow a fraction of $M$ distributions to be excluded, where each \(s_i\) quantifies the amount by which \(\mu_{\mathbb{P}_i}\) lies outside the ball. The objective in~\eqref{eq:soft_minimum_ball} trades off minimizing the ball volume against penalizing violations from exterior distributions, where the hyperparameter \(\nu\in(0,1)\) is used to regulate this balance. By writing the Lagrangian and invoking the Karush-Kuhn-Tucker (KKT) conditions, we explicitly derive the dual formulation of the primal problem~\eqref{eq:soft_minimum_ball}, which involves a finite number of decision variables \(\boldsymbol{\alpha}\):
\begin{subequations}\label{eq:dual_SMC_quadratic}
    \begin{align}
        \max_{\boldsymbol{\alpha}\in\mathbb{R}^M} ~&
        -\boldsymbol{\alpha}^\top\mathbf{K}\boldsymbol{\alpha}+\boldsymbol{\alpha}^\top{\rm diag}(\mathbf{K})\label{eq:objective_function}\\
        \text{s.t.}~~&0\leq \alpha_i\leq \frac{1}{M\nu},\quad i\in[M]\label{eq:QP_alpha_value_range}\\
        &\sum_{i=1}^M\alpha_i=1,\label{eq:QP_sum_alpha}
    \end{align}
\end{subequations}
where \(\mathbf{K}\in\mathbb{R}^{M\times M}\) is the Gram matrix with entries \(K_{ij} = \langle\mu_{\mathbb{P}_i},\mu_{\mathbb{P}_j}\rangle_{\mathcal{H}_k}\) computed via~\eqref{eq:inner_product_distribution}. Under Assumption~\ref{ass:kernel}, problem \eqref{eq:dual_SMC_quadratic} is a convex quadratic program (QP) that can be efficiently solved via sequential minimization optimization~\cite{SMO}. Here, we emphasize that this dual QP may admit a continuum of optimal solutions because the Gram matrix $\mathbf{K}$ is positive semidefinite and may not be strictly positive definite. However, this general case was not addressed by prior work~\cite{muandet2013one}. Next, we give a formal treatment by characterizing the subtlety of non-unique solutions and presenting a rigorous procedure to determine the optimal primal pair $(\mu_c^*, R^*)$.

\begin{proposition} \label{prop:discussion_dual}
    Let \(\mathcal A^*\) be the set of all optimal solutions to the dual problem~\eqref{eq:dual_SMC_quadratic}. Then the following holds:
    \begin{enumerate}[label=(\alph*), leftmargin=*]
    \item If \eqref{eq:dual_SMC_quadratic} has a unique optimal solution $\boldsymbol{\alpha}^*$, i.e., $\mathcal{A}^*=\{\boldsymbol{\alpha}^*\}$, then an optimal solution to~\eqref{eq:soft_minimum_ball} is given by
    \begin{equation} \label{eq:mu_c}
        \mu_c^* = \sum_{i=1}^M \alpha_i^* \mu_{\mathbb{P}_i}
    \end{equation}
    and
    \begin{equation} \label{eq:R_square}
        (R^*)^2 = \max_{\alpha_{i'}^* < 1/(M\nu),\,i'\in[M]} \left( K_{i'i'} - 2\sum_{i=1}^M \alpha_i^* K_{ii'} + \sum_{i_1=1}^M \sum_{i_2=1}^M \alpha_{i_1}^* \alpha_{i_2}^* K_{i_1i_2} \right).
    \end{equation}
    
    \item If \eqref{eq:dual_SMC_quadratic} admits multiple optimal solutions, i.e., $|\mathcal{A}^*| > 1$, then the primal problem~\eqref{eq:soft_minimum_ball} has a unique optimal pair $(\mu_c^*, R^*)$. Moreover, there always exists an optimal solution $\boldsymbol{\alpha}^* \in \mathcal{A}^*$ with an index $i \in [M]$ satisfying $0 < \alpha_i^* < \frac{1}{M\nu}$. Then given such an $\boldsymbol{\alpha}^*$, the optimal primal pair $(\mu_c^*, R^*)$ can be determined by~\eqref{eq:mu_c} and~\eqref{eq:R_square}.
\end{enumerate}
\end{proposition}

By Proposition~\ref{prop:discussion_dual}, solving the dual problem~\eqref{eq:dual_SMC_quadratic} determines the ball of interest:
\begin{equation} \label{eq:SMC_M}
    \mathcal{B}_{\nu} = \left\{\mu\in\mathcal{H}_k\,\middle|\,\|\mu-\mu_c^*\|_{\mathcal{H}_k}\leq R^*\right\}.
\end{equation}
We then define the uncertainty set \(\mathcal{U}_{\nu}^{\rm SMC}\) as a family of distributions whose embeddings live in \(\mathcal{B}_{\nu}\):
\begin{equation} \label{eq:SMC_ambiguity_set}
    \mathcal{U}_{\nu}^{\text{SMC}} = \left\{\mathbb{P}\in\mathcal{P}(\Xi)\,\middle|\,\mathbb{E}_{\mathbb{P}}\left\{\phi(\boldsymbol{\xi})\right\}=\mu_{\mathbb{P}}\in\mathcal{B}_{\nu}\right\}.
\end{equation}
Thanks to its learning capability, $\mathcal{U}_{\nu}^{\text{SMC}}$ is intended as a meaningful characterization of uncertainty that the decision will actually face. The size of \(\mathcal{U}_{\nu}^{\rm SMC}\) is governed by a pre-specified trade-off parameter \(\nu\). The primal problem \eqref{eq:soft_minimum_ball} implies that a smaller \(\nu\) yields a potentially larger ball in the RKHS and consequently, a more conservative uncertainty set. Indeed, it suffices to select $\nu\in(0,1)$ following a clear interpretation.

\begin{lemma}[\(\nu\)-property~\cite{scholkopf2002learning, muandet2013one}] \label{lemma:nu}
    The parameter \(\nu \in (0,1)\) upper-bounds the fraction of exterior embeddings that lie outside the RKHS ball $\mathcal{B}_{\nu}$ and satisfy $s_i^*>0$, $\alpha_i^* = 1/M\nu$.
\end{lemma}
In words, \(\mathcal{U}_{\nu}^{\rm SMC}\) is guaranteed to contain at least \((1-\nu)\times 100\%\) of \(M\) source distributions, offering an interpretable data-centric mechanism to manage conservatism. Meanwhile, this flexibility in regulating conservatism confers robustness against some outlier-like distributions deviating apparently from the majority. By excluding a fraction of source distributions from the RKHS ball, \(\mathcal{U}_{\nu}^{\rm SMC}\) can mitigate the adverse effect caused by some extremal distributions. However, this flexibility is not possessed by the convex hull uncertainty set \(\mathcal{U}^{\rm co}\), which by design contains extremal distributions and is thus prone to conservatism. Moreover, we introduce a slack parameter \(\gamma>0\) into \(\mathcal{U}_{\nu}^{\rm SMC}\) as an additional means of regulating the robustness:
\begin{equation} \label{eq:SMC_ambiguity_set_gamma}
    \mathcal{U}_{\nu,\gamma}^{\text{SMC}} = \left\{\mathbb{P}\in\mathcal{P}(\Xi)\,\middle|\,\mathbb{E}_{\mathbb{P}}\left\{\phi(\xi)\right\}=\mu_{\mathbb{P}}, ~\mu_{\mathbb{P}}\in\mathcal{B}_{\nu,\gamma}\right\},
\end{equation}
where \(\mathcal{B}_{\nu,\gamma} = \{\mu\in\mathcal{H}_k\,|\,\|\mu-\mu_c^*\|_{\mathcal{H}_k}\leq \sqrt{(R^*)^2+\gamma}\}\) is an RKHS ball larger than $\mathcal{B}_{\nu}$. By tuning $\gamma>0$, one can control the size of the uncertainty set more flexibly, which is useful especially if $M$ is small. In particular, when \(M=1\), the original ball $\mathcal{B}_{\nu}$ in RKHS may collapse to a singleton with \((\mu_c^*,R^*)=(\mu_{\mathbb P_1},0)\). This results in a lack of robustness, which can be mitigated by further imposing a positive \(\gamma\).

In practice, the perfect knowledge of $M$ source distributions \(\{\mathbb{P}_i\}_{i=1}^M\) is rarely attainable. Therefore, we consider their finite-sample approximations in a data-driven context. Formally, we assume that all \(M\) datasets \(\mathcal{D}_1,\mathcal{D}_2,\cdots, \mathcal{D}_M\) are mutually independent, and for each source distribution \(\mathbb{P}_i\) a dataset \(\mathcal{D}_i=\{\boldsymbol{\xi}_j^{(i)}\}_{j=1}^{n_i}\) consisting of \(n_i\) i.i.d. samples is collected. The empirical counterparts of \(\mathbb{P}_i\) and \(\mu_{\mathbb{P}_i}\) can be calculated by replacing the integrals with finite sums over these samples \cite{muandet2013one}, i.e., \(\hat{\mathbb{P}}_i=\frac{1}{n_i}\sum_{j=1}^{n_i}\delta_{\boldsymbol{\xi}_j^{(i)}}\), and \(
\hat{\mu}_{\mathbb{P}_i}=\frac{1}{n_i}\sum_{j=1}^{n_i}\phi(\boldsymbol{\xi}_{j}^{(i)})\). As such, the empirical Gram matrix \(\hat{\mathbf{K}}\) can be readily computed without numerical integration. Its entries admit finite-sample approximations \(\hat{K}_{i_1i_2} = \langle \hat{\mu}_{\mathbb{P}_{i_1}},\hat{\mu}_{\mathbb{P}_{i_2}} \rangle_{\mathcal{H}_k}=  \frac{1}{n_{i_1} n_{i_2}}\sum_{j_1=1}^{n_{i_1}}\sum_{j_2=1}^{n_{i_2}}k(\boldsymbol{\xi}_{j_1}^{(i_1)},\boldsymbol{\xi}_{j_2}^{(i_2)})\). By substituting \(\hat{\mathbf{K}}\) for \(\mathbf{K}\), we arrive at a sample-based approximation of the dual problem~\eqref{eq:dual_SMC_quadratic}, which can also have non-unique optimizers. Invoking Proposition~\ref{prop:discussion_dual} enables to uniquely determine a dual optimizer \(\hat{\boldsymbol{\alpha}}^*\), along with the corresponding optimal center \(\hat{\mu}_c^*\) and radius \(\hat{R}^*\). Introducing a prescribed parameter \(\gamma>0\), we obtain an enlarged RKHS ball \(\hat{\mathcal B}_{\nu,\gamma}=\{\mu\in\mathcal H_k\,|\,\|\mu-\hat\mu_c^*\|_{\mathcal H_k}\leq \sqrt{(\hat{R}^*)^2+\gamma}\}\), which yields a data-driven uncertainty set:
\begin{equation} \label{eq:SMC_ambiguity_set_gamma_empirical}
    \widehat{\mathcal{U}}_{\nu,\gamma}^{\text{SMC}} = \left\{\mathbb{P}\in\mathcal{P}(\Xi)\,\middle|\,\mathbb{E}_{\mathbb{P}}\left\{\phi(\xi)\right\}=\mu_{\mathbb{P}}, ~\mu_{\mathbb{P}}\in\hat{\mathcal{B}}_{\nu,\gamma}\right\}.
\end{equation}
Note that \(\widehat{\mathcal{U}}_{\nu,\gamma}^{\rm SMC}\) admits a non-parametric representation through RKHS inner products, which is computationally tractable by the kernel trick:
\begin{equation}\label{eq:empirical_kernelized_SMC_uncertainty_set}
    \widehat{\mathcal{U}}_{\nu,\gamma}^{\rm SMC} = \left\{\mathbb{P}\in\mathcal{P}(\Xi)\;\middle|\; 
    \begin{aligned}
        &\langle\mu_{\mathbb{P}},\mu_{\mathbb{P}}\rangle_{\mathcal{H}_k}-2\sum_{i\in {\rm SM}}\sum_{j=1}^{n_i}\frac{\hat{\alpha}_{i}^*}{n_i}\langle\mu_{\mathbb{P}},\phi(\boldsymbol{\xi}_j^{(i)})\rangle_{\mathcal{H}_k}+\sum_{i\in \rm SM}\sum_{i'\in\rm SM}\alpha_i\alpha_i'K_{ii'}\leq (\hat{R}^*)^2+\gamma\\
        &\mu_{\mathbb{P}} = \mathbb{E}_{\mathbb{P}}[\phi(\boldsymbol{\xi})]
    \end{aligned}
    \right\}.
\end{equation}
Inserting the learned uncertainty set \(\widehat{\mathcal{U}}_{\nu,\gamma}^{\rm SMC}\) into problem~\eqref{eq:main_data_driven}, we arrive at the following robustified stochastic optimization problem, which yields a worst-case surrogate for Problem~\ref{Q} and also the mathematical backbone of our RooD-SO scheme:
\begin{equation} \label{eq:Kernel_SMM_DRO_Primary_finite}
    \min_{\mathbf{x} \in \mathcal{X}} \; \sup_{\mathbb{P} \in \widehat{\mathcal{U}}_{\nu,\gamma}^{\rm SMC}} 
   \; \mathbb{E}_{\boldsymbol{\xi} \sim \mathbb{P}} \left\{ f(\mathbf{x},\boldsymbol{\xi}) \right\}.
\end{equation}
This min-max problem can be interpreted as a two-player zero-sum game, where the decision-maker seeks an optimal decision with the lowest loss in the worst case, whereas a distributional adversary selects the least-favorable distribution in \(\widehat{\mathcal{U}}_{\nu,\gamma}^{\rm SMC}\) to describe the randomness of an unobserved distribution in a pessimistic sense. 

\begin{example}
    \begin{figure}[!ht]
        \centering
        \includegraphics[width=0.32\linewidth]{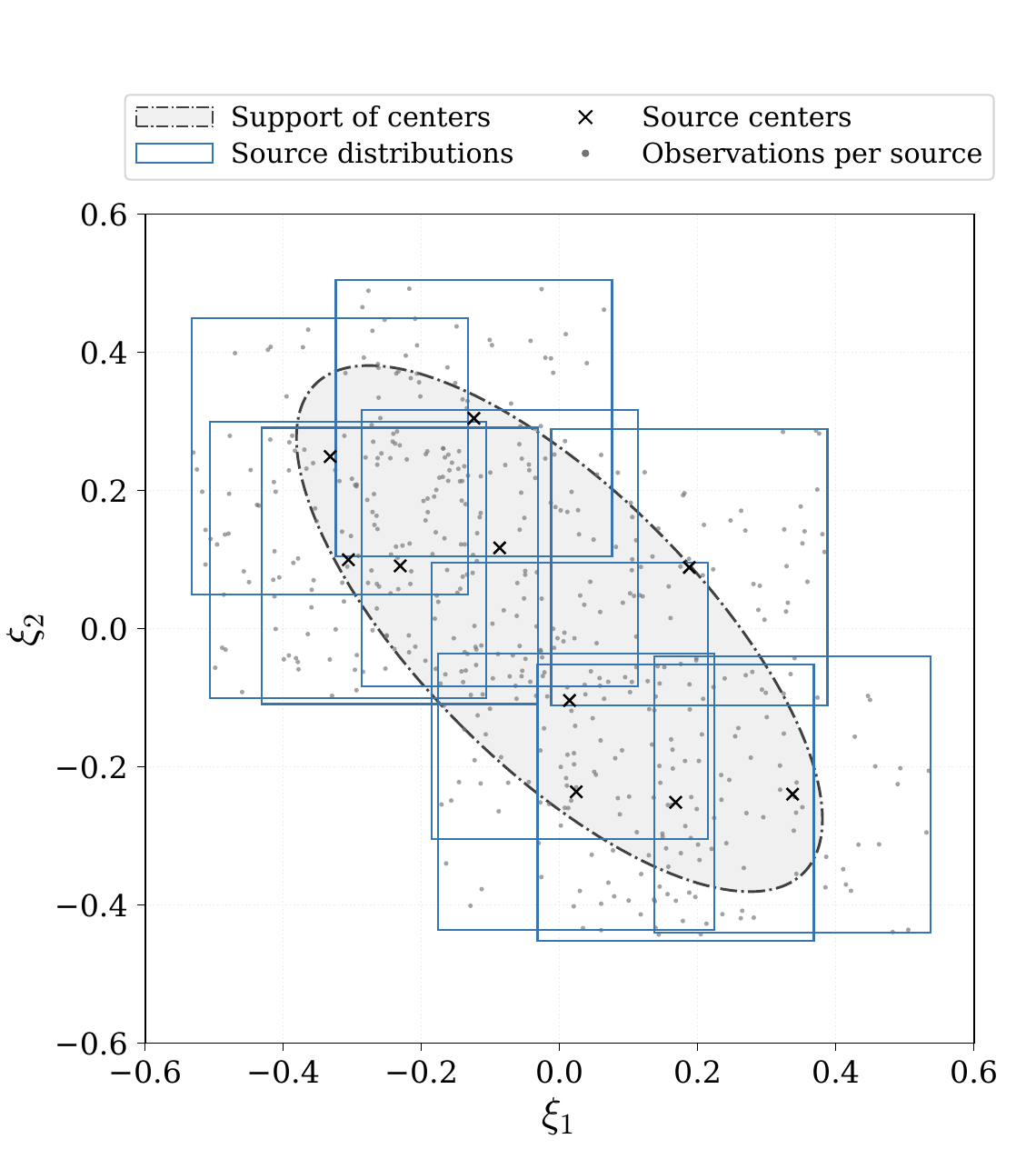}
        \caption{Distribution and data generation of a 2-D example}
        \label{fig:01_source_distributions}
    \end{figure}
    We use a 2-D meta-distribution example to intuitively illustrate the efficacy of using learning-based techniques to characterize distributional uncertainty. As shown in Figure~\ref{fig:01_source_distributions}, every plausible distribution is assumed to be a uniform distribution ${\rm Unif}(\mathbf c+[-0.2,0.2]^2)$ over a square area, whose center \(\mathbf c\in\mathbb R^2\) is randomly drawn from another top-layer uniform distribution ${\rm Unif}(\mathcal C)$ over a tilted ellipse \(\mathcal C=\left\{(z_1,z_2)\in\mathbb R^2:2(z_1-z_2)^2+12.5(z_1+z_2)^2\leq 1\right\}\). This naturally defines a meta-distribution and enables two-layer random sampling. 

    \begin{figure}[!ht]
    \centering
        \begin{subfigure}[b]{0.32\linewidth}
            \centering
            \includegraphics[width=\linewidth]{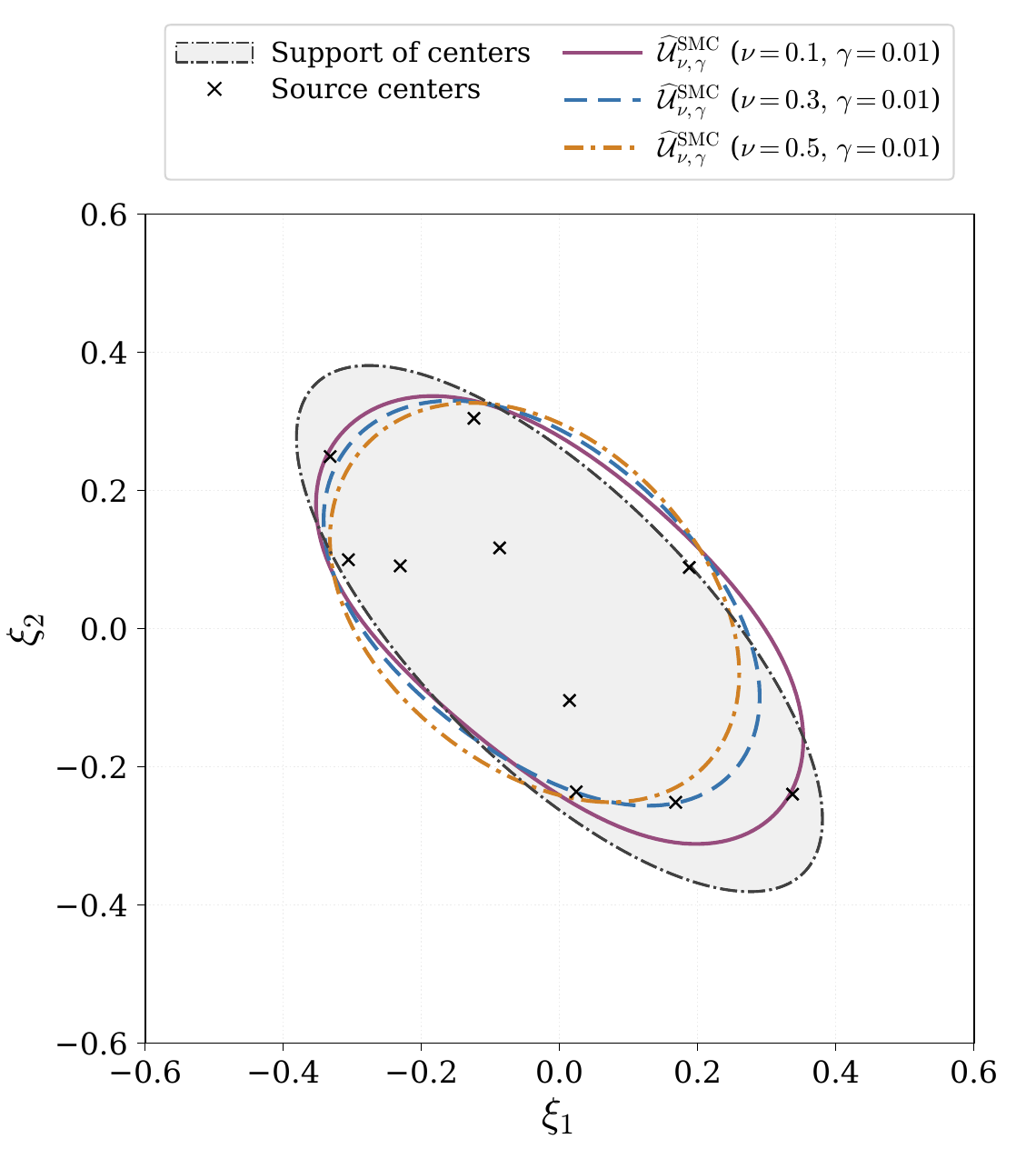}
            \caption{\(\widehat{\mathcal U}_{\nu,\gamma}^{\rm SMC}\) with varying \(\nu\)}
            \label{fig:06_roods_nu}
          \end{subfigure}
        \begin{subfigure}[b]{0.32\linewidth}
            \centering
            \includegraphics[width=\linewidth]{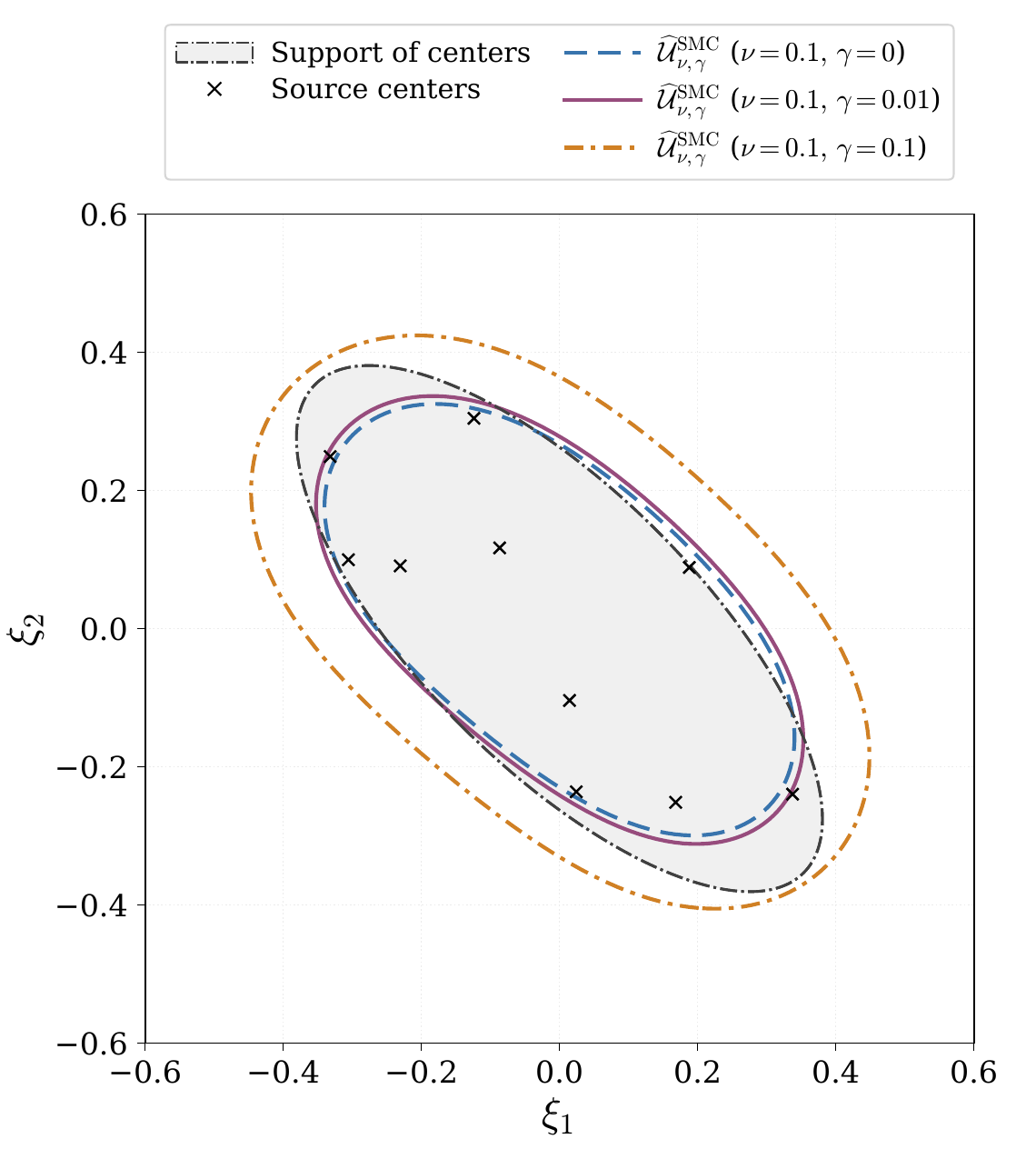}
            \caption{\(\widehat{\mathcal U}_{\nu,\gamma}^{\rm SMC}\) with varying \(\gamma\)}
            \label{fig:05_roods_gamma}
          \end{subfigure}
        \begin{subfigure}[b]{0.32\linewidth}
            \centering
            \includegraphics[width=\linewidth]{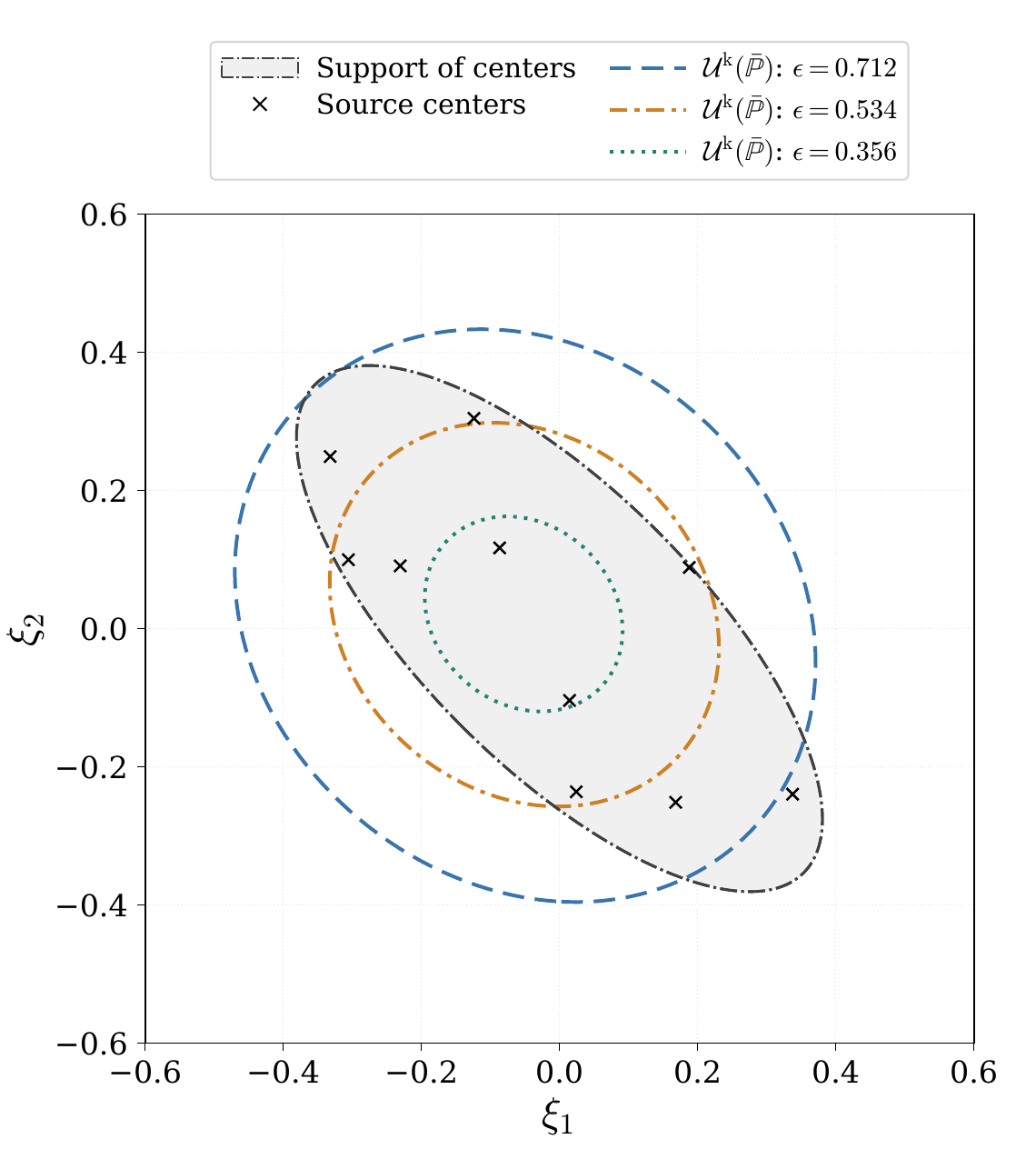}
            \caption{\(\mathcal U^{\rm k}\) centered at \(\bar{\mathbb P}\)}
            \label{fig:07_kernel_dro}
          \end{subfigure}

          \begin{subfigure}[b]{0.32\linewidth}
            \centering
            \includegraphics[width=\linewidth]{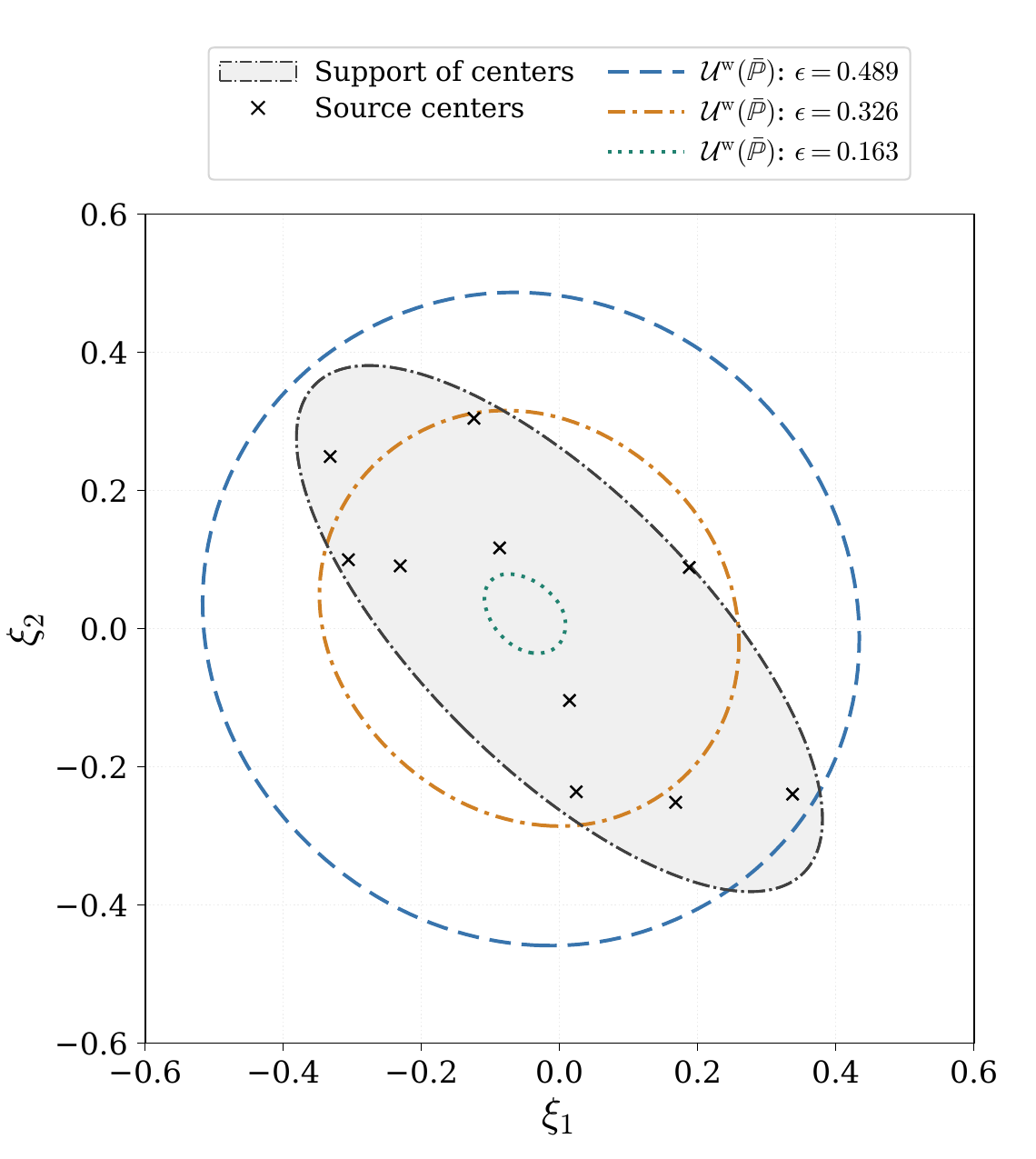}
            \caption{\(\mathcal U^{\rm w}\) centered at \(\bar{\mathbb P}\)}
            \label{fig:02_dro_w}
          \end{subfigure}
          \begin{subfigure}[b]{0.32\linewidth}
            \centering
            \includegraphics[width=\linewidth]{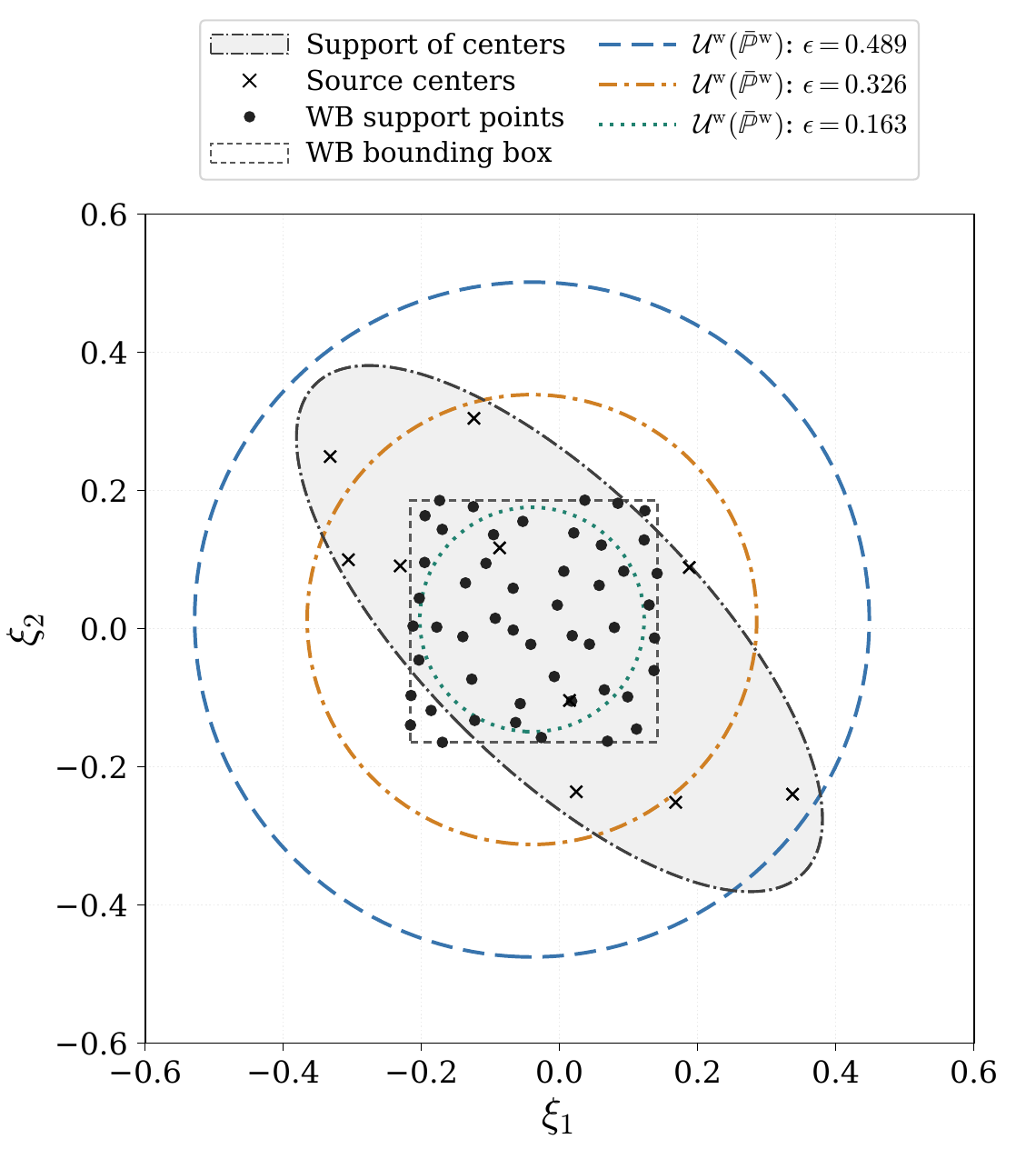}
            \caption{\(\mathcal U^{\rm w}\) centered at \(\bar{\mathbb P}^{\rm w}\)}
            \label{fig:03_dro_w_saa_wb}
          \end{subfigure}
            \begin{subfigure}[b]{0.32\linewidth}
            \centering
            \includegraphics[width=\linewidth]{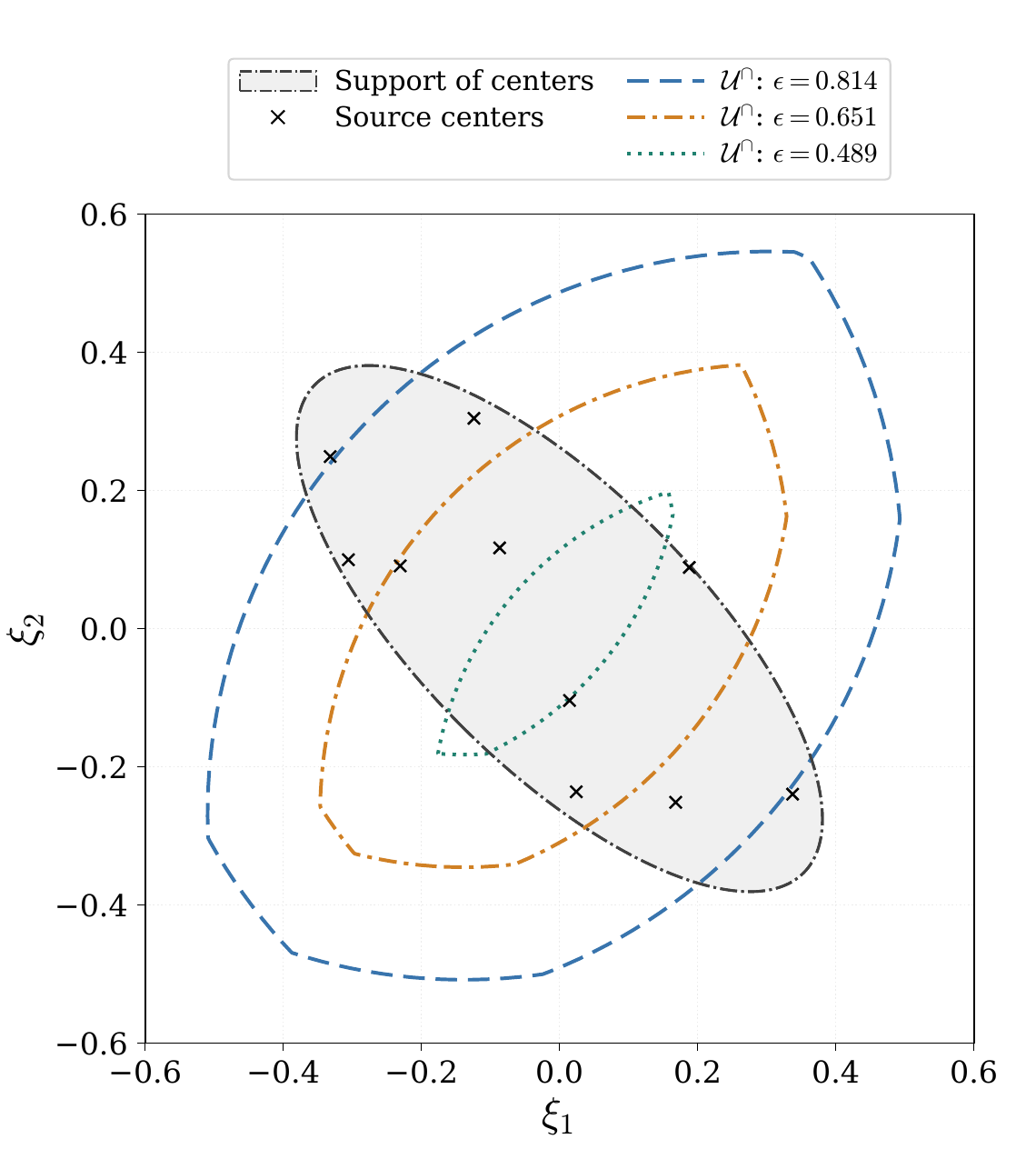}
            \caption{\(\mathcal U^{\cap}\)}
            \label{fig:04_dro_intx}
          \end{subfigure}
      \caption{Illustrative example: Visualization of several uncertainty-set construction methods.}
      \label{fig:example}
    \end{figure}
    
    We independently sample $M=10$ source distributions with random centers $\{\mathbf c_i\}_{i=1}^{10}$ and draw $n=50$ independent observations from each distribution. Then we construct a series of uncertainty sets via our proposed method alongside the previously reviewed alternatives. Given an uncertainty set, it is unlikely to visualize all distributions within. Because of this, we focus on a subset thereof, including all candidate uniform distributions ${\rm Unif}(\mathbf c+[-0.2,0.2]^2)$, and display the enclosures of their centers in 2-D space. This clearly illustrates whether the uncertainty set effectively captures the support of the distribution-generating meta-distribution.\endnote{The visualization of uncertainty sets in Figure \ref{fig:example} is achieved by scanning a fine grid over a broad domain for centers of the uniform-distribution family under consideration. Each grid point specifies a candidate distribution, which is included in the plotted enclosure if it is enclosed by the corresponding uncertainty set. For the Wasserstein-based sets, we compare the Wasserstein distance to the reference distribution with its prescribed radius. For RooD-SO, \eqref{eq:empirical_kernelized_SMC_uncertainty_set} can be used to check the inclusion. With the Gaussian kernel used, the required kernel inner products $\langle\mu_{\mathbb{P}},\mu_{\mathbb{P}}\rangle_{\mathcal{H}_k}$ and $\langle\mu_{\mathbb{P}},\phi(\boldsymbol{\xi}_j^{(i)})\rangle_{\mathcal{H}_k}$ admit analytical expressions for uniform distributions. The plotted enclosures are therefore grid approximations of the centers of this distribution family enclosed by uncertainty sets.} Figures~\ref{fig:06_roods_nu} and~\ref{fig:05_roods_gamma} show the results of \(\widehat{\mathcal U}_{\nu,\gamma}^{\rm SMC}\) in~\eqref{eq:empirical_kernelized_SMC_uncertainty_set} under varying values of $\nu$ and $\gamma$. At first glance, \(\widehat{\mathcal U}_{\nu,\gamma}^{\rm SMC}\) reasonably captures the geometry of the underlying distribution of source centers, offering empirical evidence of encouraging learning capability. Besides, these two subfigures showcase the effect of both \(\nu\) and \(\gamma\) on manipulating the conservatism. Clearly, a larger $\nu$ allows more sources to lie outside the fitted RKHS ball, while a larger $\gamma$ increases the RKHS ball radius and inflates the uncertainty set. In particular, the results highlight the role of $\gamma > 0$ in obtaining a safe coverage of \(\mathcal C\). Figure~\ref{fig:07_kernel_dro} shows the results of the RKHS balls centered at the embedding of the pooled empirical distribution \(\bar{\mathbb{P}}\) under different radii \(\epsilon\). It can be seen that such a trivial construction cannot compactly capture the underlying geometry of the distribution of \(\mathbf c\) under different values of \(\epsilon\). For the sake of safe coverage, a large radius is required to inflate the set but causes an excess of implausible candidate uniform distributions to be included. Recall that both \(\widehat{\mathcal U}_{\nu,\gamma}^{\rm SMC}\) and \(\mathcal U^{\rm k}\) take the form of RKHS ball. Thus, Figures~\ref{fig:06_roods_nu}--(c) highlight the efficacy of data-driven kernel learning in capturing the support of the meta-distribution.
    
    In Figures~\ref{fig:02_dro_w}--(f), we show the results of Wasserstein balls \(\mathcal U^{\rm w}\) centered at \(\bar{\mathbb{P}}\), Wasserstein balls around the WB of \(\{\hat{\mathbb P}_i\}\), and the intersections of Wasserstein balls \(\mathcal U^{\cap}\), respectively, where these sets exhibit varying degrees of mismatch with the geometry of \(\mathcal C\). For two Wasserstein sets in Figures~\ref{fig:02_dro_w} and~\ref{fig:03_dro_w_saa_wb}, perturbations via mass transportation are essentially isotropic in the distribution space, thereby failing to capture the anisotropic geometry of $\mathcal{C}$. As shown in Figure~\ref{fig:04_dro_intx}, the intersection-based uncertainty set \(\mathcal U^{\cap}\) manifests a significant misalignment with \(\mathcal C\) and captures more uniform distributions than necessary for large ball radii $\epsilon$. We remark that in this example the convex-hull uncertainty set \(\mathcal U^{\rm co}\) in~\eqref{eq:Group_DRO_ambiguity_set} cannot be visualized. This is because every distribution \(\mathbb P\in{\mathcal U}^{\rm co}\) is finitely supported on the observed samples, and thus \(\mathcal U^{\rm co}\) does not include any candidate uniform distribution \({\rm Unif}(\mathbf c+[-0.2,0.2]^2)\).

\end{example}

\subsection{Finite-sample error analysis}
\label{sec:data_driven_setup}
The use of SMC has been extensively exploited from a machine learning perspective, mostly for unsupervised group anomaly detection~\cite{muandet2013one, guevara2015support}, while the errors induced by finite-sample approximation have not been formally analyzed. In this subsection, we fill this knowledge gap with a systematic theoretical analysis of finite-sample approximation errors. In particular, finite-sample approximation error enters \(\hat{\mathbf{K}}\), and eventually contaminates \(\hat{\boldsymbol{\alpha}}^*\) as well as the learned uncertainty set. 
\begin{proposition}[Error bound for empirical Gram matrix] \label{prop:operator_norm_bound_Khat_K}
For any \(\delta\in(0,1)\), with probability at least \(1-\delta\) over \(\mathcal{D}_i\sim\mathbb{P}_i^{n_i},~i\in[M]\), the following holds:
\begin{equation}\label{eq:main_Khat_K_norm_bound}
\|\hat{\mathbf{K}}-\mathbf{K}\|_2\leq\sqrt{\frac{8MC^2(1+2\log(M/\delta))}{\bar{n}}}+\frac{2C(1+2\log(M/\delta))}{\bar{n}},
\end{equation}
where \(\bar{n}:=(\sum_{i=1}^M1/n_i)^{-1}\).
\end{proposition}
Proposition~\ref{prop:operator_norm_bound_Khat_K} derives a concentration inequality for \(\hat{\mathbf{K}}\) under the randomness induced by drawing \(n_i\) samples from each \(\mathbb{P}_i\), implying that the spectrum of \(\hat{\mathbf{K}}\) is close to that of \(\mathbf{K}\) with high probability. More precisely, since \(|\boldsymbol{\alpha}^\top(\hat{\mathbf{K}}-\mathbf{K})\boldsymbol{\alpha}|\leq \|\hat{\mathbf{K}}-\mathbf{K}\|_2\|\boldsymbol{\alpha}\|_2^2\), the quadratic term in the objective of~\eqref{eq:dual_SMC_quadratic} differs by at most \(\mathcal{O}(\|\hat{\mathbf{K}}-\mathbf{K}\|_2)\). Now we are in a position to establish the consistency of the (possibly non-unique) dual optimizer \(\hat{\boldsymbol{\alpha}}^*\).
\begin{proposition}[Error bound for empirical optimizer] \label{prop:consistency_optimizer}
For any \(\delta\in(0,1)\), with probability at least \(1-\delta\) over \(\mathcal{D}_i\sim\mathbb{P}_i^{n_i}\) for \(i\in[M]\), the following holds:
\begin{equation}\label{eq:alpha_optimizer_consistency}
\operatorname{dist}(\hat{\boldsymbol{\alpha}}^*,\mathcal A^*)\leq \frac{4C}{\rho}\left[\left(C_\alpha\sqrt{8M}+\sqrt{2}\right)\sqrt{\frac{1+2\log(2M/\delta)}{\bar{n}}}+ \frac{2C_\alpha(1+2\log(2M/\delta))}{\bar{n}}\right],
\end{equation}
where \(C_\alpha=\sqrt{1\wedge(M\nu)^{-1}}\) and \(\rho>0\) is a constant depending on \(\mathbf K\) and \(\nu\). 
\end{proposition}

An important implication is that when \(\bar{n}\) grows to infinity, namely \(n_i\rightarrow\infty,~\forall i\), we have \(\operatorname{dist}(\hat{\boldsymbol\alpha}^*,\mathcal A^*)=\mathcal O_{\mathbb P}\!\left(\bar n^{-1/2}\right)\). In other words, any dual optimizer under data-driven approximation converges in probability to the population optimal
solution set \(\mathcal A^*\) that need not be a singleton. We also remark that \(\rho\) is related to the quadratic growth of the dual objective function~\eqref{eq:objective_function} away from its optimal solution set~\cite{drusvyatskiy2018error}, while a smaller \(\rho\) yields a larger perturbation bound.

We close this subsection with an analysis of how finite-sample estimation errors propagate to the learned uncertainty set \(\widehat{\mathcal{U}}_{\nu,\gamma}^{\rm SMC}\). It is crucial to quantify these discrepancies because geometric changes in the ball can noticeably affect the expressive power of the uncertainty set. In the following, we establish the statistical properties of these discrepancies by deriving high-probability bounds for both \(\|\hat{\mu}_c^*-\mu_c^*\|_{\mathcal{H}_k}\) and \(|(\hat{R}^*)^2-(R^*)^2|\). We then show that these bounds are nontrivial given \(M\) distributions.
\begin{theorem}[Convergence rate of the estimated center and radius]
\label{thm:finite_sample_guar_mu_c_R_square}
For a fixed \(M\), as $\bar n\to\infty$,
\[
\|\hat{\mu}_c^*-\mu_c^*\|_{\mathcal H_k} = \mathcal{O}_{\mathbb P}\!\left(\bar n^{-1/2}\right), \qquad \left|(\hat R^*)^2-(R^*)^2\right| = \mathcal{O}_{\mathbb P}\!\left(\bar n^{-1/2}\right).
\]
\end{theorem}
Theorem~\ref{thm:finite_sample_guar_mu_c_R_square} indicates that for a fixed number $M$ of source distributions, the estimated center $\hat{\mu}_c^*$ and squared radius $(\hat R^*)^2$ both attain the standard root-$\bar n$ convergence rate. Moreover, the rates of the empirical mean embeddings, the estimated optimizer, and the resulting center and radius have the same order, implying that errors may propagate but do not compound in a way that worsens the rate.


\subsection{Out-of-distribution generalization guarantees}
\label{sec:theoretical_analysis}

Thus far, we have studied the randomness in sampling $n_i$ data points from each source distribution $\mathbb{P}_i$ and its impact on the learned uncertainty set. Beyond that, there is a higher layer of randomness in sampling $M$ distributions $\mathbb{P}_1, \mathbb{P}_2, \cdots, \mathbb{P}_M$ from the meta-distribution $\Psi$. To formally account for such two-layer randomness, we draw inspiration from statistical learning theory to establish out-of-distribution generalization guarantees for the uncertainty set as well as the resulting solution. To ease analysis, we first restrict our attention to the top-layer randomness in $\Psi$. Assume that \(M\) source distributions $\mathbb{P}_1, \mathbb{P}_2, \cdots, \mathbb{P}_M$, as independent samples from $\Psi$, are perfectly known such that an ``ideal" uncertainty set \(\mathcal{U}_{\nu,\gamma}^{\rm SMC}\) can be constructed without finite-sample approximations. Borrowing ideas from kernel machines~\cite{scholkopf2002learning}, we present an out-of-distribution probabilistic guarantee that a new distribution \(\mathbb{P}'\sim\Psi\) falls within \(\mathcal{U}_{\nu,\gamma}^{\rm SMC}\). 
\begin{theorem}[Out-of-distribution generalization guarantee for \(\mathcal{U}_{\nu,\gamma}^{\rm SMC}\)]\label{thm:SMC_generalization_error_bound}
    Assume \(\mathbb{P}_1,\cdots,\mathbb{P}_M\) are drawn i.i.d. from a meta-distribution \(\Psi\). Then, for any \(\delta\in(0,1)\) and fixed \(\gamma>0\), the following inequality holds with probability at least \(1-\delta\) over the draws of \(\{\mathbb{P}_i\}_{i=1}^M\):
    \begin{equation}\label{ineq:generalization_population}
        {\rm Pr}_{\mathbb{P}'\sim\Psi}\left\{\mathbb{P'}\notin \mathcal{U}_{\nu,\gamma}^{\rm SMC}\right\} \leq \nu+\frac{6C}{\gamma\sqrt{M}}+\sqrt{\frac{\log(1/\delta)}{2M}}.
    \end{equation}
\end{theorem}
The first term on the right-hand side of \eqref{ineq:generalization_population} upper-bounds the fraction of \(M\) source distributions that are excluded from \(\mathcal{U}_{\nu,\gamma=0}^{\rm SMC}\) in virtue of Lemma~\ref{lemma:nu}. Meanwhile, the second term captures the effect of \(\gamma\), revealing the relation between the confidence level of the theorem and the size of the expanded uncertainty set. Meanwhile, the last two terms hinge on \(M\) and vanish at rate \(\mathcal{O}(M^{-1/2})\), stemming from the uncertainty in observing finite source distributions from $\Psi$. Next, we further account for the lower-layer randomness arising from finite-sample observations by considering the learned uncertainty set \(\widehat{\mathcal{U}}_{\nu,\gamma}^{\rm {SMC}}\)~\eqref{eq:empirical_kernelized_SMC_uncertainty_set}.

\begin{proposition}[Out-of-distribution generalization guarantee for \(\widehat{\mathcal{U}}_{\nu,\gamma}^{\rm SMC}\)]\label{prop:empirical_generalization_error_bound}
    Assume $M$ source distribution \(\mathbb{P}_1,\cdots,\mathbb{P}_M\) are drawn i.i.d. from a meta-distribution \(\Psi\), and \(n_i\) observations are sampled i.i.d. from each \(\mathbb{P}_i\). For any \(\delta\in(0,1)\) and fixed \(\gamma>0\), the following inequality holds with probability \(1-\delta\) over \(\{\mathbb P_i\}_{i=1}^M\sim\Psi^M\) and \(\mathcal{D}_i\sim\mathbb{P}_i^{n_i}\) for \(i\in[M]\):
    \begin{equation}\label{ineq:finite_sample_generalization_error_bound}
        {\rm Pr}_{\mathbb{P}'\sim\Psi}\left\{\mathbb{P}'\notin \widehat{\mathcal{U}}_{\nu,\gamma}^{\rm {SMC}}\right\}\leq\nu +\frac{6C}{\gamma\sqrt M}+\sqrt{\frac{\log(1/\delta)}{2M}}+\frac{4 C}{\gamma M}\sum_{i=1}^M\frac{1}{\sqrt{n_i}}.
    \end{equation}
\end{proposition}

Relative to \eqref{ineq:generalization_population}, the out-of-distribution generalization error bound~\eqref{ineq:finite_sample_generalization_error_bound} for \(\widehat{\mathcal U}_{\nu,\gamma}^{\rm SMC}\) incurs an additional term related to the price of learning from finite samples. In a balanced sampling regime \(n_1=\cdots=n_M=n\), the last term in~\eqref{ineq:finite_sample_generalization_error_bound} reduces to \(4C/(\gamma\sqrt n)\) and thus decays at a rate of \(\mathcal O(n^{-1/2})\).  

Theorem~\ref{thm:SMC_generalization_error_bound} and Proposition~\ref{prop:empirical_generalization_error_bound} sheds light on the roles of $\nu$ and $\gamma$ in securing the out-of-distribution performance. It is worth noting that the bounds \eqref{ineq:generalization_population} and \eqref{ineq:finite_sample_generalization_error_bound} can be conservative, especially when $M$ is small, and do not by themselves provide a practical rule for calibrating $\nu$ and $\gamma$. Rather, they clarify the trade-off between generalization and robustness, thereby establishing theoretical support for using $\nu$ and $\gamma$ to adjust the conservatism of the uncertainty set. Next, we proceed to the asymptotic behavior of \(\widehat{\mathcal U}_{\nu,\gamma}^{\rm SMC}\).

\begin{proposition}[Asymptotic behavior of \(\widehat{\mathcal{U}}_{\nu,\gamma}^{\rm SMC}\)] \label{prop:asymptotic_behavior_U_empirical}
Assume $M$ source distribution \(\mathbb{P}_1,\cdots,\mathbb{P}_M\) are drawn i.i.d. from a meta-distribution \(\Psi\), and \(n_i\) observations are sampled i.i.d. from each \(\mathbb{P}_i\). Let \(R_M:=\Psi^M\otimes\mathbb{P}_1^{n_{1,M}}\otimes\cdots\otimes\mathbb{P}_{M}^{n_{M,M}}\) be the resulting joint law, and \(\gamma_M:=2\left(\frac{1}{M}\sum_{i=1}^M\frac{1}{\sqrt{n_{i,M}}}+\frac{1}{\sqrt{M}}\right)^{1/2}\). If \(\gamma_M\to0\) as \(M\to \infty\), then for every \(\epsilon>0\),
\begin{equation}\label{eq:asymptotic_SMC_empirical}
    \lim_{M\rightarrow \infty}{\rm Pr}_{R_M}\left\{{\rm Pr}_{\mathbb{P}'\sim\Psi}\left\{\mathbb{P}'\in \widehat{\mathcal{U}}_{\nu,\gamma_M}^{\rm {SMC}}\right\}\geq1-\nu-\epsilon\right\}=1.
\end{equation}
\end{proposition}

Proposition~\ref{prop:asymptotic_behavior_U_empirical} demonstrates that as \(M\to\infty\), \(\widehat{\mathcal{U}}_{\nu,\gamma_M}^{\rm SMC}\) asymptotically achieves a coverage probability of at least \(1-\nu\), in probability over the sampled data. Notably, this result does not require every \(n_{i,M}\) to diverge, and it suffices that the average \(\sum_{i=1}^{M}n_{i,M}^{-1/2}/M\) converges to zero as \(M\) goes to infinity. Leveraging these coverage guarantees, we establish out-of-distribution performance guarantees for the robust worst-case objective under the two-layer randomness. 

\begin{proposition}[Nonasymptotic out-of-distribution probabilistic bound for solutions] 
    \label{prop:nonasymptotic_out_of_distribution_guarantee_DRO_finite_sample}
    For any \(\delta\in(0,1)\) and fixed \(\gamma>0\), the following inequality holds with probability \(1-\delta\) over \(\{\mathbb P_i\}_{i=1}^M\sim\Psi^M\) and \(\mathcal{D}_i\sim\mathbb{P}_i^{n_i}\) for \(i\in[M]\):
    \begin{equation}
        {\rm Pr}_{\mathbb{P}'\sim\Psi}\left\{\forall\,\mathbf{x}\in\mathcal{X}:\mathbb{E}_{\mathbb{P}'}\left\{f(\mathbf{x},\boldsymbol{\xi})\right\}\leq \sup_{\mathbb{P}\in\widehat{\mathcal{U}}_{\nu,\gamma}^{\rm SMC}}\mathbb{E}_{\mathbb{P}}\left\{f(\mathbf{x},\boldsymbol{\xi})\right\}\right\}\geq 1-\nu -\frac{6C}{\gamma\sqrt M}-\sqrt{\frac{\log(1/\delta)}{2M}}-\frac{4 C}{\gamma M}\sum_{i=1}^M\frac{1}{\sqrt{n_i}}.
    \end{equation}
\end{proposition}
Proposition~\ref{prop:nonasymptotic_out_of_distribution_guarantee_DRO_finite_sample} shows that, with high probability over the \(M\) source distributions and their samples, the worst-case expected loss provides a uniform probabilistic upper bound on the out-of-distribution expected loss. This uniformity over \(\mathbf x\in\mathcal X\) is essential because the resulting optimizer is inherently data-dependent. Consequently, the optimal value provides a probabilistic upper bound on the out-of-distribution expected loss of the corresponding optimizer. The following proposition establishes the asymptotic validity of this performance certificate under certain sampling conditions. 
\begin{proposition}[Asymptotic out-of-distribution probabilistic bound for solutions] \label{prop:out_of_sample_guarantee_DRO_finite_sample}
    Assume $M$ source distribution \(\mathbb{P}_1,\cdots,\mathbb{P}_M\) are drawn i.i.d. from a meta-distribution \(\Psi\), and \(n_i\) observations are sampled i.i.d. from each \(\mathbb{P}_i\). Let \(R_M:=\Psi^M\otimes\mathbb{P}_1^{n_{1,M}}\otimes\cdots\otimes\mathbb{P}_{M}^{n_{M,M}}\) be the resulting joint law, and \(\gamma_M:=2\left(\frac{1}{M}\sum_{i=1}^M\frac{1}{\sqrt{n_{i,M}}}+\frac{1}{\sqrt{M}}\right)^{1/2}\). If \(\gamma_M\to0\) as \(M\to \infty\), then for every \(\epsilon>0\),
    \begin{equation}
        \lim_{M\rightarrow \infty}{\rm Pr}_{R_M}\left\{{\rm Pr}_{\mathbb{P}'\sim\Psi}\left\{\forall\,\mathbf{x}\in\mathcal{X}:\mathbb{E}_{\mathbb{P}'}\left\{f(\mathbf{x},\boldsymbol{\xi})\right\}\leq \sup_{\mathbb{P}\in\widehat{\mathcal{U}}_{\nu,\gamma_M}^{\rm SMC}}\mathbb{E}_{\mathbb{P}}\left\{f(\mathbf{x},\boldsymbol{\xi})\right\}\right\}\geq 1-\nu-\epsilon\right\}=1.
    \end{equation}
\end{proposition}
Intuitively, under the stated sampling condition, the robust optimal value provides an asymptotically valid \((1-\nu)\)-level out-of-distribution performance certificate for the optimizer.


\section{Solution strategy}
\label{sec:solution}
The robustified stochastic optimization problem~\eqref{eq:Kernel_SMM_DRO_Primary_finite} plays a central role in the proposed RooD-SO scheme. In this section, we delve into the computational aspect and design an efficient solution strategy.

\subsection{Dual reformulation and safe approximation}
We first reformulate the primal min-max problem~\eqref{eq:Kernel_SMM_DRO_Primary_finite} using duality theory. To our knowledge, \cite{zhu2021kernel} developed a unified kernel DRO framework that reinterprets several classical DRO formulations through constructing ambiguity sets in an RKHS. Their idea is useful for deriving a dual reformulation of problem~\eqref{eq:Kernel_SMM_DRO_Primary_finite} because \(\hat{\mathcal{B}}_{\nu,\gamma}\) is expressible in the same form as the RKHS-norm-ball ambiguity set in \cite{zhu2021kernel}. The following lemma presents the generalized duality result established in~\cite{zhu2021kernel}, which lays the foundation for our subsequent development.

\begin{lemma}[Generalized duality theorem \cite{zhu2021kernel}]
    Assume that \(f\) is a proper, upper semi-continuous function. Let \(\mathcal{U} := \{\mathbb{P}\in\mathcal{P}(\Xi):\mathbb{E}_{\mathbb{P}}\left\{\phi(\boldsymbol{\xi})\right\}=\mu_{\mathbb{P}}, \mu_{\mathbb{P}}\in\mathcal{M}\}\), where \(\mathcal{M} \subset\mathcal{H}_k\) is a closed, convex set such that \(\text{ri}\,(\mathcal{U})\neq \varnothing\). Then,
    \begin{equation}
        \begin{aligned}
            \sup_{\mathbb{P},\mu_{\mathbb{P}}}&~\mathbb{E}_{\mathbb{P}}\left\{f(\boldsymbol{\xi})\right\}\\
            {\rm s.t.}&~\mathbb{E}_{\mathbb{P}}\left\{\phi(\boldsymbol{\xi})\right\}=\mu_{\mathbb{P}},~\mu_{\mathbb{P}}\in\mathcal{M},~\mathbb{P}\in\mathcal{P}(\Xi)
        \end{aligned}
    \end{equation}
    is equivalent to 
    \begin{equation}
        \begin{aligned}            \min_{g\in\mathcal{H}_k,\beta\in\mathbb{R}}~&h_{\mathcal{M}}^*(g)+\beta\\
            {\rm s.t.}~\quad& f(\boldsymbol{\xi})\leq g(\boldsymbol{\xi})+\beta, \quad \forall~\boldsymbol{\xi}\in \Xi.\\
        \end{aligned}
        \label{eq:generalized_duality}
    \end{equation}
    Moreover, strong duality holds.
    \label{lemma:generalized_duality}
\end{lemma}

Note that \(g\in\mathcal{H}_k\) is an introduced dual variable, whose pointwise value is characterized by the reproducing property, i.e.,  \(g(\boldsymbol{\xi})=\langle g,\phi(\boldsymbol{\xi})\rangle_{\mathcal{H}_k}\). In spirit, solving for the worst-case distribution is tantamount to finding a surrogate function \(g+\beta\) that majorizes \(f\) and then minimizing the associated loss over \(g\) and \(\beta\)~\cite{zhu2021kernel}. Invoking Lemma~\ref{lemma:generalized_duality}, we arrive at a dual reformulation of problem~\eqref{eq:Kernel_SMM_DRO_Primary_finite}.

\begin{proposition} [Dual reformulation] \label{prop:dual_reformulation_SMC}
    Define \(\hat R^*_\gamma:=\sqrt{(\hat R^*)^2+\gamma}\). Problem \eqref{eq:Kernel_SMM_DRO_Primary_finite} can be reformulated as the following equivalent problem: 
    \begin{subequations}\label{eq:SMC_infinite_dual_formulation}
        \begin{align}
            z_{\rm D}:=\min_{\mathbf x\in\mathcal X, g\in \mathcal{H}_k,\beta \in\mathbb{R}}~&\sum_{i\in {\rm SM}}\sum_{j=1}^{n_i}\frac{\hat{\alpha}_i^*}{n_i}g(\boldsymbol{\xi}_j^{(i)}) + \hat{R}^*_\gamma\|g\|_{\mathcal{H}_k} + \beta\\
            {\rm s.t.}~~~\quad&f(\mathbf x, \boldsymbol{\xi})\leq g(\boldsymbol{\xi})+\beta,\quad\forall\,\boldsymbol{\xi}\in \Xi.\label{ineq:infinite_constraint}
        \end{align}
\end{subequations}
\end{proposition}

Solving problem~\eqref{eq:SMC_infinite_dual_formulation} is computationally formidable due to its semi-infinite nature, including the functional form of the dual variable $g$ and infinitely many constraints \eqref{ineq:infinite_constraint}. Indeed, techniques in modern robust optimization~\cite{ben2015deriving} are not applicable to the semi-infinite constraints~\eqref{ineq:infinite_constraint}, which, for general kernels, do not admit tractable finite-dimensional reformulations. To our knowledge, the scenario programming approximation is the only existent computational expedient~\cite{calafiore2006scenario,YANG2022107953, zhu2021kernel}. It builds upon a finite set of samples \(\Xi_N:=\{\boldsymbol\zeta_1,\ldots,\boldsymbol\zeta_N\}\) drawn from \(\Xi\) and replaces the semi-infinite constraint with finite sampled ones:
\begin{equation}\label{eq:sampled_relaxed_problem}
\begin{aligned}
    z_{\rm LB}:=\min_{\mathbf x\in\mathcal X, g\in\mathcal H_k,\,\beta\in\mathbb R}\quad&\sum_{i\in{\rm SM}}\sum_{j=1}^{n_i}\frac{\hat{\alpha}_i^*}{n_i}g(\boldsymbol{\xi}_j^{(i)})+\hat R^*_\gamma\|g\|_{\mathcal H_k} +\beta\\
    \text{s.t.}~~~~~~\quad&f(\mathbf x, \boldsymbol{\zeta}_l)\leq g(\boldsymbol{\zeta}_l)+\beta,\qquad \forall l\in[N].
\end{aligned}
\end{equation}
Apparently, problem~\eqref{eq:sampled_relaxed_problem} involving $N$ sampled constraints represents a relaxation of the dual problem~\eqref{eq:SMC_infinite_dual_formulation} as well as the primal min-max problem \eqref{eq:Kernel_SMM_DRO_Primary_finite}, which confirms a lower bound on the optimal value of \eqref{eq:SMC_infinite_dual_formulation}:
\begin{equation} \label{ineq:z_LB_leq_z_D}
    z_{\rm LB}\leq z_{\rm D}.
\end{equation}
As a consequence, the optimal solution to \eqref{eq:sampled_relaxed_problem} may be infeasible for the infinitely constrained dual problem and thus no longer inherits the out-of-distribution guarantees established in Section~\ref{sec:theoretical_analysis}.

Next, we design a novel approximation scheme for the dual problem~\eqref{eq:SMC_infinite_dual_formulation} by appropriately tightening the sampled constraints in the relaxed problem~\eqref{eq:sampled_relaxed_problem}. This yields a \textit{safe approximation} of both the dual problem \eqref{eq:SMC_infinite_dual_formulation} and the primal problem \eqref{eq:Kernel_SMM_DRO_Primary_finite} while preserving tractability.
\begin{theorem} \label{thm:tightening_and_prior_bound}
    Consider the following optimization problem that tightens constraints in~\eqref{eq:sampled_relaxed_problem}:    \begin{subequations}\label{eq:tightened_problem}
        \begin{align}
            z_{\rm UB}:=\min_{\mathbf x\in\mathcal X, g\in \mathcal{H}_k,\beta \in\mathbb{R}}~&\sum_{i\in{\rm SM}}\sum_{j=1}^{n_i}\frac{\hat{\alpha}_i^*}{n_i}g(\boldsymbol{\xi}_j^{(i)}) + \hat{R}^*_\gamma\|g\|_{\mathcal{H}_k} + \beta\\
                {\rm s.t.}~~~\quad &f(\mathbf x, \boldsymbol{\zeta}_l)+\omega^f_{r_N}(\mathbf x)+\omega^k_{r_N}\|g\|_{\mathcal H_k}\leq g(\boldsymbol{\zeta}_l) +\beta, \quad \forall~l\in[N],\label{ineq:certificate_margin_finite_dimensional_constraint}
         \end{align}
    \end{subequations} 
    where \(\omega^{k}_r:=\sup_{\|\boldsymbol\xi-\boldsymbol\zeta\|\leq r} \sqrt{k(\boldsymbol\xi,\boldsymbol\xi)+k(\boldsymbol\zeta,\boldsymbol\zeta)-2k(\boldsymbol\xi,\boldsymbol\zeta)}\), \(\omega^{f}_r(\mathbf x)\) is a convex majorant of function \(\sup_{\|\boldsymbol\xi-\boldsymbol\zeta\|\leq r} |f(\mathbf x,\boldsymbol\xi)-f(\mathbf x,\boldsymbol\zeta)|\) satisfying \(\lim_{r\downarrow0}\omega^f_r(\mathbf x)=0,~\forall \mathbf x\in\mathcal X\), and \(r_N:=\sup_{\boldsymbol{\xi}\in\Xi}\min_{i\in[N]}\|\boldsymbol{\xi}-\boldsymbol{\zeta}_i\|\). Then the following statements hold:
    \begin{enumerate}[label=(\alph*), leftmargin=*]
        \item \textbf{Feasibility:} Every feasible solution to~\eqref{eq:tightened_problem} is also feasible for~\eqref{eq:SMC_infinite_dual_formulation}, and thus 
        \begin{equation} \label{ineq:zD_leq_zUB}
            z_{\rm D}\leq z_{\rm UB}.
        \end{equation}
        \item \textbf{Bounded optimality gap:} Let \(\mathbf x^*_{\rm D}\) be an optimal solution to~\eqref{eq:SMC_infinite_dual_formulation}. The optimal value \(z_{\rm UB}\) is bounded above by
        \begin{equation}
            z_{\rm UB}\leq z_{\rm D} + \omega^f_{r_N}(\mathbf x_{\rm D}^*)+\bar B_\gamma \omega^k_{r_N},
        \end{equation}
        where $\bar B_\gamma:=(\sup_{\boldsymbol{\xi}\in\Xi}f(\mathbf x^*_{\rm D}, \boldsymbol{\xi})-\inf_{\boldsymbol{\xi}\in\Xi}f(\mathbf x^*_{\rm D}, \boldsymbol{\xi})) / \hat R^*_\gamma$. \label{thm2:bounded_optimality_gap}
        \item \textbf{Consistency:} If \(\lim_{N\to\infty}r_N=0\), then:
        \begin{equation}
            \lim_{N\rightarrow\infty}z_{\rm UB}=z_{\rm D}.
        \end{equation}
    \end{enumerate}
\end{theorem}
Interestingly, this theorem shows that we recapture the feasibility of the infinitely constrained dual problem~\eqref{eq:SMC_infinite_dual_formulation} by merely adding a safety margin $\omega^f_{r_N}(\mathbf x)+\omega^k_{r_N}\|g\|_{\mathcal H_k}$ to $N$ sampled constraints in the relaxed problem \eqref{eq:sampled_relaxed_problem}. The induced optimality gap $\omega^f_{r_N}(\mathbf x_{\rm D}^*)+\bar B_\gamma \omega^k_{r_N}$ is contingent on \(r_N\), which is precisely the \textit{covering radius} of $\Xi_N$ over \(\Xi\) and thus characterizes the density of $\Xi_N$. Indeed, given $\Xi$ and \(\Xi_N\),\endnote{If \(\Xi\) is unknown, we can estimate it from data using off-the-shelf support estimation algorithms.} the value of \(r_N\) can be effectively computed using the conventional Voronoi diagram method~\cite{aurenhammer1991voronoi}. Furthermore, the optimality gap becomes narrower as the coverage density of $\Xi_N$ gets higher, so one can close the gap by choosing $\Xi_N$ dense enough such that $r_N \to 0$. 

To formulate the tightened problem \eqref{eq:tightened_problem}, it is necessary to derive both the convex majorant $\omega^f_{r}(\mathbf x)$ associated with the cost function $f(\mathbf{x},\boldsymbol{\xi})$ and the coefficient $\omega^k_{r}$, the latter of which is the modulus of continuity of the kernel function $k(\boldsymbol\xi,\boldsymbol\zeta)$. Tables~\ref{tab:convex-continuity-majorants} and~\ref{tab:kernel-continuity-bounds} present the corresponding expressions for several common classes of cost functions and three widely used kernels, respectively. 
\begin{table}[!t]
\centering
\caption{Choices of convex continuity majorants under Assumption~\ref{ass:loss_function}}
\label{tab:convex-continuity-majorants}
\small
\renewcommand{\arraystretch}{1.5}
\begin{tabularx}{\textwidth}{@{}p{0.24\textwidth}>{\raggedright\arraybackslash}X>{\raggedright\arraybackslash}p{0.24\textwidth}@{}}
\toprule
Cost Function
& Definition
& \(\omega^f_r(\mathbf x)\) \\
\midrule
Convex piecewise-linear
&\(\displaystyle f(\mathbf x,\boldsymbol\xi) = \max_{p\in[P]} \{\mathbf x^\top \mathbf A_p\boldsymbol\xi+\mathbf b_p^\top\boldsymbol{\xi}+c_p\} \)
&\(\displaystyle r\max_{p\in[P]}\|\mathbf A_p^\top\mathbf x+\mathbf b_p\|_2 \)\\
Quadratic tracking
&\(f(\mathbf x,\boldsymbol\xi)= \|A\mathbf x-\boldsymbol\xi\|_2^2\), where \(\Xi\) is compact and \(R_\Xi:=\sup_{\boldsymbol\xi\in\Xi}\|\boldsymbol\xi\|_2<\infty\).
&
\(2r\bigl(R_\Xi+\|A\mathbf x\|_2\bigr)\)\\
Newsvendor
&\(f(\mathbf x,\boldsymbol\xi)= \mathbf h^\top(\mathbf x-\boldsymbol\xi)^+ +\mathbf b^\top(\boldsymbol\xi-\mathbf x)^+\)
&\(r\|\max\{\mathbf h,\mathbf b\}\|_2\)\\
Lipschitz-continuous
&\(f(\mathbf x,\cdot)\) is Lipschitz continuous on \(\Xi\) with a finite Lipschitz constant \(L_f(\mathbf x)\) that is convex in \(\mathbf x\).
&\(rL_f(\mathbf x)\)\\

Hölder-continuous
&\(f(\mathbf x,\boldsymbol\xi)= \|\mathbf x\|_2\|\boldsymbol\xi\|_2^\alpha\), where \(0<\alpha<1\).
&\(r^\alpha\|\mathbf x\|_2\)\\
\bottomrule
\end{tabularx}
\end{table}

\begin{table}[!ht]
\centering
\caption{Moduli of continuity for common kernels under Assumption~\ref{ass:kernel}}
\label{tab:kernel-continuity-bounds}
\small
\begin{tabular}{@{}lll@{}}
\toprule
Kernel Function & Definition  & \(\omega^k_r\) \\
\midrule

RBF &\(k(\boldsymbol\xi,\boldsymbol\zeta)=\displaystyle\exp\!\left\{-\frac{\|\boldsymbol\xi-\boldsymbol\zeta\|_2^2}{2\sigma^2}\right\}\)
&\(\displaystyle\sqrt{2-2\exp\!\left\{-\frac{r^2}{2\sigma^2}\right\}}\)
\\

Laplace &\(k(\boldsymbol\xi,\boldsymbol\zeta)=\displaystyle \exp\!\left\{-\frac{\|\boldsymbol\xi-\boldsymbol\zeta\|_2}{\sigma}\right\} \)
&\(\displaystyle \sqrt{2-2\exp\!\left\{-\frac{r}{\sigma}\right\}} \)
\\

Mat\'ern & \(k(\boldsymbol\xi,\boldsymbol\zeta)=\displaystyle M_\nu\!\left(\frac{\sqrt{2\nu}\|\boldsymbol\xi-\boldsymbol\zeta\|_2}{\sigma}\right)\)
& \(\displaystyle \sqrt{2-2M_\nu\!\left(\frac{\sqrt{2\nu}\,r}{\sigma}\right)}\)
\\

\bottomrule
\end{tabular}
\medskip
\begin{minipage}{\textwidth}
\footnotesize
\textit{Notes.} \(M_\nu(t):=\frac{2^{1-\nu}}{\Gamma(\nu)}t^\nu K_\nu(t)\) for \(t>0\) and \(\nu>0\), with \(M_\nu(0):=1\), where \(\Gamma(\cdot)\) and \(K_\nu(\cdot)\) denote the gamma function and the modified Bessel function of the second kind, respectively.
\end{minipage}
\end{table}

Given \(r>0\), the convex majorants \(\omega^f_{r}(\mathbf x)\) in Table~\ref{tab:convex-continuity-majorants} are convex in \(\mathbf x\), and \(\omega^k_r\) is a nonnegative scalar. Thus, the safety margin \(\omega^f_r(\mathbf x)+\omega^k_r\|g\|_{\mathcal H_k}\) can be evaluated explicitly and its insertion into sampled constraints preserves convexity. Particularly, in the Lipschitz-continuous case, if a common Lipschitz constant \(L_f\) applies to all \(\mathbf x\in\mathcal X\), the majorant simplifies to \(\omega^f_r(\mathbf x)=rL_f\) and is irrelevant to the decision variable. 

Next we tackle how to effectively resolve the safe approximation \eqref{eq:tightened_problem}, where an obstacle lies in the infinite-dimensional decision variable \(g\in\mathcal H_k\). By the representer theorem~\cite{zhu2021kernel}, the optimizer \(g^*\) always admits a finite-dimensional kernel representation, which eventually enables a tractable convex programming reformulation of \eqref{eq:tightened_problem}. Specifically, let
\begin{equation} \label{eq:set_upsilon}
    \Upsilon := \left ( \cup_{i\in{\rm SM}} \mathcal{D}_i \right ) \cup \Xi_N,
\end{equation}
be the set of both observed samples from support measures and discretization points. We enumerate the elements of \(\Upsilon\) as \(\{\boldsymbol{\upsilon}_l\}_{l=1}^{N_\Upsilon}\), where $N_\Upsilon = \sum_{i \in {\rm SM}}n_i+N$ is the cardinality of $\Upsilon$, and assume no two points in \(\Upsilon\) coincide without loss of generality. Then the representer theorem~\cite{scholkopf2001generalized, diwale2018generalized, zhu2021kernel} ensures that the optimizer \(g^*\) of~\eqref{eq:tightened_problem} can be  parameterized as 
\begin{equation} \label{eq:empirical_representer_g_star}
    g^* = \sum_{l=1}^{N_\Upsilon}\theta_l\phi(\boldsymbol{\upsilon}_l).
\end{equation}
Plugging \eqref{eq:empirical_representer_g_star} into the infinite-dimensional dual problem~\eqref{eq:tightened_problem} and leveraging the kernel trick, we arrive at a tractable finite-dimensional convex program:
\begin{subequations} \label{eq:tightened_problem_wo_Nystrom}
    \begin{align}
        \min_{\mathbf x\in\mathcal X, \boldsymbol{\theta}\in\mathbb{R}^{N_\Upsilon},\beta \in\mathbb{R}, t\in\mathbb R_+}~&\sum_{i\in{\rm SM}}\sum_{j=1}^{n_i}\frac{\hat{\alpha}_i^*}{n_i}\sum_{l=1}^{N_\Upsilon}\theta_lk(\boldsymbol{\upsilon}_l, \boldsymbol{\xi}_j^{(i)}) + \hat{R}^*_\gamma t + \beta\\
            \text{s.t.}~~~~~\quad\quad&\sum_{l=1}^{N_\Upsilon}\sum_{l'=1}^{N_\Upsilon}\theta_l\theta_{l'}k(\boldsymbol{\upsilon}_l,\boldsymbol{\upsilon}_{l'})\leq t^2,\label{eq:tightened_problem_wo_Nystrom_conic}\\
            &f(\mathbf x, \boldsymbol{\zeta}_l)+\omega^f_{r_N}(\mathbf x)+\omega^k_{r_N}t\leq \sum_{j=1}^{N_\Upsilon}\theta_jk(\boldsymbol{\zeta}_l, \boldsymbol{\upsilon}_j) +\beta,\quad \forall~l\in[N].
    \end{align}
\end{subequations}
Under Assumption~\ref{ass:kernel}, the Gram matrix \([k(\boldsymbol\upsilon_l,\boldsymbol\upsilon_{l'})]_{l,l'=1}^{N_\Upsilon}\) is always positive semidefinite, so constraint~\eqref{eq:tightened_problem_wo_Nystrom_conic} admits a second-order cone representation for \(t\in\mathbb R_+\). Henceforth, the resulting reformulation~\eqref{eq:tightened_problem_wo_Nystrom} is a convex program that can be handled using off-the-shelf convex optimization solvers. 
\begin{remark} \label{remark:preservation_computationally_formulation_guarantee}
    The optimal solution to the constraint-tightened problem~\eqref{eq:tightened_problem} and its convex programming reformulation~\eqref{eq:tightened_problem_wo_Nystrom} inherit the desired asymptotic decision guarantee. Let \(\mathbf{x}_{\rm UB}^*\) denote the optimizer of~\eqref{eq:tightened_problem} and~\eqref{eq:tightened_problem_wo_Nystrom}. Then under the assumptions of Proposition~\ref{prop:out_of_sample_guarantee_DRO_finite_sample}, it follows that for every \(\epsilon>0\), 
    \begin{equation}
        \lim_{M\rightarrow \infty}{\rm Pr}_{R_M}\left\{{\rm Pr}_{\mathbb{P}'\sim\Psi}\left\{\mathbb{E}_{\mathbb{P}'}\left\{f(\mathbf x^*_{\rm UB},\boldsymbol{\xi})\right\}\leq z_{\rm UB}\right\}\geq 1-\nu-\epsilon\right\}=1.
    \end{equation}
\end{remark}

\begin{remark}
    Since an exact solution of the semi-infinite problem~\eqref{eq:SMC_infinite_dual_formulation} is generally unavailable, one cannot practically evaluate the optimality loss \(z_{\rm UB}-z_{\rm D}\) in Theorem~\ref{thm:tightening_and_prior_bound}(b). Instead, \(z_{\rm UB}-z_{\rm LB}\), which is computable by solving two tractable convex programs, provides a certified upper bound on \(z_{\rm UB}-z_{\rm D}\). A small gap indicates near-optimality; otherwise, one shall increase \(N\) to shrink the covering radius \(r_N\) and tighten the certificate.
\end{remark}

\subsection{Nystr\"om approximation-enabled acceleration}
Notice that problem~\eqref{eq:tightened_problem_wo_Nystrom} has a total of $N_\Upsilon+d+2$ decision variables and can appear as a large-scale convex program. To mitigate the computational burden, we develop an acceleration strategy based on Nystr\"om approximation. Specifically, let us consider a subset $\mathcal{S}$ of $\Upsilon$:
\begin{equation} \label{eq:set_Gamma}
    \mathcal{S}:=\{\boldsymbol{\gamma}_l\}_{l=1}^{N_\mathcal{S}} \subset \Upsilon,
\end{equation}
where \(N_\mathcal{S}(<N_\Upsilon)\) is the cardinality of $\mathcal{S}$. We then confine \(g\) to a lower-dimensional subspace \(\mathcal{H}_{\mathcal{S}}:={\rm span}(\left\{\phi(\boldsymbol{\gamma}_l)\right\}_{l=1}^{N_\mathcal{S}})\) and obtain a reduced-basis approximation of \eqref{eq:empirical_representer_g_star}:
\begin{equation} \label{eq:compact_representer} 
    \tilde{g} = \sum_{l=1}^{N_\mathcal{S}}\vartheta_l\phi(\boldsymbol{\gamma}_l),
\end{equation}
with coefficients \(\boldsymbol{\vartheta}=[\vartheta_1,\cdots,\vartheta_{N_\mathcal{S}}]^\top\in\mathbb{R}^{N_\mathcal{S}}\).
By plugging \eqref{eq:compact_representer} in lieu of~\eqref{eq:empirical_representer_g_star} into the tightened problem~\eqref{eq:tightened_problem}, we arrive at the following approximation of~\eqref{eq:tightened_problem_wo_Nystrom}:
\begin{equation} \label{eq:tightened_problem_Nystrom}
    \begin{aligned}
        z_{\rm UB+}:=\min_{\mathbf x\in\mathcal X, \boldsymbol{\vartheta}\in\mathbb{R}^{N_\mathcal{S}},\beta \in\mathbb{R}, t\in\mathbb R_+}&~\sum_{i\in{\rm SM}}\sum_{j=1}^{n_i}\frac{\hat{\alpha}_i^*}{n_i}\sum_{l=1}^{N_\mathcal{S}}\vartheta_lk(\boldsymbol{\gamma}_l, \boldsymbol{\xi}_j^{(i)}) + \hat{R}^*_\gamma t + \beta\\
            \text{s.t.}~~~\qquad~&\sum_{l=1}^{N_\mathcal{S}}\sum_{l'=1}^{N_\mathcal{S}}\vartheta_l\vartheta_{l'}k(\boldsymbol{\gamma}_l,\boldsymbol{\gamma}_{l'})\leq t^2,\\
            &f(\mathbf x, \boldsymbol{\zeta}_j)+\omega^f_{r_N}(\mathbf x)+\omega^k_{r_N}t\leq \sum_{l=1}^{N_\mathcal{S}}\vartheta_lk(\boldsymbol{\zeta}_j, \boldsymbol{\gamma}_l) +\beta, \quad \forall~j\in[N],
    \end{aligned}
\end{equation}
which naturally yields \(z_{\rm UB+} \geq z_{\rm UB} \). Thanks to the Nystr\"om approximation, \eqref{eq:tightened_problem_Nystrom} has \(N_\mathcal{S}+d+2\) decision variables and is thus easier to solve than \eqref{eq:tightened_problem_wo_Nystrom}. This computational benefit turns out to induce a controlled loss of optimality, as formally stated below.

\begin{proposition} \label{prop:gap_bound} 
Let \((\hat{\mathbf x}^*,\hat{g}^*, \hat{\beta}^*, \hat t^*)\) denote the optimal solution to problem~\eqref{eq:tightened_problem_wo_Nystrom}, where the optimizer admits the finite expansion \(\hat{g}^*=\sum_{l=1}^{N_\Upsilon}\theta_l^*\phi(\boldsymbol{\upsilon}_l)\) with \(\boldsymbol{\theta}^*\in\mathbb{R}^{N_\Upsilon}\). Let \(\Upsilon\) and \(\mathcal{S}\) be as in \eqref{eq:set_upsilon} and \eqref{eq:set_Gamma}, respectively, and define the Gram blocks by \(\mathbf{K}_{\mathcal{S}\Upsilon}:=[k(\boldsymbol{\gamma}_i,\boldsymbol{\upsilon}_j)]_{i=1}^{N_\mathcal{S}}{}_{j=1}^{N_\Upsilon}\in\mathbb{R}^{N_\mathcal{S}\times N_\Upsilon}, \mathbf{K}_{\mathcal{S}\mathcal{S}}:=[k(\boldsymbol{\gamma}_i,\boldsymbol{\gamma}_j)]_{i=1}^{N_\mathcal{S}}{}_{j=1}^{N_\mathcal{S}}\in\mathbb{R}^{N_\mathcal{S}\times N_\mathcal{S}}\), and \(\mathbf{K}_{\Upsilon\Upsilon}:=[k(\boldsymbol{\upsilon}_i,\boldsymbol{\upsilon}_j)]_{i=1}^{N_\Upsilon}{}_{j=1}^{N_\Upsilon}\in\mathbb{R}^{N_\Upsilon \times N_\Upsilon}\).
Then, solving problem~\eqref{eq:tightened_problem_Nystrom} under the Nystr\"om approximation yields the following optimality gap:
\begin{equation}
    z_{\rm UB+}-z_{\rm UB}\leq 2\sqrt C\|\boldsymbol{\theta}^*\|_2\sqrt{\|\mathbf{R}_{\mathcal S} \|_2},
\end{equation}
where \(\mathbf{R}_{\mathcal S} := \mathbf{K}_{\Upsilon\Upsilon}-\mathbf{K}_{\mathcal{S}\Upsilon}^\top\mathbf{K}^\dagger_{\mathcal{S}\mathcal{S}}\mathbf{K}_{\mathcal{S}\Upsilon}\) denotes the Nystr\"om residual. 
\end{proposition}

Proposition \ref{prop:gap_bound} indicates that the loss of optimality rests on how well the subset $\mathcal{S}$ approximates $\Upsilon$. In machine learning and statistics, a large number of random sampling and clustering algorithms are available for selecting a subset from a large dataset. To reduce approximation errors, we design a two-step subsampling heuristic by exploiting the structural characteristics of $\Upsilon$ in \eqref{eq:set_upsilon}. Notice that the data points in $\mathcal{D}_i$ can manifest spatial similarity, while the discretization points in $\Xi_N$ are uniformly distributed over \(\Xi\). Henceforth, we first subsample each \(\mathcal D_i\) and merge the resulting subsets with \(\Xi_N\) to form an intermediate dataset $\mathcal{S}'$. Then we progressively subsample $\mathcal{S}'$ until its relative spectral-norm residual of the Nystr\"om approximation fulfills \(\rho_{\mathcal S}:=\|\mathbf R_\mathcal S\|_2/{\|\mathbf K_{\Upsilon\Upsilon}\|_2}\leq \varepsilon\), where \(\varepsilon > 0\) is a prescribed tolerance. For efficient subsampling, we suggest using the kernel herding algorithm~\cite{Chen2010herding, bach2017equivalence} due to its promise in approximating the underlying distribution in RKHS. In practical implementation, we perform kernel herding in batches of \(T\) iterations on each \(\mathcal D_i\), updating \(\mathcal S'\) and checking its Nystr\"om residual relative to \(\Upsilon\) after each batch, until \(\rho_{\mathcal S'}\leq \varepsilon/2\) and follow the residual-based stopping criterion when subsampling $\mathcal{S}'$. Finally, we present the following algorithm that streamlines the entire proposed RooD-SO framework.

\begin{algorithm}[!ht]
  \caption{Data-driven solution procedure for RooD-SO}
  \label{alg:1}
  \begin{algorithmic}[1]
    \Require \(M\) sources of data \(\{\mathcal{D}_{i}\}_{i=1}^M\) where \(\mathcal{D}_i=\{\boldsymbol{\xi}_j^{(i)}\}_{j=1}^{n_i}\), first-stage herding step \(T\), Nystr\"om residual threshold \(\varepsilon\), number of discretization points \(N\), regularization parameter \(\nu\), enlargement parameter \(\gamma\), compact support set \(\Xi\)
    \State Calculate the Gram matrix \(\hat{\mathbf{K}}\) using datasets \(\{\mathcal D_i\}_{i=1}^M\).
    \State Solve the dual problem~\eqref{eq:dual_SMC_quadratic}, and derive the optimizer \(\hat{\boldsymbol{\alpha}}^*\) as well as the radius \(\hat{R}^*\) as in~\eqref{eq:R_square}.
    \State Generate \(N\) discretization points \(\Xi_N=\{\boldsymbol{\zeta}_j\}_{j=1}^N\) using space-filling design methods over \(\Xi\), such as Latin hypercube sampling or regular grid design. 
    \State Calculate \(r_N\), evaluate \(\omega^k_{r_N}\) and \(\omega^f_{r_N}(\mathbf x)\).
    \State Apply kernel herding to each \(\mathcal D_i\), \(i\in\mathrm{SM}\), in batches of \(T\) iterations, updating \(\mathcal S'\) by merging the accumulated subsets with \(\Xi_N\) after each round, until \(\rho_{\mathcal S'}\leq\varepsilon/2\).
    \State Perform kernel herding on \(\mathcal{S}'\) to select points until \(\rho_{\mathcal S}\leq \varepsilon\), and obtain \(\mathcal{S}\).
    \State Formulate the Nystr\"om approximation-based problem~\eqref{eq:tightened_problem_Nystrom}, and then solve it using off-the-shelf convex programming solvers.
  \end{algorithmic}
\end{algorithm}

\begin{remark}
    The solution returned by Algorithm~\ref{alg:1} retains the asymptotic out-of-distribution performance guarantee established in Proposition~\ref{prop:out_of_sample_guarantee_DRO_finite_sample}. Because $z_{\rm UB+} \geq z_{\rm UB}$, an argument analogous to Remark~\ref{remark:preservation_computationally_formulation_guarantee} implies that for every \(\epsilon>0\),
    \begin{equation}
        \lim_{M\rightarrow \infty}{\rm Pr}_{R_M}\left\{{\rm Pr}_{\mathbb{P}'\sim\Psi}\left\{\mathbb{E}_{\mathbb{P}'}\left\{f(\mathbf{x}_{\rm UB+}^*,\boldsymbol{\xi})\right\}\leq z_{\rm UB+}\right\}\geq 1-\nu-\epsilon\right\}=1,
    \end{equation}
    where \(\mathbf{x}_{\rm UB+}^*\) denotes the optimizer of problem~\eqref{eq:tightened_problem_Nystrom} under Nystr\"om approximation. 
\end{remark}

\section{Computational experiments}
\label{sec:experiments}

In this section, we carry out comprehensive computational experiments on robust decision-making tasks under unobserved target distributions, using only multiple data sources as prior information. To implement our proposed method \textbf{RooD-SO}, the RBF kernel \(k(\boldsymbol{\xi},\boldsymbol{\xi}')=\exp\left\{-\|\boldsymbol{\xi}-\boldsymbol{\xi}'\|^2/{(2\sigma^2)}\right\}\) is employed, where the bandwidth parameter \(\sigma\) is determined using the median heuristic~\cite{gretton2012kernel}. Besides, we adopt a variety of known data-driven decision-making models for comprehensive comparisons. The following data-centric optimization methods based on distribution averaging are implemented: (i) \textbf{SAA}~\cite{shapiro2021lectures}: the classical SAA scheme applied to the pooled dataset \(\mathcal{D}^{\cup}\); (ii) \textbf{DRO-W}~\cite{Kuhn2019WassersteinDistance}: DRO based on a 1-Wasserstein ball centered at the empirical distribution \(\hat{\mathbb{P}}^{\cup}\) defined by \(\mathcal D^\cup\); (iii) \textbf{SAA-WB}~\cite{cuturi2014fast}: an SAA variant based on the WB of \(\{\hat{\mathbb{P}}_i\}_{i=1}^M\); (iv) \textbf{DRO-K}~\cite{zhu2021kernel}: DRO based on an RKHS-ball centered at the empirical distribution \(\hat{\mathbb{P}}^{\cup}\). Besides, the following DRO models explicitly handling inter-source heterogeneity are implemented: (v) \textbf{DRO-Intx}~\cite{rychener2024wasserstein}: DRO with an uncertainty set defined as the intersection of \(M\) 1-Wasserstein balls, each centered at a source distribution~\eqref{eq:intersection_ambiguity_set}; (vi) \textbf{DRO-Conv}~\cite{zhang2023stochastic}: DRO whose uncertainty set is the convex hull of \(M\) distributions, i.e., \(\mathcal{U}^{\rm co}\) in~\eqref{eq:Group_DRO_ambiguity_set}.

We conduct comparisons on the multi-item newsvendor problem~\cite{turken2012multi} and the portfolio optimization problem~\cite{krokhmal2002portfolio}. All the experiments are conducted using Python 3.13.13 and COPT v8.0.6, implemented on a rack server with dual AMD EPYC 9334 32-core processors @ up to 3.91 GHz and 256 GB RAM.

\subsection{Newsvendor problem}
We consider a \(d\)-item newsvendor problem where a retailer commits to an order vector \(\mathbf{x}\in \mathbb{R}^d\) prior to observing stochastic demand \(\boldsymbol{\xi}\in \mathbb{R}^d\) of a target product. The following cost function is employed to balance sales, holding, and backorder costs:
\begin{equation}\label{eq:multi_newsvendor}
    f(\mathbf{x},\boldsymbol{\xi})=-\mathbf{p}^\top\min\{\mathbf{x},\boldsymbol{\xi}\}+\mathbf{c}^\top\mathbf{x}+\mathbf{h}^\top(\mathbf{x}-\boldsymbol{\xi})^++\mathbf{b}^\top(\boldsymbol{\xi}-\mathbf{x})^+,
\end{equation}
where parameters \(\mathbf{p}\), \(\mathbf{c}\), \(\mathbf{h}\), and \(\mathbf{b}\) denote the per-unit selling price, ordering cost, holding cost, and backorder cost for each item, respectively. It is assumed that there are no historical data on demand $\boldsymbol{\xi}$ for the target product, but we have access to datasets \(\mathcal{D}_1,\cdots,\mathcal{D}_M\) with respect to \(M\geq 2\) relevant products. 

\subsubsection{Case studies on synthetic data}
\label{sec:synthetic_newsvendor}
We consider a simple two-item instance, i.e., \(d=2\), with \(\mathbf{p} = [5,4],~\mathbf{c}=[2,1.5],~ \mathbf{h}=[1, 1.6]\), and \(\mathbf{b}=[1, 0.8]\). We generate heterogeneous source distributions using a hierarchical mixture model. For each source \(i\), we independently draw a latent source type $\mathbf{z}_i$ with $\Pr\{\mathbf{z}_i=(0,0)\}=0.8$ and $\Pr\{\mathbf{z}_i=(-1,0)\}=\Pr\{\mathbf{z}_i=(0,-1)\}=0.1$, and an independent perturbation \(\boldsymbol{\delta}_i\sim{\rm Unif}([-0.02, 0.02]^2)\). We then define the mixing probabilities for source $i$ as \(\mathbf q_i=0.5\boldsymbol{1}+0.4(\mathbf z_i-\mathbb E\{\mathbf z_i\})+\boldsymbol{\delta}_i\), where \(\mathbb E\{\mathbf z_i\}=(-0.1,-0.1)\). Conditional on \(\mathbf q_i\), the demand distribution of source \(i\) is \(\mathbb P_i:=\bigotimes_{j=1}^2(1-q_{ij}){\rm Unif}([2,3])+q_{ij}{\rm Unif}([3.4,4.4])\), and we draw \(n_i\) conditionally independent samples from \(\mathbb P_i\) to form \(\mathcal D_i\). 

For a fair comparison, all hyperparameters are carefully calibrated. For all conservatism-related parameters, including parameters \(\nu\) and \(\gamma\) in RooD-SO and the radii of Wasserstein balls in DRO-W and DRO-Intx, an independent validation split is used to select their values. We perform five-fold cross-validation by partitioning the available data from each source into five subsets of equal size. In each fold, four subsets are used for training and the remaining subset for validation. We conduct a grid search and select the hyperparameters achieving the best average validation performance across the five folds. For RooD-SO, we search over the interval $(0,1)$ for $\nu$ with a step size of $0.05$ and select $\gamma$ from $\{0.001,0.01,0.1,1\}$. We select the radii of the Wasserstein balls in DRO-W and the RKHS balls in DRO-K from $\{0, 0.001, 0.002, 0.005, 0.01, 0.02, 0.05, 0.1, 0.2,0.4\}$. When \(M=2\) in DRO-Intx, we set the Wasserstein radii to \(\epsilon_1=\lambda_1(1+\lambda_2)\mathcal{W}_1(\hat{\mathbb{P}}_1,\hat{\mathbb{P}}_2)\) and \(\epsilon_2=(1-\lambda_1)(1+\lambda_2)\mathcal{W}_1(\hat{\mathbb{P}}_1,\hat{\mathbb{P}}_2)\), where \(\lambda_1\) ranges from 0 to 1 in steps of 0.1, and \(\lambda_2\in\{0.002,0.005,0.01,0.02\}\), as in~\cite{rychener2024wasserstein}. For \(M>2\), it is computationally prohibitive to tune \(M\) radii jointly, so in order to ensure a nonempty intersection of \(M\) Wasserstein balls, we heuristically set the radius of each Wasserstein ball to the largest pairwise 1-Wasserstein distance among \(\{\hat{\mathbb{P}}_i\}_{i=1}^M\). To calculate the barycenter in SAA-WB, we follow~\cite{cuturi2014fast} by aligning the number of barycenter support points \(n_W\) with the average sample size of the \(M\) source datasets, i.e., \(n_W = \lfloor \frac{1}{M}\sum_{i=1}^M n_i \rfloor\). For RooD-SO, we set \(N=2000\), \(T=50\), and \(\varepsilon=1\%\). By Table~\ref{tab:convex-continuity-majorants}, the newsvendor cost function admits the explicit convex continuity majorant \(\omega_r^f(\mathbf x)= r\|\max(\mathbf p+\mathbf h,\mathbf b)\|_2\). For all methods, we set a maximal solution time of 1,800 seconds. 


\begin{figure}[!ht]
    \begin{subfigure}[b]{0.32\linewidth}
        \centering
                \includegraphics[width=\linewidth]{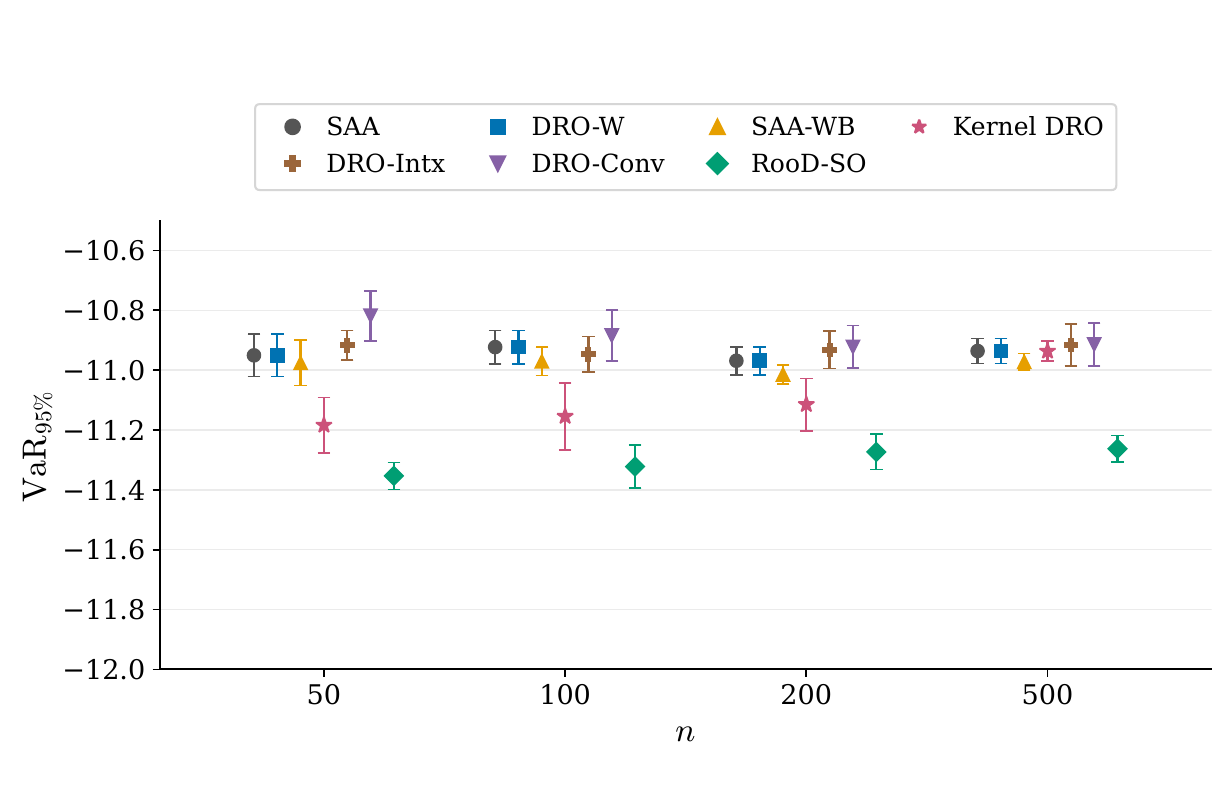}
        \caption{\(M=2\)}
        \label{fig:rank500_r040_M2}
      \end{subfigure}
      \hspace{0.02\linewidth}
  \begin{subfigure}[b]{0.32\linewidth}
    \centering
    \includegraphics[width=\linewidth]{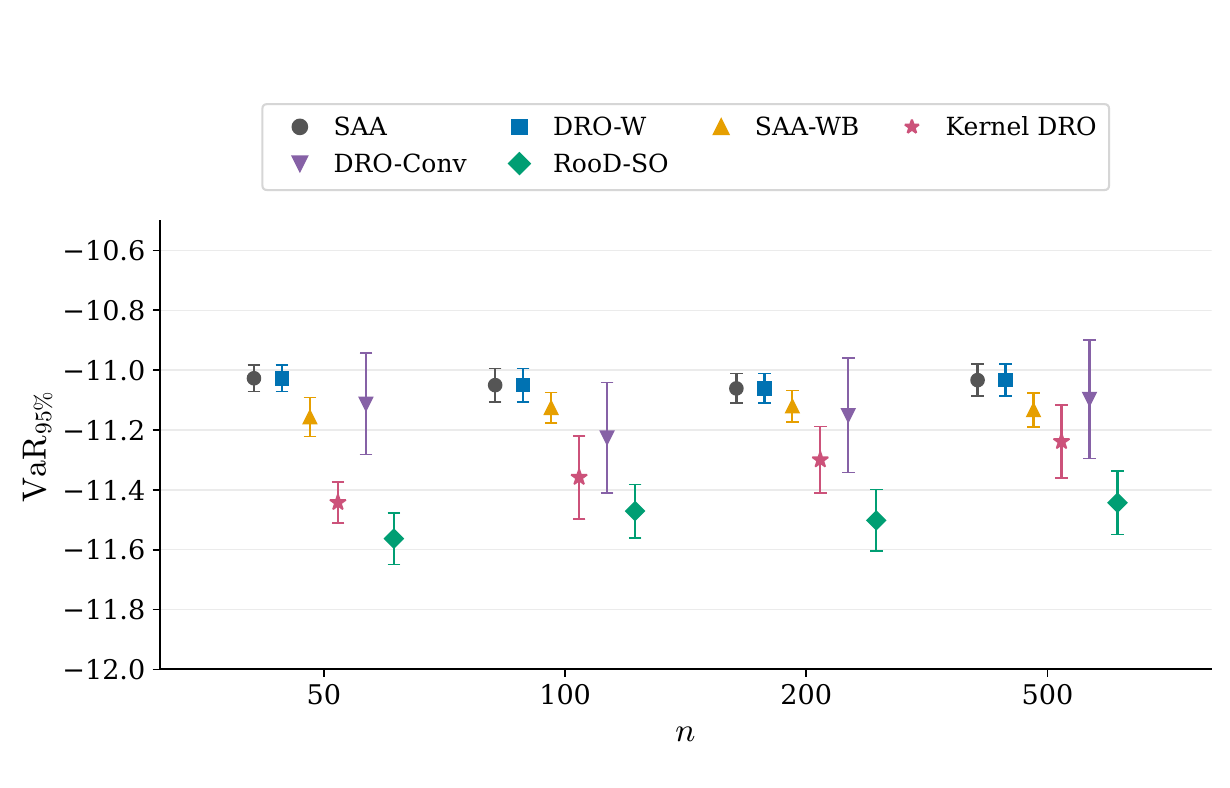}
    \caption{\(M=5\)}
    \label{fig:rank500_r040_M5}
  \end{subfigure}
    \begin{subfigure}[b]{0.32\linewidth}
    \centering
    \includegraphics[width=\linewidth]{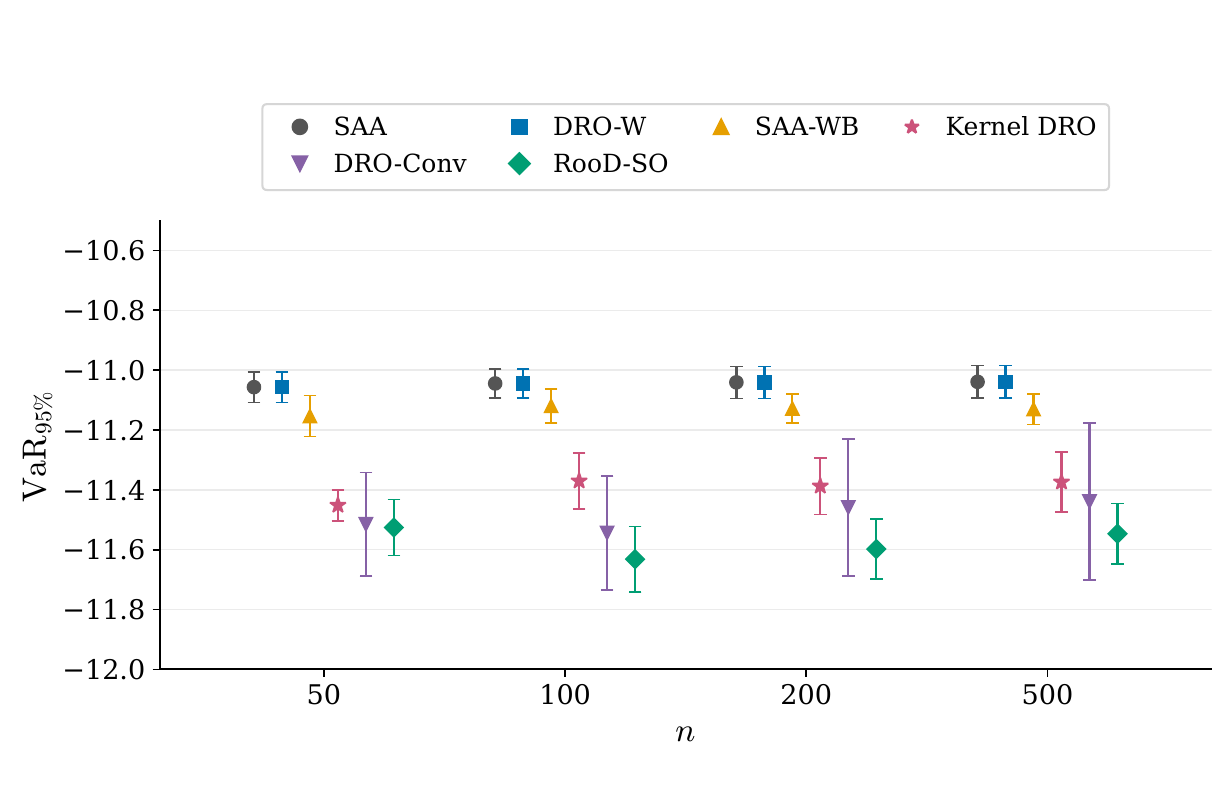}
    \caption{\(M=10\)}
    \label{fig:rank500_r040_M10}
  \end{subfigure}
  \caption{Comparisons of out-of-distribution performance across different \((M,n)\) configurations.}
  \label{fig:Var95_M}
\end{figure} 

We consider a balanced sampling regime by setting \(n_1=\cdots=n_M=n \in \{50,100,200,500\}\) and let $M\in\{2,5,10\}$, yielding 12 configurations of \((M,n)\). To evaluate the out-of-distribution robustness, we generate $10{,}000$ test distributions using the sampling procedure as aforementioned, and draw $100{,}000$ samples from each distribution. Figure~\ref{fig:Var95_M} reports the out-of-distribution performance across different $(M,n)$ configurations, evaluated by Value-at-Risk
at the $95\%$ level ($\mathrm{VaR}_{95\%}$) over test distributions through 10 independent replications. It is seen that our RooD-SO achieves the best overall out-of-distribution performance. This excellence is attributed to the modeling power of \(\mathcal{U}_{\nu,\gamma}^{\rm SMC}\) brought by learning and the robustness provided by both \(\nu\) and \(\gamma\). In contrast, methods based on distribution averaging, i.e., SAA, DRO-W, and SAA-WB, exhibit inferior performance compared to those explicitly considering source heterogeneity. Note that SAA and DRO-W are tied for the worst, which is consistent with the theoretical finding in \cite[Remark 6.7]{mohajerin2018data} that DRO-W boils down to SAA for the newsvendor problem. Although DRO-Conv explicitly encodes source heterogeneity, it suffers from inflexibility in calibrating the robustness and thus shows heavier conservatism. As for DRO-Intx, we do not report its results for \(M \in \{ 5,10 \} \) due to its prohibitive computational cost, which grows exponentially with \(M\) owing to the Cartesian product across sources, i.e., \(\mathcal{O}(\prod_{i=1}^M n_i)\)~\cite{rychener2024wasserstein}.


\begin{figure}[!ht]
    \begin{subfigure}[b]{0.32\linewidth}
        \centering
                \includegraphics[width=\linewidth]{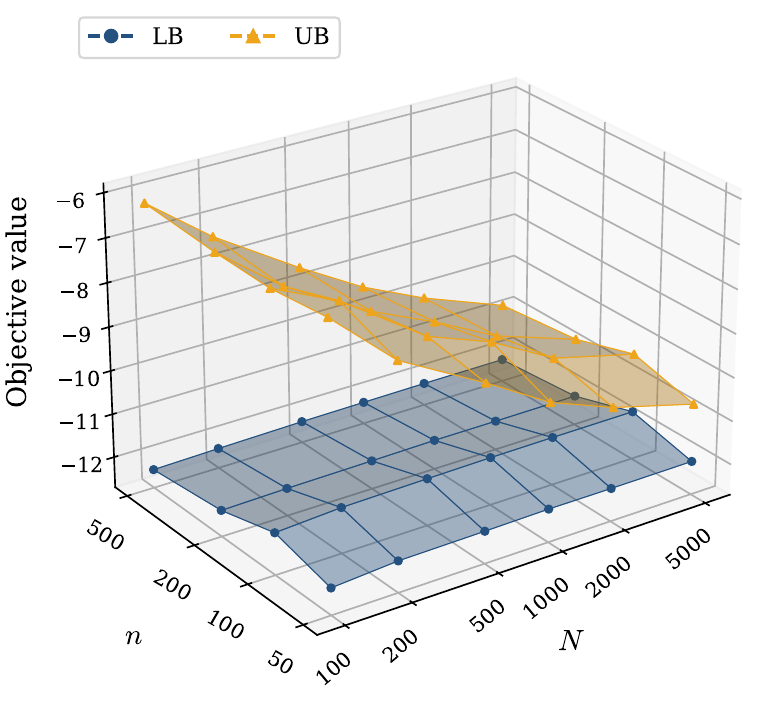}
        \caption{\(M=2\)}
        \label{fig:objective_M2}
      \end{subfigure}
      \hspace{0.02\linewidth}
  \begin{subfigure}[b]{0.32\linewidth}
    \centering
    \includegraphics[width=\linewidth]{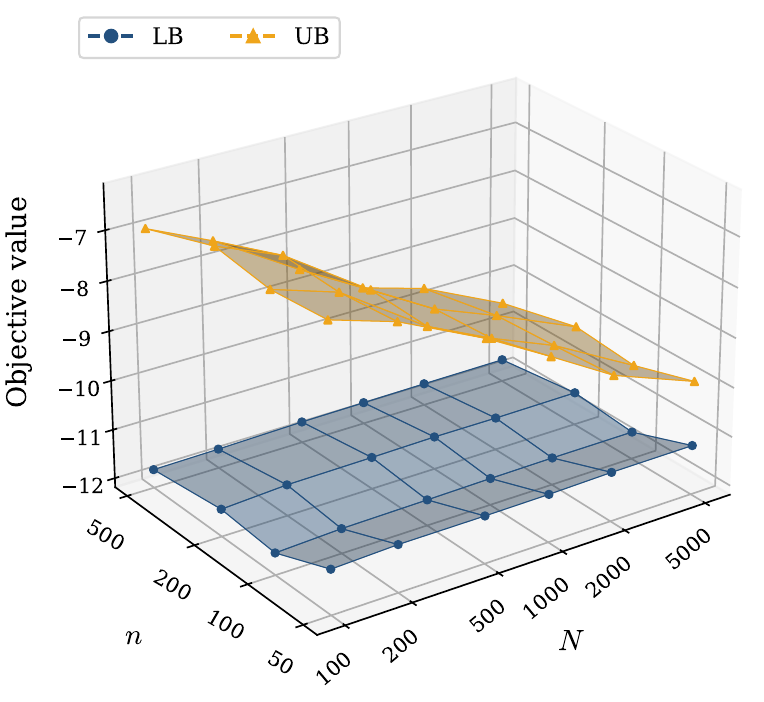}
    \caption{\(M=5\)}
    \label{fig:objective_M5}
  \end{subfigure}
    \begin{subfigure}[b]{0.32\linewidth}
    \centering
    \includegraphics[width=\linewidth]{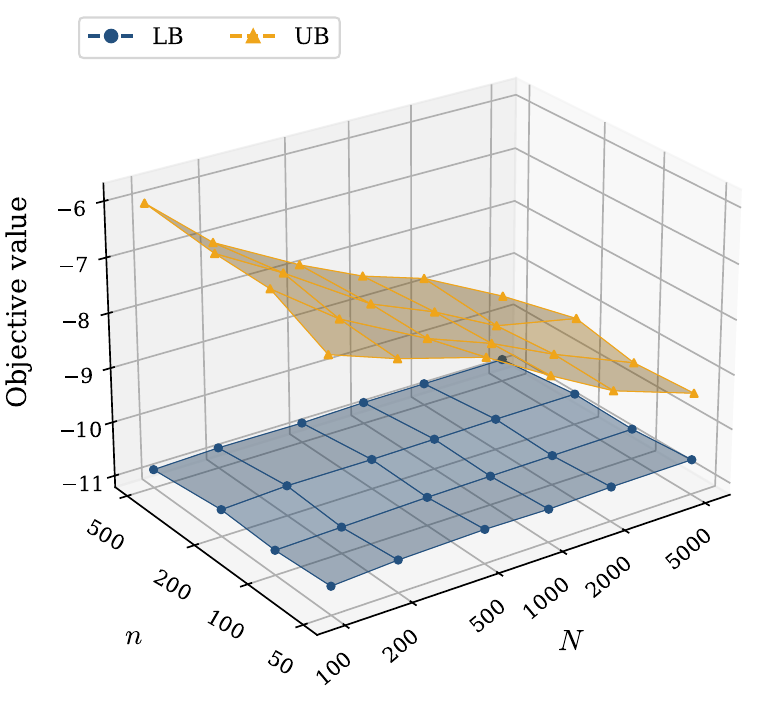}
    \caption{\(M=10\)}
    \label{fig:objective_M10}
  \end{subfigure}
  \caption{Objective values of LB and UB across different \((M,n,N)\) configurations.}
  \label{fig:objective_N}
\end{figure}

Since the solution of RooD-SO problems hinges on safe constraint-tightening approximation, we investigate how the discretization size $N$ affects optimality loss by varying \(N\) from \(100\) to \(5{,}000\) under the previous setting. Figure~\ref{fig:objective_N} depicts the objective value \(z_{\rm UB}\) of the tightened formulation~\eqref{eq:tightened_problem_wo_Nystrom} without Nystr\"om approximation. As \(N\) increases, \(z_{\rm UB}\) decreases evidently, implying an effective reduction in the conservatism associated with the margin and the loss of optimality. Since it is not possible to solve the semi-infinite program~\eqref{eq:SMC_infinite_dual_formulation} and obtain its optimal value $z_{\rm D}$, $z_{\rm LB}$ is adopted as a computable lower bound and reported in Figure~\ref{fig:objective_N}. It can be seen that as $N$ grows, there is a narrower gap between \(z_{\rm UB}\) and $z_{\rm LB}$ that upper bounds the true optimality loss. 

\begin{figure}[!ht]
    \begin{subfigure}[b]{0.32\linewidth}
        \centering
        \includegraphics[width=\linewidth]{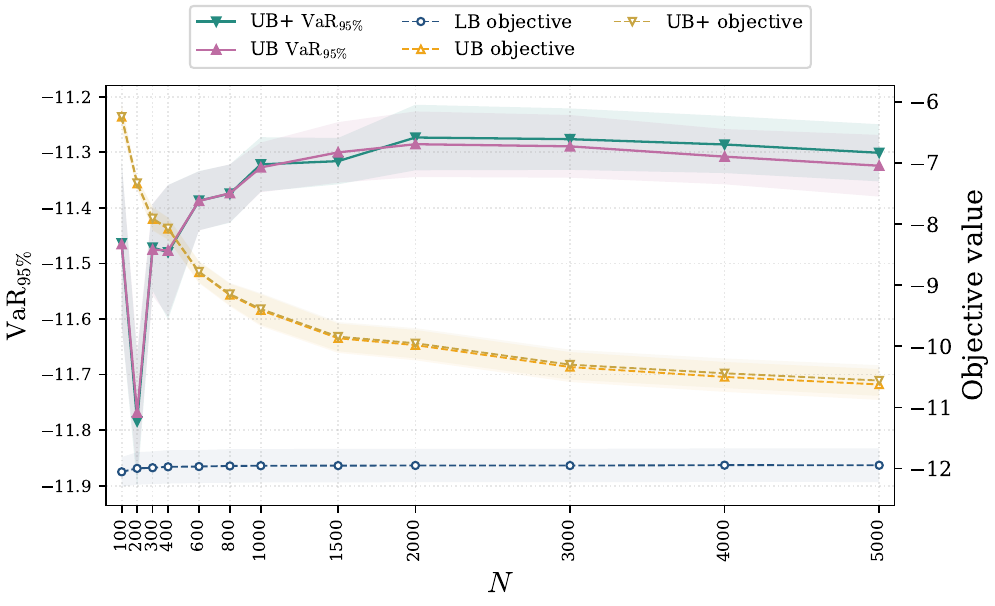}
        \caption{\(M=2\)}
        \label{fig:ooD_M2}
      \end{subfigure}
      \hspace{0.02\linewidth}
  \begin{subfigure}[b]{0.32\linewidth}
    \centering
    \includegraphics[width=\linewidth]{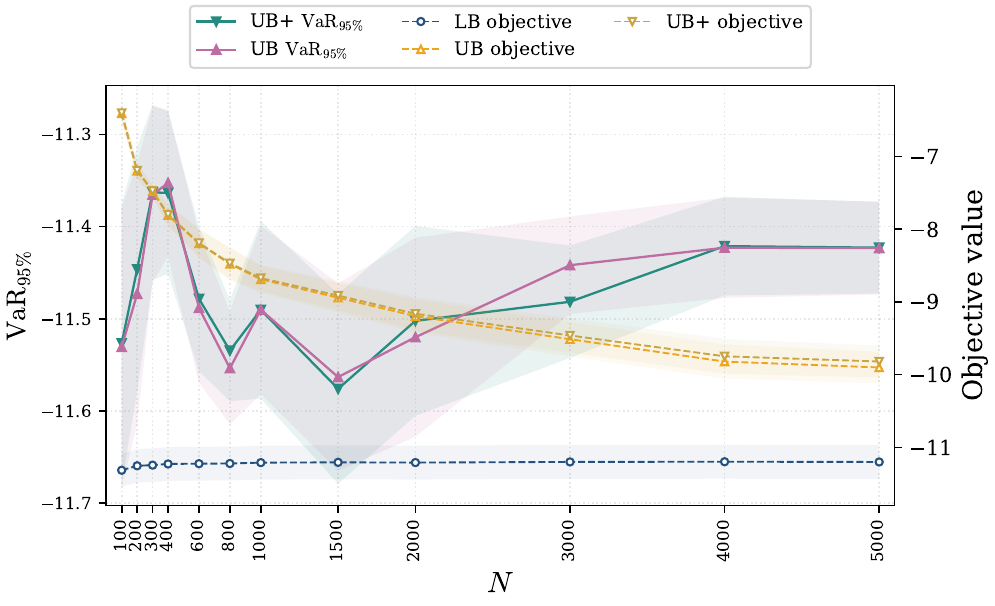}
    \caption{\(M=5\)}
    \label{fig:ooD_M5}
  \end{subfigure}
    \begin{subfigure}[b]{0.32\linewidth}
    \centering
    \includegraphics[width=\linewidth]{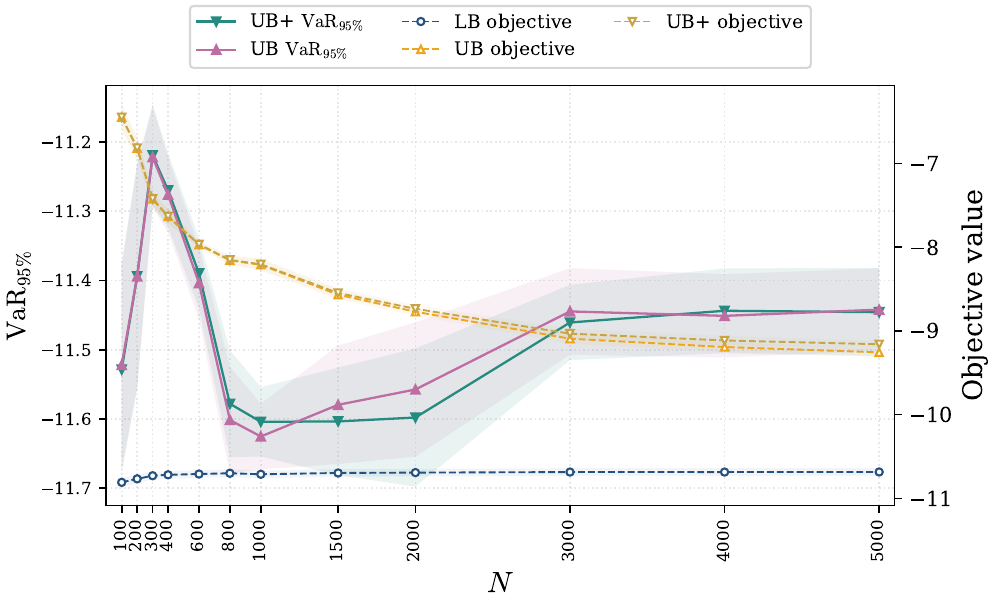}
    \caption{\(M=10\)}
    \label{fig:ooD_M10}
  \end{subfigure}
  \caption{Out-of-distribution performance and objective values for LB, UB, and UB+ across different \((M,N)\) configurations with \(n=200\).}
  \label{fig:ooD_M}
\end{figure}

\begin{figure}[!ht]
    \begin{subfigure}[b]{0.32\linewidth}
        \centering
        \includegraphics[width=\linewidth]{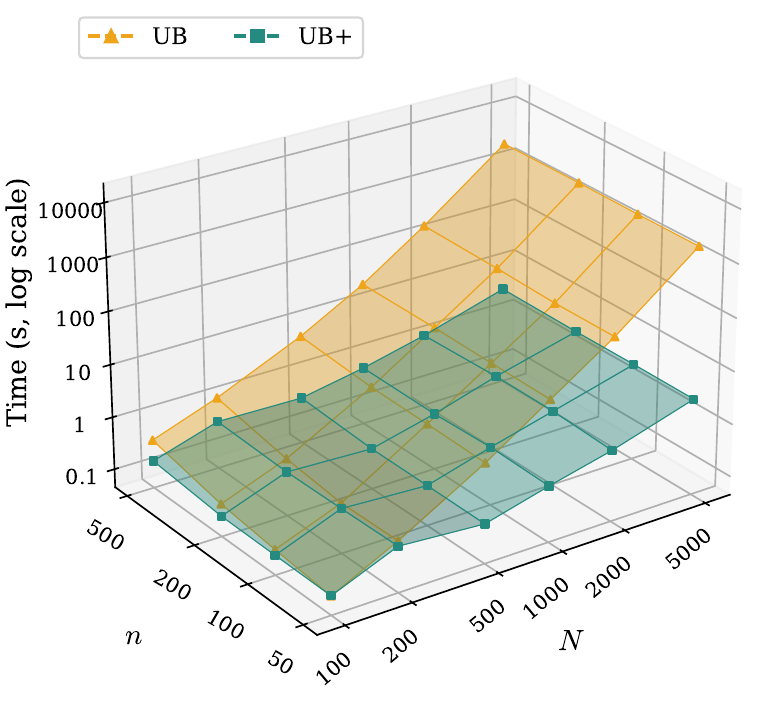}
        \caption{\(M=2\)}
        \label{fig:effect_of_N_M2}
      \end{subfigure}
      \hspace{0.02\linewidth}
  \begin{subfigure}[b]{0.32\linewidth}
    \centering
    \includegraphics[width=\linewidth]{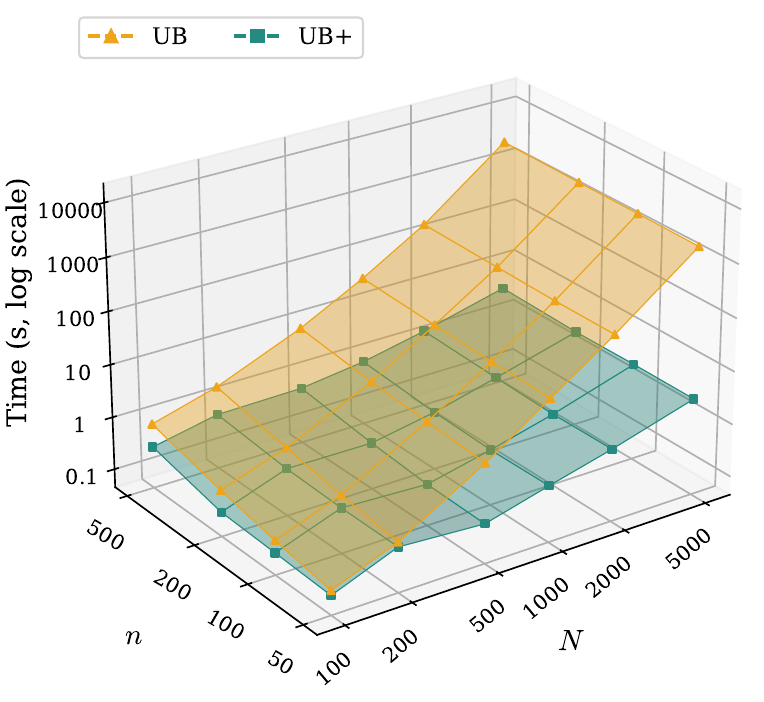}
    \caption{\(M=5\)}
    \label{fig:effect_of_N_5}
  \end{subfigure}
    \begin{subfigure}[b]{0.32\linewidth}
    \centering
    \includegraphics[width=\linewidth]{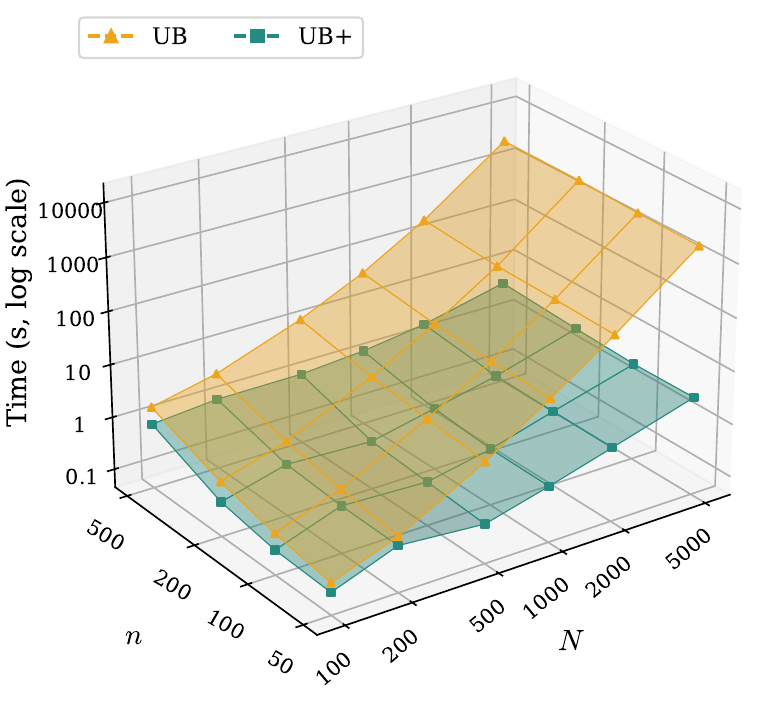}
    \caption{\(M=10\)}
    \label{fig:effect_of_N_M10}
  \end{subfigure}
  \caption{Computational times of UB and UB+ across different \((M,n,N)\) configurations.}
  \label{fig:effect_of_N}
\end{figure}

In practice, the solution process can be accelerated using the Nystr\"om approximation~\eqref{eq:tightened_problem_Nystrom}, so we further examine its impact on solution quality and out-of-distribution performance. For clearer exposition, we focus on the \(n=200\) setting. Alongside \(z_{\rm LB}\) and \(z_{\rm UB}\), Figure~\ref{fig:ooD_M} depicts the objective value \(z_{\rm UB+}\) obtained with Nystr\"om approximation~\eqref{eq:tightened_problem_Nystrom}, and VaR$_{95\%}$ achieved by the decisions from the two tightened formulations, evaluated on $10{,}000$ test distributions with $100{,}000$ samples each. We notice that the UB and UB+ formulations yield highly similar objective values and comparable out-of-distribution performance under varying $N$. This indicates that the Nystr\"om approximation does not practically compromise out-of-distribution performance. Meanwhile, a stable out-of-distribution performance of the robust solution can be achieved by choosing \(N\) sufficiently large. Furthermore, we examine the computational speedup incurred by Nystr\"om approximation. Figure~\ref{fig:effect_of_N} depicts the runtimes with and without the Nystr\"om approximation, where substantial speedups can be observed, especially under large $N$. This is because for the multi-item newsvendor problem, the resulting convex program \eqref{eq:tightened_problem_wo_Nystrom} is a second-order cone program, whose solution complexity via interior-point methods grows cubically with the number of decision variables~\cite{POTRA2000281}, so the acceleration benefit brought by Nystr\"om approximation appears to be more pronounced for larger $N$. Moreover, for all test configurations, the tightened formulation~\eqref{eq:tightened_problem_Nystrom} under Nystr\"om approximation can be solved within 5s. These results together show that the Nystr\"om approximation considerably reduces computational cost while maintaining satisfactory solution quality.

\subsubsection{Out-of-distribution performance on \emph{Coffee Shop Sales} dataset}
\label{sec:experiment_coffee_sales}
Next, we turn to a more realistic coffee shop retail setting, which echoes the instance in Figure~\ref{fig:top}. We consider a cold-start inventory planning problem for a newly launched coffee product, where the initial order quantities are constrained by both an upfront budget and a total storage capacity limit~\cite{buzacott2004inventory,chao2008dynamic}. These constraints are practically relevant because a retailer may adopt a cautious initial inventory commitment to limit early-stage downside exposure, analogous to the survival-oriented inventory policies of start-up firms \cite{archibald2002should}. Specifically, the feasible set of initial stocking decisions is defined as \(\mathcal X=\left\{ \mathbf{x}\in\mathbb{R}^{d} \,\middle|\, \boldsymbol{\kappa}^{\top}\mathbf{x}\leq B,\boldsymbol{1}^{\top}\mathbf{x}\leq Q_{\max},\mathbf{x}\geq\boldsymbol{0}\right\}\).

We use the publicly available \emph{Coffee Shop Sales} dataset, which comprises transaction records for Maven Roasters, a coffee retailer with three locations in New York City.\endnote{The \emph{Coffee Shop Sales} dataset is available at \url{https://mavenanalytics.io/data-playground/coffee-shop-sales}.} We focus on five coffee products, each characterized by its product type and base coffee origin. For each product, we compute average daily sales per week for small, regular, and large cup sizes, which form a three-dimensional vector of uncertainty \(\boldsymbol{\xi}=(\xi_{\rm S}, \xi_{\rm M}, \xi_{\rm L})\). Considering data scarcity in this case, we evaluate all methods in a leave-one-out fashion: each product is treated in turn as the unseen target domain, and the remaining four coffee products serve as available relevant data sources. In each run, an order decision for the target coffee product is informed by historical sales of other products, and subsequently tested on the held-out target product data in a zero-shot context. Further experimental details are provided in Appendix~\ref{app:coffee_experiment_details}.

\begin{table}[!ht]
\centering
\caption{Costs across Five Leave-One-Out Targets}
\label{tab:cost_coffee_5sets}
\resizebox{\linewidth}{!}{%
\begin{tabular}{ccccccccccccccc}
\toprule
\multirow[b]{2}{*}{\raisebox{-0.6ex}{\shortstack{Target\\[-1pt]Product}}} & \multicolumn{2}{c}{SAA} & \multicolumn{2}{c}{DRO-W} & \multicolumn{2}{c}{SAA-WB} & \multicolumn{2}{c}{DRO-K} & \multicolumn{2}{c}{DRO-Conv} & \multicolumn{2}{c}{RooD-SO}\\
\cmidrule(lr){2-3} \cmidrule(lr){4-5} \cmidrule(lr){6-7} \cmidrule(lr){8-9} \cmidrule(lr){10-11} \cmidrule(lr){12-13} 
& Cost & Rank & Cost & Rank & Cost & Rank & Cost & Rank & Cost & Rank & Cost & Rank \\
\midrule
Product 1 & -56.15 & 4 & -56.15 & 4 & -56.27 & 2 & -55.32 & 6 & \textbf{-56.36} & \textbf{1} & -56.20 & 3 \\
Product 2 & -55.60 & 4 & -55.60 & 4 & -55.66 & 3 & -55.75 & 2 & -55.43 & 6 & \textbf{-55.76} & \textbf{1} \\
Product 3 & -54.58 & 4 & -54.58 & 4 & -54.62 & 3 & \textbf{-54.79} & \textbf{1} & -54.47 & 6 & -54.72 & 2 \\
Product 4 & \textbf{-53.79} & \textbf{1} & \textbf{-53.79} & \textbf{1} & -53.57 & 5 & -53.48 & 6 & -53.60 & 3 & -53.58 & 4 \\
Product 5 & -54.17 & 3 & -54.17 & 3 & -53.92 & 6 & \textbf{-54.26} & \textbf{1} & -54.09 & 5 & -54.24 & 2 \\
\midrule
Avg. Rank & \multicolumn{2}{c}{3.20} & \multicolumn{2}{c}{3.20} & \multicolumn{2}{c}{3.80} & \multicolumn{2}{c}{3.20} & \multicolumn{2}{c}{4.20} & \multicolumn{2}{c}{2.40} \\
\bottomrule
\end{tabular}}
\begin{tablenotes}[flushleft]
\footnotesize
\item \textit{Notes.} For DRO-K and RooD-SO, we use a regular grid design approach to choose \(N=4000\) discretization points.
\end{tablenotes}
\end{table}

As summarized in Table~\ref{tab:cost_coffee_5sets}, our proposed RooD-SO achieves the best average rank. These results suggest that the generalization capabilities and adjustable conservatism of the data-driven uncertainty set contribute to its favorable overall out-of-distribution performance. Although DRO-K achieves the lowest costs for Products 3 and 5, it ranks last for Products 1 and 4. This variability contrasts with the more consistent performance of RooD-SO across targets. Recall that DRO-W reduces to SAA under the newsvendor loss, thereby providing no robustness. Although DRO-Conv explicitly describes distributional heterogeneity by constructing a convex hull of all source distributions, it lacks the flexibility of regulating conservatism. We also remark that the results for DRO-Intx are not reported here because its optimization process always terminates prematurely due to an \texttt{Out-of-Memory} error.

\subsection{Portfolio optimization on real-world ETF data}
In this subsection, we consider a classical risk-averse portfolio optimization problem cast in a Conditional Value-at-Risk (CVaR) formulation \cite{rockafellar2000optimization, krokhmal2002portfolio}. The cost function is defined with an auxiliary variable \(\tau\in \mathbb{R}\) introduced:
\begin{equation}
    \begin{aligned}
        f((\mathbf{x},\tau), \boldsymbol{\xi}) &= -\langle \mathbf{x}, \boldsymbol{\xi}\rangle +\delta_1\left[\tau + \frac{1}{\delta_2}(-\langle \mathbf{x}, \boldsymbol{\xi}\rangle-\tau)^+\right].
    \end{aligned}
    \label{eq:portfolio}
\end{equation}
Since this is a convex piecewise-linear cost function, its convex continuity majorant \(\omega_r^f(\mathbf x)\) can be derived explicitly according to Table~\ref{tab:convex-continuity-majorants}. The decision variable \(\mathbf{x}\in\mathbb{R}^{d}\) pertains to the portfolio asset allocations, with its feasible set defined as \(\mathcal{X} = \{\mathbf{x}\in \mathbb{R}^{d}\,|\, \boldsymbol{1}^\top\mathbf{x}=1, \mathbf{x}\geq \boldsymbol{0}\}\), ensuring full investment and no short-selling. The uncertain asset returns are represented by \(\boldsymbol{\xi}\in\mathbb{R}^d\). We assume no access to historical return data \(\boldsymbol{\xi}\), but datasets \(\mathcal{D}_1,\cdots,\mathcal{D}_M\) from \(M\geq 2\) sources are available. 

We compare the methods for portfolio optimization using real-world Exchange-Traded Fund (ETF) data. The portfolios are constructed with exposure to three distinct sectors: Utilities, Information Technology, and Healthcare. This sectoral combination spans rate-sensitive defensive (Utilities), innovation-driven growth (Information Technology), and demand-stable defensive (Healthcare) sectors, thereby providing diverse exposures with distinct macroeconomic sensitivities. In this experiment, we construct nine sectoral ETF proxy sets for these three sectors, all denominated in U.S. dollars, with each set annotated by the geographic exposure of its underlying holdings and the exchange where it is listed. These nine ETF proxy sets feature heterogeneity in portfolio design choices, including factor tilts, weighting schemes, size exposures, selection methodologies, as well as geographic dimensions reflected in exposure region and listing venue. Historical price data are obtained from the Bloomberg Terminal for July 2020 to October 2025, and we compute monthly log-returns for each ETF over the sample period. The experiments are carried out in a leave-one-out manner similar to Section~\ref{sec:experiment_coffee_sales}.\endnote{Sectoral ETF data were collected through Bloomberg Terminal (institutional access).} Details of the data and computational setup are provided in
Appendix~\ref{app:portfolio_optimization_details}.

\begin{table}[t]
\centering
\caption{Average Sharpe Ratios across Nine Leave-One-Out Targets}
\label{tab:sharpe_loso_9sets}
\resizebox{\linewidth}{!}{%
\begin{tabular}{ccccccccccccccc}
\toprule
\multirow[b]{2}{*}{\raisebox{-0.6ex}{{\shortstack{Target\\[-1pt]Set}}}} & \multicolumn{2}{c}{SAA} & \multicolumn{2}{c}{DRO-W} & \multicolumn{2}{c}{SAA-WB} & \multicolumn{2}{c}{DRO-K} & \multicolumn{2}{c}{DRO-Conv} & \multicolumn{2}{c}{RooD-SO}\\
\cmidrule(lr){2-3} \cmidrule(lr){4-5} \cmidrule(lr){6-7} \cmidrule(lr){8-9} \cmidrule(lr){10-11} \cmidrule(lr){12-13}
& Avg. & Rank & Avg. & Rank & Avg. & Rank & Avg. (Std) & Rank  & Avg. & Rank & Avg. (Std) & Rank \\
\midrule
Set 1 & 0.250 & 4 & 0.265 & 3 & 0.237 & 5 & 0.284 (0.002) & 2 & 0.226 & 6 & \textbf{0.297 (0.005)} & 1 \\
Set 2 & 0.242 & 4 & 0.262 & 3 & 0.229 & 5 & 0.268 (0.002) & 2 & 0.222 & 6 & \textbf{0.282 (0.005)} & 1 \\
Set 3 & 0.265 & 4 & 0.280 & 3 & 0.264 & 5 & 0.291 (0.003) & 2 & 0.239 & 6 & \textbf{0.304 (0.007)} & 1 \\
Set 4 & 0.247 & 4 & 0.256 & 3 & 0.238 & 5 & 0.273 (0.002) & 2 & 0.228 & 6 & \textbf{0.286 (0.006)} & 1 \\
Set 5 & 0.063 & 5 & 0.074 & 4 & 0.062 & 6 & 0.076 (0.002) & 3 & \textbf{0.100} & 1 & 0.083 (0.005) & 2 \\
Set 6 & 0.224 & 5 & 0.221 & 6 & 0.231 & 4 & 0.234 (0.001) & 3 & \textbf{0.288} & 1 & 0.240 (0.003) & 2 \\
Set 7 & 0.146 & 4 & 0.159 & 3 & 0.141 & 5 & 0.195 (0.003) & 2 & 0.134 & 6 & \textbf{0.221 (0.010)} & 1 \\
Set 8 & 0.314 & 3 & 0.311 & 4 & 0.293 & 6 & 0.343 (0.003) & 2 & 0.298 & 5 & \textbf{0.361 (0.007)} & 1 \\
Set 9 & 0.259 & 3 & 0.259 & 3 & 0.255 & 5 & 0.290 (0.003) & 2 & 0.251 & 6 & \textbf{0.302 (0.007)} & 1 \\
\midrule
Avg. Rank & \multicolumn{2}{c}{4.00} & \multicolumn{2}{c}{3.56} & \multicolumn{2}{c}{5.11} & \multicolumn{2}{c}{2.22} & \multicolumn{2}{c}{4.78} & \multicolumn{2}{c}{1.22}\\
\bottomrule
\end{tabular}}
\begin{tablenotes}[flushleft]
\footnotesize
\item \textit{Notes.} For DRO-K and RooD-SO, Latin hypercube sampling is used for discretization. Sampling the supplementary discretization points introduces additional randomness; we therefore also report the standard deviation of the Sharpe ratio over 20 replications.
\end{tablenotes}
\end{table}

The computational results in Table~\ref{tab:sharpe_loso_9sets} showcase the excellent performance of RooD-SO, which achieves the highest Sharpe ratios in 7 out of 9 settings. Meanwhile, DRO-Conv demonstrates highly polarized performance. Although it secures the top rank in the remaining two cases, it consistently falls to the last place in the other seven. Such a highly variable performance stems from the inherent conservatism of the convex-hull construction of the uncertainty set. By contrast, RooD-SO inherits the powerful learning capability of SMC, ensuring effective out-of-distribution generalization across diverse scenarios. Furthermore, as an alternative capable of modulating conservatism, DRO-K stands out as the most competitive baseline. However, it underperforms RooD-SO, primarily due to the lack of an explicit description of distributional heterogeneity. In this case, DRO-Intx also terminates prematurely due to an \texttt{Out-of-Memory} error.

\section{Concluding remarks}
\label{sec:conclusion}
In this work, we proposed the RooD-SO framework to inform safeguarding decisions under a completely unobserved distribution, empowered by multiple relevant data sources. To encode distributional heterogeneity, we established a novel meta-distribution modeling philosophy, where both the unobserved target distribution and all observed source distributions are treated as independent realizations from an unknown meta-distribution. Grounded in this view, we presented a learning-theoretic uncertainty-set construction, which flexibly captures the structured uncertainty, and derived a worst-case formulation. The meta-distribution assumption entails two layers of randomness, arising from both the sampling of source distributions and the drawing of observations within each source. To account for both, we established out-of-distribution generalization guarantees, thereby providing theoretical justifications for the proposed framework. On the computational side, we developed a tractable constraint-tightening scheme and further introduced a Nystr\"om approximation-based acceleration strategy to improve scalability. We proved that these computational approximations preserve the theoretical guarantees of the original formulation and derived explicit upper bounds on their respective optimality losses. To assess the out-of-distribution robustness of our framework, we performed a systematic experimental evaluation. Numerical results demonstrate that, in practical scenarios where only a small or moderate number of data sources are available, our framework achieves superior out-of-distribution generalization and outperforms various data-driven decision-making schemes in unseen environments.

Several promising topics are worthy of future investigation. First, it would make sense to describe target-distribution uncertainty using higher-order or more informative statistical profiles of the meta-distribution. Besides, incorporating contextual information~\cite{sadana_survey_2024} into our considered setting seems to be a promising direction along two dimensions. By augmenting each source with its associated covariates, one can learn a covariate-to-distribution mapping for conditional estimation, enabling out-of-distribution decision-making using contextual information. Alternatively, our framework can be extended to cases where each source distribution is itself conditional on covariates.


\theendnotes

\bibliography{sn-bibliography}
\bibliographystyle{plain}

\newpage

\begin{appendices}

\section{Proofs}\label{Appendix:proofs}

\subsection{Additional notation}
Before proceeding, we introduce some additional notation used only in the appendix. The notation \(\|Z\|_{\infty}\) denotes the essential supremum norm of a random variable \(Z\), defined as \(\|Z\|_{\infty}:=\inf\{C\geq 0\,|\,{\rm Pr}(\|Z\|\leq C)=1\}\). We use \(\chi_{\mathcal{Y}}\) to denote the \(\{0,1\}\)-valued indicator function of a set \(\mathcal{Y}\). We denote by \({\rm Proj}_{\mathcal{H}}\) the orthogonal projector onto \(\mathcal{H}\). For any real numbers \(a, b\in\mathbb{R}\), we write \(a \wedge b:=\min\{a, b\}\). For $a \in \mathbb{R}$, we write $a_+:=\max\{a,0\}$; positive parts and minima of vectors are taken componentwise. For any bounded operator \(\Delta:\mathcal{H}_1\rightarrow\mathcal{H}_2\), where \(\mathcal{H}_1\) and \(\mathcal{H}_2\) are Hilbert spaces, we define its induced operator norm \(\|\Delta\|_{\rm op}:=\sup_{\|g\|_{\mathcal{H}_1}=1}\|\Delta(g)\|_{\mathcal{H}_2}\), and use \(\Delta^*\) to denote its adjoint operator. 



\subsection{Supporting lemmas}
We recall a known result regarding the upper bound on \(\mathbb{E}\{\|\hat{\mu}_{\mathbb{P}}-\mu_{\mathbb{P}}\|_{\mathcal{H}_k}^2\}\) as well as the \(\sqrt{n}\)-consistency of \(\hat{\mu}_{\mathbb{P}}\)~\cite{tolstikhin2017minimax}. 
\begin{lemma}[\(\sqrt{n}\)-consistency of \(\hat{\mu}_{\mathbb{P}}\),~\cite{tolstikhin2017minimax}]\label{lemma:empirical_kernel_mean_embedding_consistency}
    Let \(\{\boldsymbol{\xi}_j\}_{j=1}^n\) be random samples drawn i.i.d. from \(\mathbb{P}\in\mathcal{P}(\Xi)\). Then, the estimation error of kernel mean embedding admits the following second-moment bound:
    \begin{equation}\label{ineq:2_moment_bound_E_empirical}
        \left(\mathbb{E}\left\{\|\hat{\mu}_\mathbb{P}-\mu_{\mathbb{P}}\|_{\mathcal{H}_k}\right\}\right)^2\leq \mathbb{E}\left\{\|\hat{\mu}_{\mathbb{P}}-\mu_{\mathbb{P}}\|_{\mathcal{H}_k}^2\right\}\leq \frac{C}{n}.
    \end{equation}
    Moreover, for any \(\delta\in(0,1)\), the following inequality holds with probability at least \(1-\delta\):
    \begin{equation*}
        \|\hat{\mu}_{\mathbb{P}}-\mu_{\mathbb{P}}\|_{\mathcal{H}_k}\leq \sqrt{\frac{C}{n}} + \sqrt{\frac{2C\log(1/\delta)}{n}}.
    \end{equation*}
\end{lemma}

\begin{lemma}[Moment error bound for empirical mean embedding] \label{lemma:moment_bound_E_empirical}
    Let \(\{\boldsymbol{\xi}_j\}_{j=1}^n\) be random samples drawn i.i.d. from \(\mathbb{P}\in\mathcal{P}(\Xi)\). Then, for any \(q>0\), the following holds:
    \begin{equation}\label{ineq:q_moment_bound_E_empirical}
        \mathbb{E}\left\{\|\hat{\mu}_{\mathbb{P}}-\mu_{\mathbb{P}}\|_{\mathcal{H}_k}^q\right\}\leq\sqrt{4\pi q}\left(\frac{4Cq}{en}\right)^{q/2}e^{1/(6q)}.
    \end{equation}
\end{lemma}

\begin{proof}
    Since \(\sup_{\boldsymbol{\xi}\in \Xi}k(\boldsymbol{\xi},\boldsymbol{\xi})\leq C<\infty\), we have \(\|\phi(\boldsymbol{\xi})-\mu_{\mathbb{P}}\|_{\mathcal{H}_k}\leq \|\phi(\boldsymbol{\xi})\|_{\mathcal{H}_k}+\|\mu_{\mathbb{P}}\|_{\mathcal{H}_k}\leq 2\sqrt{C}\). Given \(\{\boldsymbol{\xi}_{j}\}_{j=1}^n\overset{\rm i.i.d.}{\sim}\mathbb{P}\), we define \(\psi_j:=\phi(\boldsymbol{\xi}_j)-\mu_{\mathbb{P}}\), which satisfies \(\mathbb{E}_{\mathbb{P}}\left\{\psi_j\right\}=0\) and \(\|\psi_j\|_{\mathcal{H}_k}\leq 2\sqrt{C}\). 

    By \cite[Theorem 3.5]{pinelis1994optimum}, which applies to martingales in \(2\)-smooth separable Banach spaces, we recall that any Hilbert space is 2-smooth with 2-smoothness constant \(D=1\), and that the partial sums of independent zero-mean random vectors form a martingale. Hence, \cite[Theorem 3.5]{pinelis1994optimum} is applicable to our setting. Define the martingale differences \(d_j:=\psi_j/n, j\in[n],\)
    and the associated martingale \(y_m:=\sum_{j=1}^{m}d_j, m\in[n]\).
    Let \(\mathcal{F}_j:=\sigma(\boldsymbol{\xi}_1,\cdots,\boldsymbol{\xi}_j)\) denote the \(\sigma\)-algebra generated by the first \(j\) samples. Then, \(\{d_j,\mathcal{F}_j\}\) is a martingale difference sequence in \(\mathcal{H}_k\). Moreover, since \(\|d_j\|_{\mathcal{H}_k}\leq 2\sqrt{C}/n\), the following holds:
    \begin{equation*}
        \sum_{j=1}^n\|d_j\|_{\infty}^2\leq n\cdot \frac{4C}{n^2}=\frac{4C}{n}=:b_*^2.
    \end{equation*}
    Invoking \cite[Theorem 3.5]{pinelis1994optimum} with \(D=1\) gives, for all \(t\geq 0\), 
    \begin{equation*}
        {\rm Pr}_{\mathbb{P}}\left\{\bar{y}\geq t\right\}\leq 2\exp\left\{-\frac{t^2}{2D^2b_*^2}\right\}=2\exp\left\{-\frac{nt^2}{8C}\right\},
    \end{equation*}
    where \(\bar{y}:=\sup_{m\leq n}\|y_m\|_{\mathcal{H}_k}\). Since \(\|(1/n)\sum_{j=1}^n\psi_j\|_{\mathcal{H}_k}=\|y_n\|_{\mathcal{H}_k}\leq \bar{y}\), we obtain
    \begin{equation*}
        \left\{\left\|\frac{1}{n}\sum_{j=1}^n\psi_j\right\|_{\mathcal{H}_k}\geq t\right\}\subseteq\left\{\bar{y}\geq t\right\},
    \end{equation*}
    and
    \begin{equation*}
        {\rm Pr}_{\mathbb{P}^n}\left\{\left\|\frac{1}{n}\sum_{j=1}^{n}\psi_j\right\|_{\mathcal{H}_k}\geq t\right\}\leq {\rm Pr}_{\mathbb{P}^n}\left\{\bar{y}\geq t\right\}\leq2\exp\left\{-\frac{nt^2}{8C}\right\}.
    \end{equation*}
    For \(q>0\), by the layer cake representation, we have the following tail-integration formula:
    \begin{equation*}
        \mathbb{E}\left\{\|y_n\|_{\mathcal{H}_k}^q\right\}=q\int_{0}^{\infty}t^{q-1}{\rm Pr}_{\mathbb{P}^n}\left\{y_n\geq t\right\}\mathrm{d}t\leq2q\int_{0}^{\infty}t^{q-1}\exp\left\{-\frac{nt^2}{8C}\right\}\mathrm{d}t=q\left(\frac{8C}{n}\right)^{q/2}\Gamma\left(\frac{q}{2}\right).
    \end{equation*}
    By a non-asymptotic Stirling bound for the Gamma function~\cite{gordon1994stochastic}, for any \(x>0\), \(\Gamma(x)\leq \sqrt{2\pi}x^{x-\frac{1}{2}}e^{-x+\frac{1}{12x}}\). Hence, with \(x=q/2\), \(q\Gamma\left(\frac{q}{2}\right)\leq \sqrt{4\pi q}\left(\frac{q}{2e}\right)^{q/2}e^{1/(6q)}\). Accordingly, we obtain
    \begin{equation*}
        \mathbb{E}\left\{\|y_n\|_{\mathcal{H}_k}^q\right\}\leq\sqrt{4\pi q}\left(\frac{4Cq}{en}\right)^{q/2}e^{1/(6q)},
    \end{equation*}
    which completes the proof.    
\end{proof}

Lemma~\ref{lemma:moment_bound_E_empirical} implies that \(\mathbb{E}\{\|\hat{\mu}_{\mathbb{P}}-\mu_{\mathbb{P}}\|_{\mathcal{H}_k}^q\}\lesssim n^{-q/2}\). Note that previous work~\cite{tolstikhin2017minimax} has considered the case with $q=2$ and established that \(\mathbb{E}\{\|\hat{\mu}_{\mathbb{P}}-\mu_{\mathbb{P}}\|_{\mathcal{H}_k}^2\} \le C/n\). By contrast, our result is more conservative in the case of \(q=2\), where \eqref{ineq:q_moment_bound_E_empirical} becomes \(\mathbb{E}\{\|\hat{\mu}_{\mathbb{P}}-\mu_{\mathbb{P}}\|_{\mathcal{H}_k}^2\}\leq 16.04C/n\). This is because our general \(q\)-moment inequality is technically built upon a uniform bound for the Gamma function. 

\subsection{Proof of Proposition~\ref{prop:discussion_dual}}
For part (a), given the unique dual optimizer \(\boldsymbol{\alpha}^*\), the KKT characterization yields \(\mu_c^*=\sum_{i=1}^M\alpha_i^*\mu_{\mathbb P_i}\). To determine an optimal radius, consider the two index sets
\[
    {\rm BSM}:=\{i\,|\,0<\alpha_i^*<1/M\nu\},\quad {\rm ESM}:=\{i\,|\,\alpha_i^*=1/M\nu\}.
\]
For any \(i\in{\rm BSM}\), complementary slackness implies \(s_i^*=0\), and thus \((R^*)^2=\|\mu_{\mathbb P_i}-\mu_c\|_{\mathcal H_k}^2\). Thus, if BSM is nonempty, the optimal radius is uniquely determined. However, BSM may be empty, in which case the primal problem~\eqref{eq:soft_minimum_ball} may admit non-unique optimal radii. To justify the radius selection, let \(d_i^*:=\|\mu_{\mathbb P_i}-\mu_c\|_{\mathcal H_k}^2\). The KKT conditions imply 
\[
    \max_{i:\alpha_i^*<1/M\nu}d_i^*\leq (R^*)^2\leq \min_{i:\alpha_i^*>0}d_i^*.
\]
Under the conditions in (a), \(\max_{i:\alpha_i^*<1/M\nu}d_i^*\) is well-defined and strictly positive. Set \((R^*)^2=\max_{i:\alpha_i^*<1/M\nu}d_i^*\) and \(s_i^*=\max\{0,d_i^*-\max_{i:\alpha_i^*<1/M\nu}d_i^*\}\). The primal feasibility still holds. Moreover, if \(\alpha_i^*<1/M\nu\), then \(s_i^*=0\); if \(\alpha_i^*>0\), then \(d_i^*=s_i^*+\max_{i:\alpha_i^*<1/M\nu}d_i^*\). Hence both complementary slackness conditions hold. Together with the center formula and dual feasibility, this verifies all KKT conditions and establishes optimality. This selection gives the smallest optimal radius. All measures indexed by ESM lie outside the ball, while all remaining measures lie on or inside it. This formula~\eqref{eq:R_square} also applies when BSM is nonempty, since in that case \(\max_{i:\alpha_i^*<1/M\nu}d_i^*=\min_{i:\alpha_i^*>0}d_i^* = (R^*)^2\). This proves part (a).

For part (b), where ~\eqref{eq:dual_SMC_quadratic} has multiple optimizers, we first prove the uniqueness of the optimal center. Let $K_{ij}=\langle\mu_{\mathbb P_i},\mu_{\mathbb P_j}\rangle_{\mathcal H_k}$, and denote the dual objective by
\[
    D(\boldsymbol\alpha):=-\boldsymbol\alpha^\top \mathbf K\boldsymbol\alpha +\operatorname{diag}(\mathbf K)^\top\boldsymbol\alpha.
\]
Take any $\boldsymbol\alpha,\widetilde{\boldsymbol\alpha}\in\mathcal A^*$, and let $D^*$ be the optimal dual value. Since the dual feasible set is convex, $(\boldsymbol\alpha+\widetilde{\boldsymbol\alpha})/2$ is feasible. Moreover,
\[
    D\!\left(\frac{\boldsymbol\alpha+\widetilde{\boldsymbol\alpha}}{2}\right) = \frac{D(\boldsymbol\alpha)+D(\widetilde{\boldsymbol\alpha})}{2}+\frac14(\boldsymbol\alpha-\widetilde{\boldsymbol\alpha})^\top \mathbf K(\boldsymbol\alpha-\widetilde{\boldsymbol\alpha}).
\]
By optimality of $\boldsymbol\alpha$ and ${\widetilde{\boldsymbol\alpha}}$, and positive semidefiniteness of $\mathbf K$, it follows that
\[
    D^* \ge D\!\left(\frac{\boldsymbol\alpha+\widetilde{\boldsymbol\alpha}}{2}\right) = D^*+\frac14(\boldsymbol\alpha-\widetilde{\boldsymbol\alpha})^\top \mathbf K(\boldsymbol\alpha-\widetilde{\boldsymbol\alpha}) \ge D^*.
\]
Consequently,
\[
    0 = (\boldsymbol\alpha-\widetilde{\boldsymbol\alpha})^\top \mathbf K(\boldsymbol\alpha-\widetilde{\boldsymbol\alpha})=\left\|\sum_{i=1}^M(\boldsymbol\alpha_i-\widetilde{\boldsymbol\alpha}_i)\mu_{\mathbb P_i}\right\|_{\mathcal H_k}^2.
\]
Thus $T\boldsymbol\alpha=T\widetilde{\boldsymbol\alpha}$, where $T\boldsymbol\alpha:=\sum_{i=1}^M\alpha_i\mu_{\mathbb P_i}$. By strong duality, the KKT identity $\mu_c^*=T\boldsymbol\alpha^*$ holds for any primal optimizer and any dual optimizer. Therefore, all optimal primal centers coincide, and every dual optimizer induces this same center. 

Next, we establish the uniqueness of the optimal radius \(R^*\). We first show that there exists a dual optimal solution with a nonempty BSM set. Suppose, for contradiction, that for every \(\boldsymbol{\alpha}\in\mathcal A^*\), the corresponding BSM set is empty. Then, every coordinate of every dual optimizer belongs to \(\{0,1/M\nu\}\). Since \(|\mathcal A^*|>1\), choose two distinct dual optimizers $\boldsymbol\alpha,\widetilde{\boldsymbol\alpha}\in\mathcal A^*$. There exists an index \(j\) such that \(\alpha_j\neq \tilde\alpha_j\), and hence \(\{\alpha_j,\tilde\alpha_j\}=\{0,1/M\nu\}\). By convexity of the feasible set and concavity of \(D\), \(\boldsymbol\alpha^\lambda \triangleq \lambda\boldsymbol\alpha+(1-\lambda)\widetilde{\boldsymbol\alpha}\) is also an optimal solution for any \(\lambda\in(0,1)\). However, \(0<\alpha_j^\lambda<1/M\nu\), thereby contradicting the assumption that every dual optimizer has an empty BSM set. Hence, there exist \(\boldsymbol{\alpha}'\in\mathcal A^*\) and an index \(i_{\rm R}\) such that \(0<\alpha'_{i_{\rm R}}<1/M\nu\). Fix this dual optimizer. For any optimal primal solution, the KKT conditions give \((R^*)^2=\|\mu_{\mathbb P_{i_{\rm R}}}-\mu_c\|_{\mathcal H_k}^2\). Since the optimal center is unique, the optimal radius is also unique. Moreover, for every \(i\) with \(\alpha'_i<1/M\nu\), the KKT conditions imply \(\|\mu_{\mathbb P_{i}}-\mu_c\|_{\mathcal H_k}^2\leq (R^*)^2\), with equality at \(i=i_{\rm R}\). Therefore, we can uniquely determine the optimal radius by~\eqref{eq:R_square} as well. This proves part (b).

\subsection{Proof of Proposition~\ref{prop:operator_norm_bound_Khat_K}}
\label{appendix:proof_thm_K_Khat_op}
Given \(M\) distributions \(\{\mathbb{P}_i\}_{i=1}^M\), their true kernel mean embeddings \(\{\mu_{\mathbb{P}_i}\}_{i=1}^M\) and the corresponding empirical estimates \(\{\hat{\mu}_{\mathbb{P}_i}\}\), we denote the estimation error as \(\Delta_{i}:=\hat{\mu}_{\mathbb{P}_i}-\mu_{\mathbb{P}_i}\) for all \(i\in[M]\). For \(\boldsymbol{v}=(v_1,\cdots,v_M)^\top\in\mathbb{R}^M\), we define \(\psi=(\mu_{\mathbb{P}_1},\dots,\mu_{\mathbb{P}_M}):\mathbb{R}^M\rightarrow \mathcal{H}_k\) with \(\psi(\boldsymbol{v})=\sum_{i=1}^Mv_i\mu_{\mathbb{P}_i}\), and \(\Delta=(\Delta_1,\dots,\Delta_M):\mathbb{R}^M\rightarrow\mathcal{H}_k\) with \(\Delta(\boldsymbol{v})=\sum_{i=1}^Mv_i\Delta_i\).  Then, the population and empirical Gram matrices can be written as \(\mathbf{K}=\boldsymbol{\psi}^* \boldsymbol{\psi}\) and \(\hat{\mathbf{K}}=(\boldsymbol{\psi}+\Delta)^*(\boldsymbol{\psi}+\Delta)\).
Algebraically, \(\hat{\mathbf{K}}-\mathbf{K}=(\boldsymbol{\psi}+\Delta)^*(\boldsymbol{\psi}+\Delta)-\boldsymbol{\psi}^* \boldsymbol{\psi}
= \boldsymbol{\psi}^*\Delta+\Delta^*\boldsymbol{\psi} + \Delta^*\Delta\), and thus
\begin{equation} \label{eq:decompose_Khat_K}
    \|\hat{\mathbf{K}}-\mathbf{K}\|_2 = \|\boldsymbol{\psi}^*\Delta+\Delta^*\boldsymbol{\psi} + \Delta^*\Delta\|_2\leq 2\|\boldsymbol{\psi}\|_{\rm op}\|\Delta\|_{\rm op}+\|\Delta\|_{\rm op}^2.
\end{equation}
To derive an upper bound for \(\|\hat{\mathbf{K}}-\mathbf{K}\|_2\), we first bound \(\|\Delta\|_{\rm op}\):
\begin{equation*}
    \begin{aligned}
        \|\Delta\|_{\rm op}=\sup_{\|v\|_2=1}\left\|\sum_{i=1}^Mv_i\Delta_i\right\|_{\mathcal{H}_k}\leq\sup_{\|v\|_2=1}\sum_{i=1}^M|v_i|\|\Delta_i\|_{\mathcal{H}_k}\leq\left(\sum_{i=1}^M \|\Delta_i\|_{\mathcal{H}_k}^2\right)^{1/2}.
    \end{aligned}
\end{equation*}
Then, we bound each \(\|\Delta_i\|_{\mathcal{H}_k}\).
By Lemma~\ref{lemma:empirical_kernel_mean_embedding_consistency},
for any \(\delta_i\in(0,1)\),
\begin{equation*}
    {\rm Pr}_{\mathbb{P}_i^{n_i}}\left\{\|\Delta_i\|_{\mathcal{H}_k}
\leq \sqrt{\frac{C}{n_i}}+\sqrt{\frac{2C\log(1/\delta_i)}{n_i}}\right\}\ge 1-\delta_i.
\end{equation*}
By the Bonferroni inequality, setting \(\delta_i=\delta/M\) and taking a union bound gives
\begin{equation}\label{eq:Delta_i_2}
    {\rm Pr}_{\mathbb{P}_1^{n_1}\otimes\cdots\otimes\mathbb{P}_M^{n_M}}\left\{\sum_{i=1}^M \|\Delta_i\|_{\mathcal{H}_k}^2
\le
2C \sum_{i=1}^M\!\Big(\frac{1}{n_i}+\frac{2\log(M/\delta)}{n_i}\Big)\right\}\geq 1-\delta.
\end{equation}
Next, we establish an upper bound for \(\|\boldsymbol{\psi}\|_{\text{op}}\). Since \(\mathbf{K}=\boldsymbol{\psi}^*\boldsymbol{\psi}\), we have
\begin{equation} \label{eq:psi_op_2_equal_K_op}
    \|\boldsymbol{\psi}\|_{\text{op}}^2 =\lambda_{\max}(\boldsymbol{\psi}^*\boldsymbol{\psi})=\lambda_{\max}(\mathbf{K}) = \|\mathbf{K}\|_2.
\end{equation}
By the inequality \(\lambda_{\max}(\mathbf{K})\leq {\rm tr}(\mathbf K)=\sum_{i=1}^M\|\mu_{\mathbb{P}_i}\|_{\mathcal{H}_k}^2\), we obtain
\begin{equation*}
    \|\boldsymbol{\psi}\|_{\rm op}\leq \left(\sum_{i=1}^M\|\mu_{\mathbb{P}_i}\|_{\mathcal{H}_k}^2\right)^{1/2}\leq \sqrt{M}\max_i\|\mu_{\mathbb{P}_i}\|_{\mathcal{H}_k}.
\end{equation*}
Moreover, since \(\|\mu_{\mathbb{P}_i}\|_{\mathcal{H}_k}^2=\mathbb{E}_{\boldsymbol{\xi},\boldsymbol{\zeta}\sim\mathbb{P}_i}\{ k(\boldsymbol{\xi},\boldsymbol{\zeta}) \} \leq C\), we have
\begin{equation} \label{eq:psi_op}
    \|\boldsymbol{\psi}\|_{\rm op}\leq \sqrt{MC}.
\end{equation}
Plugging \eqref{eq:Delta_i_2} and \eqref{eq:psi_op} into \eqref{eq:decompose_Khat_K} yields~\eqref{eq:main_Khat_K_norm_bound}. Thus, the result follows.

\subsection{Proof of Proposition~\ref{prop:consistency_optimizer}}
\label{appendix:proof_prop_consistency_optimizer}
Let
\[
    \mathcal A_\nu:=\left\{\boldsymbol\alpha\in\mathbb R^M: 0\leq\alpha_i\leq\frac{1}{M\nu},\ i\in[M],\quad\mathbf 1^\top\boldsymbol\alpha=1\right\}.
\]
The population dual problem can be rewritten as
\begin{equation}
    \min_{\boldsymbol\alpha\in\mathcal A_\nu}F(\boldsymbol\alpha),\qquad F(\boldsymbol\alpha):=\boldsymbol\alpha^\top\mathbf K\boldsymbol\alpha-\operatorname{diag}(\mathbf K)^\top\boldsymbol\alpha.
\end{equation}
Define \(F^*:=\min_{\boldsymbol\alpha\in \mathcal A_\nu}F(\boldsymbol\alpha)\). We have \(\mathcal A^*=\arg\min_{\boldsymbol\alpha\in\mathcal A_\nu}F(\boldsymbol\alpha)\). Since $\mathcal A_\nu$ is nonempty and compact and $F$ is continuous, $\mathcal A^*$ is nonempty. Fix $\boldsymbol\alpha_0\in\mathcal A^*$ and define \(\mathbf L:=\mathbf K^{1/2}\), \(g_0:=\nabla F(\boldsymbol\alpha_0)\). The first-order optimality condition gives
\[
    s(\boldsymbol\alpha):=g_0^\top(\boldsymbol\alpha-\boldsymbol\alpha_0)\ge0, \quad \forall\boldsymbol\alpha\in \mathcal A_\nu.
\]
Moreover, the quadratic structure of $F$ implies
\[
    F(\boldsymbol\alpha)-F^* = \|\mathbf L(\boldsymbol\alpha-\boldsymbol\alpha_0)\|_2^2+s(\alpha).
\]
Consequently,
\[
    \mathcal A^* =\left\{\boldsymbol\alpha\in \mathcal A_\nu: \mathbf L(\boldsymbol\alpha-\boldsymbol\alpha_0)=0,\ g_0^\top(\boldsymbol\alpha-\boldsymbol\alpha_0)=0 \right\}.
\]
By the Hoffman's error bound~\cite{beck2017linearly}, there exists a finite constant $H>0$ such that
\[
\operatorname{dist}(\boldsymbol\alpha,\mathcal A^*)^2\le H^2\left( \|\mathbf L(\boldsymbol\alpha-\boldsymbol\alpha_0)\|_2^2+s(\boldsymbol\alpha)^2 \right),\qquad \forall\boldsymbol\alpha\in \mathcal A_\nu.
\]

Since $\mathcal A_\nu$ is compact, $B:=\max_{\boldsymbol\alpha\in \mathcal A_\nu}s(\boldsymbol\alpha)$ is finite and nonnegative. Thus $s(\boldsymbol\alpha)^2\le B s(\boldsymbol\alpha)$, and hence
\[
    \operatorname{dist}(\boldsymbol\alpha,\mathcal A^*)^2 \le H^2\max\{1,B\}\bigl(F(\boldsymbol\alpha)-F^*\bigr).
\]
Setting $\rho:=2/(H^2\max\{1,B\})>0$ yields 
\begin{equation} \label{eq:Falpha_minus_F_rho}
    F(\boldsymbol{\alpha})-F^*\geq \frac{\rho}{2}\operatorname{dist}(\boldsymbol{\alpha},\mathcal A^*)^2,\quad \forall \boldsymbol{\alpha}\in\mathcal A_\nu.
\end{equation}

In the finite-sample regime, the empirical dual problem is cast as
\begin{equation}
    \min_{\boldsymbol\alpha\in\mathcal A_\nu}\hat F(\boldsymbol\alpha),\qquad \hat F(\boldsymbol\alpha):=\boldsymbol\alpha^\top\hat{\mathbf K}\boldsymbol\alpha-\operatorname{diag}(\hat{\mathbf K})^\top\boldsymbol\alpha.
\end{equation}
Let \(\hat F^*:=\min_{\boldsymbol\alpha\in \mathcal A_\nu}\hat F(\boldsymbol\alpha)\) and \(\hat{\mathcal A}^*=\arg\min_{\boldsymbol\alpha\in\mathcal A_\nu}\hat F(\boldsymbol\alpha)\). For notational simplicity, let \(\mathbf q:=\operatorname{diag}(\mathbf K)\) and \(\hat{\mathbf q}:=\operatorname{diag}(\hat{\mathbf K})\). We consider \(\hat{\boldsymbol{\alpha}}^*\in\hat{\mathcal A}^*\), choose \(\boldsymbol{\alpha}^*=\operatorname{Proj}_{\mathcal A^*}(\hat{\boldsymbol{\alpha}}^*)\), and define \(\mathbf d:=\hat{\boldsymbol\alpha}^*-\boldsymbol\alpha^*\). Then we have \(\|\mathbf d\|_2=\operatorname{dist}(\hat{\boldsymbol{\alpha}}^*,\mathcal A^*)\). Due to the convexity of \(F(\cdot)\), the first-order optimality condition for the empirical problem, and~\eqref{eq:Falpha_minus_F_rho}, we have
\begin{align*}
    \frac{\rho}{2}\|\mathbf d\|_2^2\leq F(\hat{\boldsymbol{\alpha}}^*)-F^*
    &\leq\left\langle \nabla F(\hat{\boldsymbol\alpha}^*), \mathbf d \right\rangle \leq \left\langle \nabla F(\hat{\boldsymbol\alpha}^*) - \nabla\hat F(\hat{\boldsymbol\alpha}^*), \mathbf d \right\rangle = \left\langle 2(\mathbf K-\hat{\mathbf K}) \hat{\boldsymbol\alpha}^* + \hat{\mathbf q}-\mathbf q, \mathbf d \right\rangle\\
    &\leq \left( 2\|\hat{\mathbf K}-\mathbf K\|_2 \|\hat{\boldsymbol\alpha}^*\|_2 + \|\hat{\mathbf q}-\mathbf q\|_2 \right) \|\mathbf d\|_2.
\end{align*}
If \(\mathbf d=\mathbf 0\), the desired result is immediate. Otherwise, dividing both sides by \(\rho\|\mathbf d\|_2/2\) yields
\begin{equation}
\label{eq:restricted_qp_perturbation}
    \operatorname{dist}(\hat{\boldsymbol\alpha}^*,\mathcal A^*)\leq \frac{4}{\rho}\left( \|\hat{\mathbf K}-\mathbf K\|_2 \|\hat{\boldsymbol\alpha}^*\|_2 + \frac12\|\hat{\mathbf q}-\mathbf q\|_2 \right).
\end{equation}
For any feasible \(\boldsymbol\alpha\), \(\|\boldsymbol\alpha\|_2^2 \leq \|\boldsymbol\alpha\|_\infty  \|\boldsymbol\alpha\|_1 \leq 1\wedge\frac{1}{M\nu}\). Define \(C_\alpha^2:=1\wedge\frac{1}{M\nu}\), and we have
\begin{equation} \label{eq:alpha_feasible_norm}
    \|\hat{\boldsymbol{\alpha}}^*\|_2\leq C_\alpha.
\end{equation}
Applying Proposition~\ref{prop:operator_norm_bound_Khat_K}
with confidence level \(\delta/2\), we obtain
\begin{equation}
\label{eq:K_concentration_for_alpha}
    \|\hat{\mathbf K}-\mathbf K\|_2\leq  \sqrt{\frac{8MC^2(1+2\log(2M/\delta))}{\bar{n}}}+\frac{2C(1+2\log(2M/\delta))}{\bar{n}},
\end{equation}
with probability at least \(1-\delta/2\). Moreover, applying the concentration inequality for empirical kernel mean embeddings with confidence level \(\delta/(2M)\) and taking a union bound gives
\[
    \|\hat\mu_{\mathbb P_i}-\mu_{\mathbb P_i}\|_{\mathcal H_k} \leq \sqrt{\frac{C}{n_i}} + \sqrt{\frac{2C\log(2M/\delta)}{n_i}}, \qquad i\in[M],
\]
with probability at least \(1-\delta/2\). On this event,
\begin{align}
    \frac12\|\hat{\mathbf q}-\mathbf q\|_2&\leq
    \frac{C}{\sqrt{\bar n}}\left(1+\sqrt{2\log(2M/\delta)}\right),
\label{eq:q_concentration_for_alpha}
\end{align}
where the first inequality follows from \(|\hat q_i-q_i|\leq 
\left|\|\hat\mu_{\mathbb P_i}\|_{\mathcal H_k}^2- \|\mu_{\mathbb P_i}\|_{\mathcal H_k}^2\right|\leq 2\sqrt{C}\|\hat\mu_{\mathbb P_i}-\mu_{\mathbb P_i}\|_{\mathcal H_k}\). Finally, combining~\eqref{eq:restricted_qp_perturbation},~\eqref{eq:alpha_feasible_norm},~\eqref{eq:K_concentration_for_alpha}, and~\eqref{eq:q_concentration_for_alpha}, and taking a union bound gives
\[
\left\| \hat{\boldsymbol\alpha}^{*} - \boldsymbol\alpha^{*} \right\|_2 \leq \frac{4C}{\rho}
\left[ \left( C_\alpha\sqrt{8M}+\sqrt{2} \right) \sqrt{\frac{1+2\log(2M/\delta)}{\bar n}} + \frac{2C_\alpha (1+2\log(2M/\delta))}{\bar n} \right]
\]
with probability at least \(1-\delta\), which completes the proof.

\subsection{Proof of Theorem~\ref{thm:finite_sample_guar_mu_c_R_square}}
\label{appendix:proof_thm_finite_sample_guar_mu_c_R_square}

Let \(\Delta_i:=\hat\mu_{\mathbb P_i}-\mu_{\mathbb P_i}\) for \(i\in[M]\), and choose \(\boldsymbol\alpha^*:=\operatorname{Proj}_{\mathcal A^*}(\hat{\boldsymbol\alpha}^*).\) By Proposition~\ref{prop:consistency_optimizer}, for fixed \(M\), \(\nu\), and kernel \(k\), \(\|\hat{\boldsymbol\alpha}^*-\boldsymbol\alpha^*\|_2 = \operatorname{dist} (\hat{\boldsymbol\alpha}^*,\mathcal A^*)  = \mathcal O_{\mathbb P}\!\left(\bar n^{-1/2}\right)\). By Proposition~\ref{prop:discussion_dual}, this choice satisfies \(\mu_c^* = \sum_{i=1}^M\alpha_i^*\mu_{\mathbb P_i}\). Moreover, the concentration inequality for empirical kernel mean embeddings yields $\|\Delta_i\|_{\mathcal H_k}=\mathcal{O}_{\mathbb P}(n_i^{-1/2})$ for each $i\in[M]$. Since $M$ is fixed, it follows that
\[
\left(\sum_{i=1}^M\|\Delta_i\|_{\mathcal H_k}^2\right)^{1/2}=\mathcal{O}_{\mathbb P}\!\left(\bar n^{-1/2}\right).
\]
Using $\mu_c^*=\sum_{i=1}^M\alpha_i^*\mu_{\mathbb P_i}$ and $\hat\mu_c^*=\sum_{i=1}^M\hat\alpha_i^*\hat\mu_{\mathbb P_i}$, we obtain
\begin{align}
\|\hat\mu_c^*-\mu_c^*\|_{\mathcal H_k}
&\leq \left\|\sum_{i=1}^M(\hat\alpha_i^*-\alpha_i^*)\mu_{\mathbb P_i}\right\|_{\mathcal H_k}+\left\|\sum_{i=1}^M\hat\alpha_i^*\Delta_i\right\|_{\mathcal H_k} \notag\\
&\leq \sqrt{\lambda_{\max}(\mathbf K)}\|\hat{\boldsymbol\alpha}^*-\boldsymbol\alpha^*\|_2+\|\hat{\boldsymbol\alpha}^*\|_2\left(\sum_{i=1}^M\|\Delta_i\|_{\mathcal H_k}^2\right)^{1/2}.
\end{align}
Since $\hat{\boldsymbol\alpha}^*\in\mathcal A_\nu$, $\|\hat{\boldsymbol\alpha}^*\|_2\leq\|\hat{\boldsymbol\alpha}^*\|_1=1$. Therefore,
\[
\|\hat\mu_c^*-\mu_c^*\|_{\mathcal H_k}=\mathcal{O}_{\mathbb P}\!\left(\bar n^{-1/2}\right).
\]

We next consider the estimated radius. For $\mathbf d=(d_1,\ldots,d_M)\in\mathbb R_+^M$, define
\[
Q_\nu(\mathbf d):=\min\operatorname*{arg\,min}_{t\geq0}\left\{t+\frac{1}{M\nu}\sum_{i=1}^M(d_i-t)_+\right\}.
\]
For a fixed center, eliminating the slack variables from the primal SMC problem shows that the squared radius is given by this order-statistic map. Consequently, \((R^*)^2=Q_\nu((\|\mu_{\mathbb P_i}-\mu_c^*\|_{\mathcal H_k}^2)_{i=1}^M)\), and analogously for $(\hat R^*)^2$. Since $Q_\nu$ is an order-statistic map, \(|Q_\nu(\mathbf a)-Q_\nu(\mathbf b)|\leq\|\mathbf a-\mathbf b\|_\infty\). Define $d_i:=\|\mu_{\mathbb P_i}-\mu_c^*\|_{\mathcal H_k}^2$ and $\hat d_i:=\|\hat\mu_{\mathbb P_i}-\hat\mu_c^*\|_{\mathcal H_k}^2$. It follows that
\[
|(\hat R^*)^2-(R^*)^2|\leq\max_{i\in[M]}|\hat d_i-d_i|.
\]
By Assumption~\ref{ass:kernel}, $\|\mu_{\mathbb P_i}\|_{\mathcal H_k}$, $\|\hat\mu_{\mathbb P_i}\|_{\mathcal H_k}$, $\|\mu_c^*\|_{\mathcal H_k}$, and $\|\hat\mu_c^*\|_{\mathcal H_k}$ are all bounded above by $\sqrt C$. Hence, for every $i\in[M]$,
\begin{align}
|\hat d_i-d_i|
&\leq \left(\|\hat\mu_{\mathbb P_i}-\hat\mu_c^*\|_{\mathcal H_k}+\|\mu_{\mathbb P_i}-\mu_c^*\|_{\mathcal H_k}\right)\left\|(\hat\mu_{\mathbb P_i}-\hat\mu_c^*)-(\mu_{\mathbb P_i}-\mu_c^*)\right\|_{\mathcal H_k} \notag\\
&\leq 4\sqrt C\left(\|\Delta_i\|_{\mathcal H_k}+\|\hat\mu_c^*-\mu_c^*\|_{\mathcal H_k}\right).
\end{align}
Therefore,
\[
|(\hat R^*)^2-(R^*)^2|\leq4\sqrt C\left[\max_{i\in[M]}\|\Delta_i\|_{\mathcal H_k}+\|\hat\mu_c^*-\mu_c^*\|_{\mathcal H_k}\right].
\]
Since $M$ is fixed,
\[
\max_{i\in[M]}\|\Delta_i\|_{\mathcal H_k}\leq\left(\sum_{i=1}^M\|\Delta_i\|_{\mathcal H_k}^2\right)^{1/2}=\mathcal{O}_{\mathbb P}\!\left(\bar n^{-1/2}\right).
\]
Combining this with the center bound gives
\[
|(\hat R^*)^2-(R^*)^2|=\mathcal{O}_{\mathbb P}\!\left(\bar n^{-1/2}\right).
\]
This completes the proof.

\subsection{Proof of Theorem~\ref{thm:SMC_generalization_error_bound}}
\label{appendix:proof_thm_SMC_generalization_error_bound}
    Under Assumption~\ref{ass:kernel}, we have \(\|\mu_\mathbb{P}\|_{\mathcal{H}_k}^2=\left\|\mathbb{E}_{\mathbb{P}}\left\{\phi(\boldsymbol{\xi})\right\}\right\|_{\mathcal{H}_k}^2\leq \mathbb{E}_{\mathbb{P}}\left\{\|\phi(\boldsymbol{\xi})\|_{\mathcal{H}_k}^2\right\}=\mathbb{E}_{\mathbb{P}}\left\{k(\boldsymbol{\xi},\boldsymbol{\xi})\right\}\le C\). Before proceeding, we introduce a decision function \(h(\mu)\) to encode membership in \(\mathcal{B}_{\nu,\gamma}\). In SMC theory~\cite{muandet2013one}, this decision function used for testing whether the uniquely determined by a new mean embedding \(\mu\) belongs to an RKHS-ball with center \(c\) and radius \(r\) takes the following form:
    \begin{equation*}
        h(\mu)=\text{sgn}\left(r^2-\|\mu-c\|_{\mathcal H_k}^2\right) = (\text{sgn} \circ g_{c,r})(\mu),
    \end{equation*}
    where \(g_{c,r}(\mu)=r^2-\|\mu-c\|_{\mathcal H_k}^2\). To derive an out-of-distribution generalization error bound for \(\mathcal{U}_{\nu,\gamma}^{\rm SMC}\), we begin by introducing the population-level risk function, for any \(\gamma > 0\):
    \begin{equation}\label{eq:risk_function_Rg}
        \mathcal{R}(g)={\rm Pr}(g(\mu)\leq-\gamma) = \mathbb{E}\{\chi_{\mathcal{Y}}(g(\mu))\},
    \end{equation}
    where \(\mathcal{Y}:=\{y\in\mathbb{R}:y\leq-\gamma\}\). The empirical counterpart of \(\mathcal{R}(g)\) over \(\{\mu_{\mathbb{P}_i}\}_{i=1}^M\) is defined as:
    \begin{equation*} 
        \hat{\mathcal{R}}(g)=\frac{1}{M}\sum_{i=1}^{M}\chi_{\mathcal{Y}}(g(\mu_{i})),
    \end{equation*}
   where \(\mu_i:=\mu_{\mathbb{P}_i}\) for notational simplicity. Let \(\boldsymbol{\mu}:=(\mu_1,\cdots, \mu_{M})\), and \(\hat{g}^*\) be the optimal decision function  by solving problem~\eqref{eq:dual_SMC_quadratic} to optimality. The statement of the theorem is equivalent to establishing an upper bound on \(\mathcal{R}(\hat{g}^*)\). Since the indicator cost function \(\chi_{\mathcal{Y}}(\cdot)\) is inherently discontinuous, we first replace it with a Lipschitz-continuous surrogate, specifically a truncated hinge-type loss defined as:
    \begin{equation} \label{eq:truncated_loss_function}
        \ell_{\gamma}(z):=1\wedge\left(-\frac{z}{\gamma}\right)^+.
    \end{equation}
    This function is \(\frac{1}{\gamma}\)-Lipschitz and serves as an upper bound on the indicator function: \(\chi_{\mathcal{Y}}(z)\leq \ell_\gamma(z)\). Based on this surrogate, we define the associated surrogate risk as follows:
    \begin{equation*}
        \mathcal{R}_s(g)=\mathbb{E}\left\{\ell_\gamma(g(\mu))\right\},\quad \hat{\mathcal{R}}_s(g)=\frac{1}{M}\sum_{i=1}^M\ell_\gamma(g(\mu_i)).
    \end{equation*}
    By Lemma~\ref{lemma:nu}, at most a fraction of \(\nu\) of the source embeddings lie strictly outside the learned ball. Consequently, 
    \begin{equation} \label{ineq:apply_Lemma_nu}
        \hat{\mathcal{R}}_s(\hat g^*)=\frac{1}{M}\sum_{i=1}^M\ell_\gamma(\hat g^*(\mu_i))\leq \frac{1}{M}\sum_{i=1}^M\boldsymbol{1}\{\hat g^*(\mu_i)<0\}\leq \nu.
    \end{equation}
    The deviation of the empirical risk minimizer from its population counterpart can be controlled uniformly over the function class \(\mathcal{G}\) via the following inequality:
    \begin{equation*}
        \begin{aligned}
            \mathcal{R}_s(\hat{g}^*)-\hat{\mathcal{R}}_s(\hat{g}^*) &\leq \sup_{g\in\mathcal{G}}\left\{ \mathcal{R}_s(g)-\hat{\mathcal{R}}_s(g) \right\}.\\
        \end{aligned}
    \end{equation*}
    Then, by applying McDiarmid's inequality~\cite{mcdiarmid1989method}, with the probability at least \(1-\delta\), we have
    \begin{equation*}
        \sup_{g\in\mathcal{G}}\left\{ \mathcal{R}_s(g)-\hat{\mathcal{R}}_s(g) \right\}\leq \mathbb{E}_{\boldsymbol{\mu}}\left\{\sup_{g\in\mathcal{G}}\left( \mathcal{R}_s(g)-\hat{\mathcal{R}}_s(g) \right)\right\}+\sqrt{\frac{\log(1/\delta)}{2M}}.
    \end{equation*}
    To prove the theorem, it suffices to show that, with probability at least \(1-\delta\), the following holds:
    \begin{equation}\label{ineq:R_s_final_inequality}
        \mathcal{R}(\hat{g}^*)\leq \mathcal{R}_s(\hat{g}^*)\leq \hat{\mathcal{R}}_s(\hat{g}^*)+\mathbb{E}_{\boldsymbol{\mu}}\left\{\sup_{g\in\mathcal{G}}\left( \mathcal{R}_s(g)-\hat{\mathcal{R}}_s(g) \right)\right\}+\sqrt{\frac{\log(1/\delta)}{2M}}.
    \end{equation}
    Henceforth, we turn to upper bound the term \(\mathbb{E}_{\boldsymbol{\mu}}\{\sup_{g\in\mathcal{G}}( \mathcal{R}_s(g)-\hat{\mathcal{R}}_s(g))\}\). According to the Vapnik-Chervonenkis symmetrization lemma \cite{vapnik2006estimation} and the Ledoux-Talagrand contraction inequality \cite{ledoux2013probability}, we have
    \begin{equation}\label{ineq:E_Rademacher}
        \begin{aligned}
            \mathbb{E}_{\boldsymbol{\mu}}\left\{\sup_{g\in\mathcal{G}}\left( \mathcal{R}_s(g)-\hat{\mathcal{R}}_s(g) \right)\right\}&\leq 2\,\mathbb{E}_{\boldsymbol{\mu},\boldsymbol{\sigma}}\left\{\sup_{g\in\mathcal{G}}\frac{1}{M}\sum_{i=1}^M\sigma_i(\ell_{\gamma}\circ g)(\mu_i)\right\}\\
            &\leq \frac{2}{\gamma}\,\mathbb{E}_{\boldsymbol{\mu},\boldsymbol{\sigma}}\left\{\sup_{g\in\mathcal{G}}\frac{1}{M}\sum_{i=1}^M\sigma_ig(\mu_i)\right\}= \frac{2}{\gamma}\,\mathfrak{R}_M(\mathcal{G}),
        \end{aligned}
    \end{equation}
    where \(\mathfrak{R}_M(\mathcal{G})\) is the Rademacher average of \(\mathcal{G}\) with the form:
    \begin{equation*}
        \mathfrak{R}_M(\mathcal{G}):=\mathbb{E}_{\boldsymbol{\mu}}\left\{\mathbb{E}_{\boldsymbol{\sigma}}\left\{\sup_{g\in \mathcal{G}}\frac{1}{M}\sum_{i=1}^M\sigma_ig(\mu_i)\,\middle|\, \boldsymbol{\mu}\right\}\right\}.
    \end{equation*}
    Since the learned center is a convex combination of the source embeddings, its norm is at most $\sqrt C$. Moreover, the feasible choice $c=0$, $r^2=C$, and $s_i=0$ has objective value $C$ in~\eqref{eq:soft_minimum_ball}, so every optimal squared radius is at most $C$. Hence, the learned decision function belongs to the following class:
    \begin{equation} \label{eq:function_class_G}
        \mathcal{G} = \{g_{c,r}(\mu)=r^2-\|\mu-c\|_{\mathcal H_k}^2: \|c\|_{\mathcal{H}_k}\leq \sqrt{C}, 0\leq r^2\leq C\}.
    \end{equation} 
    Next, we compute the upper bound of \(\mathfrak{R}_{M}(\mathcal{G})\):
    \begin{equation*}
        \begin{aligned}
            \mathfrak{R}_{M}(\mathcal{G}) &= \mathbb{E}_{\boldsymbol{\sigma},\boldsymbol{\mu}}\left\{\sup_{\|c\|_{\mathcal{H}_k}\leq \sqrt{C}, 0\leq r^2\leq C} \frac{1}{M}\sum_{i=1}^M\sigma_i(r^2-\|\mu_i-c\|_{\mathcal H_k}^2)\right\}\\
            &\leq\mathbb{E}_{\boldsymbol{\sigma},\boldsymbol{\mu}}\left\{\sup_{\|c\|_{\mathcal H_k}\leq \sqrt{C}}\left(\frac{1}{M}\sum_{i=1}^M2\sigma_i\langle \mu_i,c\rangle_{\mathcal H_k}\right)\right\}+\frac{C}{M}\cdot\mathbb{E}_{\boldsymbol{\sigma}}\left\{\left|\sum_{i=1}^M\sigma_i\right|\right\}\\
            &\leq 2\sqrt C\cdot \mathbb{E}_{\boldsymbol{\sigma},\boldsymbol{\mu}}\left\{\left\|\frac{1}{M}\sum_{i=1}^M\sigma_i\mu_i\right\|_{\mathcal{H}_k}\right\}+ \frac{C}{M}\cdot\mathbb{E}_{\boldsymbol{\sigma}}\left\{\left|\sum_{i=1}^M\sigma_i\right|\right\},\\
        \end{aligned}
    \end{equation*}
    where the first inequality follows from 
    \[
        \mathbb E_{\boldsymbol \sigma}\!\left\{-\frac1M\sum_{i=1}^M\sigma_i\|\mu_i\|_{\mathcal H_k}^2\right\}=0,
    \]
    and 
    \[
        \frac{1}{M}\sum_{i=1}^M\sigma_i(r^2-\|c\|_{\mathcal H_k}^2)\leq \frac{C}{M}\left|\sum_{i=1}^M\sigma_i\right|.
    \]
    The second inequality follows from 
    \[
        2\left\langle \frac1M\sum_{i=1}^M\sigma_i\mu_i,c\right\rangle_{\mathcal H_k}\le 2\sqrt C\,\left\|\frac1M\sum_{i=1}^M\sigma_i\mu_i\right\|_{\mathcal H_k}.
    \]
    To further bound the expected norm term above, we obtain
    \begin{equation}\label{ineq:mu_nodata_1}
        \begin{aligned}
            \mathbb{E}_{\boldsymbol{\sigma},\boldsymbol{\mu}}\left\{\left\|\frac{1}{M}\sum_{i=1}^M\sigma_i\mu_i\right\|_{\mathcal{H}_k}\right\} 
            \leq \frac{1}{M}\left(\mathbb{E}_{\boldsymbol{\sigma},\boldsymbol{\mu}}\left\{\left\|\sum_{i=1}^M\sigma_i\mu_i\right\|_{\mathcal{H}_k}^2\right\}\right)^{1/2}=\frac{1}{M}\left(\mathbb{E}_{\boldsymbol{\mu}}\left\{\mathbb{E}_{\boldsymbol{\sigma}}\left\{\left\|\sum_{i=1}^M\sigma_i\mu_i\right\|_{\mathcal{H}_k}^2\,\middle|\,\boldsymbol{\mu}\right\}\right\}\right)^{1/2}.
        \end{aligned}
    \end{equation}
    With \(\sigma_i\) i.i.d. Rademacher and independent of \(\boldsymbol{\mu}\), it follows that 
    \begin{equation}\label{ineq:mu_nodata_2}
        \mathbb{E}_{\boldsymbol{\sigma}}\left\{\left\|\sum_{i=1}^M\sigma_i\mu_i\right\|^2\right\}=\mathbb{E}_{\boldsymbol{\sigma}}\left\{\left\langle\sum_{i=1}^M\sigma_i\mu_i,\sum_{i=1}^M\sigma_i\mu_i\right\rangle_{\mathcal{H}_k}\right\}=\sum_{i=1}^M\sum_{j=1}^M\mathbb{E}_{\boldsymbol{\sigma}}\left\{\sigma_i\sigma_j\right\}\langle\mu_i,\mu_j\rangle_{\mathcal{H}_k}=\sum_{i=1}^M\|\mu_i\|^2_{\mathcal{H}_k}.
    \end{equation}
   Accordingly, we derive
    \begin{equation}
            \mathfrak{R}_{M}(\mathcal{G}) 
            \leq \frac{2\sqrt C}{M}\sqrt{\sum_{i=1}^M\mathbb{E}_{\boldsymbol{\mu}}\left\{\|\mu_i\|_{\mathcal{H}_k}^2\right\}}+\frac{C}{M}\sqrt{\mathbb{E}_{\sigma}\left\{ \left(\sum_{i=1}^M\sigma_i\right)^2 \right\}}\leq \frac{3C}{\sqrt{M}}.
        \label{ineq:SMC_RC_upperbound2}
    \end{equation}
    Combining~\eqref{ineq:apply_Lemma_nu}, \eqref{ineq:E_Rademacher} and~\eqref{ineq:SMC_RC_upperbound2} with~\eqref{ineq:R_s_final_inequality}, it follows that
    \begin{equation*}
        \mathcal{R}(\hat{g}^{*})\leq \nu+\frac{6C}{\gamma\sqrt{M}}+\sqrt{\frac{\log(1/\delta)}{2M}}.
    \end{equation*}
    Finally, we obtain the claimed result of the theorem.


\subsection{Proof of Proposition~\ref{prop:empirical_generalization_error_bound}}
\label{appendix:proof_prop_empirical_generalization_error_bound}
    Given the stated settings, we write \(\mu_{i}\) for the kernel mean embedding of each true distribution \(\mathbb{P}_i\) and \(\hat{\mu}_{i}\) for the corresponding empirical embedding computed based on \(\mathcal{D}_i\). We further define \(\boldsymbol{\mu}:=(\mu_1,\cdots, \mu_{M})\) and \(\hat{\boldsymbol{\mu}}:=(\hat{\mu}_1,\cdots, \hat{\mu}_{M})\) for subsequent analysis. 

    To prove this proposition, we follow a similar strategy to that of Theorem~\ref{thm:SMC_generalization_error_bound}; the proof of the theorem appears in Appendix~\ref{appendix:proof_thm_SMC_generalization_error_bound}. We define \(\mathcal{Y}:=\{y\in\mathbb{R}:y\leq-\gamma\}\) and the population-level risk function \(\mathcal{R}(g)\) as in~\eqref{eq:risk_function_Rg}, whereas here the empirical counterpart of \(\mathcal{R}(g)\) over \(\hat{\boldsymbol{\mu}}\) is defined as:
    \begin{equation*}
        \hat{\mathcal{R}}(g)=\frac{1}{M}\sum_{i=1}^{M}\chi_{\mathcal{Y}}(g(\hat{\mu}_{i})).
    \end{equation*}
    Let \(\hat{g}^*\) be an optimal decision function obtained by solving problem~\eqref{eq:dual_SMC_quadratic} with \(\boldsymbol{\mu}\) replaced by \(\hat{\boldsymbol{\mu}}\). The statement of this proposition is thus equivalent to establishing an upper bound on \(\mathcal{R}(\hat{g}^*)\). To this end, we use the same truncated hinge-type cost function \(\ell_\gamma(z)\) for any \(\gamma>0\), as defined in~\eqref{eq:truncated_loss_function}. Hence, the corresponding surrogate risk is defined as:
    \begin{equation*}
        \mathcal{R}_s(g)=\mathbb{E}_{\mathbb P}\left\{\ell_\gamma(g(\mu))\right\},\quad \bar{\mathcal R}_s(g)=\frac{1}{M}\sum_{i=1}^M\mathbb{E}_{(\mathbb P_i,\mathcal D_i)}\left\{\ell_\gamma(g(\hat\mu_i))\right\},\quad \hat{\mathcal{R}}_s(g)=\frac{1}{M}\sum_{i=1}^M\ell_\gamma(g(\hat{\mu}_i)).
    \end{equation*}
    The difference between the population and empirical surrogate risks of the learned decision function can be controlled uniformly over $\mathcal G$ as follows:
    \begin{equation*}
        \begin{aligned}
            \mathcal{R}_s(\hat{g}^*)-\hat{\mathcal{R}}_s(\hat{g}^*) \leq \sup_{g\in\mathcal{G}}\left|\mathcal{R}_s(g)-\bar{\mathcal{R}}_s(g)\right|+\sup_{g\in\mathcal{G}}\left\{\bar{\mathcal{R}}_s(g)-\hat{\mathcal{R}}_s(g)\right\},
        \end{aligned}
    \end{equation*}
    where the function class \(\mathcal{G}\) is defined in \eqref{eq:function_class_G}. To bound the first term, we first note that every \(g\in\mathcal G\) is \(4\sqrt C\)-Lipschitz on the closed ball of radius \(\sqrt C\) in \(\mathcal H_k\). Indeed, for any \(\mu,\nu\in\mathcal H_k\) satisfying \(\|\mu\|_{\mathcal H_k}\le\sqrt C\) and \(\|\nu\|_{\mathcal H_k}\le\sqrt C\), the definition of \(\mathcal G\) implies \(\|c\|_{\mathcal H_k}\le\sqrt C\), and hence \(|g(\mu)-g(\nu)|\leq |\langle \mu-\nu,\mu+\nu-2c \rangle_{\mathcal H_k}|\leq 4\sqrt C\|\mu-\nu\|_{\mathcal H_k}\). Thus, for every \(g\in\mathcal G\), we have
    \begin{equation*}
        \left| \mathcal R_s(g)-\bar{\mathcal R}_s(g)\right|=\frac1M \sum_{i=1}^M\mathbb E_{\mathbb P_i}\left\{\ell_\gamma(g(\mu_i))-\mathbb{E}_{\mathcal D_i}\{\ell_{\gamma}(g(\hat\mu_i))\}\right\}\leq \frac{4\sqrt C}{\gamma M}\sum_{i=1}^M\mathbb{E}\|\mu_i-\hat\mu_i\|_{\mathcal H_k}.
    \end{equation*}
    Moreover, by Lemma~\ref{lemma:empirical_kernel_mean_embedding_consistency}, \(\mathbb E \|\mu_i-\hat\mu_i\|_{\mathcal H_k}\leq \sqrt{C/n_i}\) gives 
    \begin{equation}\label{ineq: first_term}
    \sup_{g\in\mathcal G}\left|\mathcal R_s(g)-\bar{\mathcal R}_s(g)\right|\leq \frac{4 C}{\gamma M}\sum_{i=1}^M\frac{1}{\sqrt{n_i}}.
    \end{equation}

    As for the second term, it measures the average empirical generalization gap on the \(M\) observed sources. Note that the two-layer data-generating processes are independent across \(i\). The distributions \(\{\mathbb{P}_i\}\) are independently generated, and the samples are then drawn independently from their respective distributions. Therefore, \(\hat{\mu}_1,\cdots,\hat{\mu}_M\) are independent random variables, and we can still apply McDiarmid's inequality~\cite{mcdiarmid1989method}. Then, with the probability at least \(1-\delta\), we have
    \begin{equation*}
        \sup_{g\in\mathcal{G}}\left\{ \bar{\mathcal{R}}_s(g)-\hat{\mathcal{R}}_s(g) \right\}\leq \mathbb{E}_{\boldsymbol{\hat{\mu}}}\left\{\sup_{g\in\mathcal{G}}\left[ \bar{\mathcal{R}}_s(g)-\hat{\mathcal{R}}_s(g) \right]\right\}+\sqrt{\frac{\log(1/\delta)}{2M}}.
    \end{equation*}
    Henceforth, we turn to upper bound the term \(\mathbb{E}_{\hat{\boldsymbol{\mu}}}\left\{\sup_{g\in\mathcal{G}}\left| \bar{\mathcal{R}}_s(g)-\hat{\mathcal{R}}_s(g) \right|\right\}\). According to the Vapnik-Chervonenkis symmetrization lemma \cite{vapnik2006estimation} and the Ledoux-Talagrand contraction inequality \cite{ledoux2013probability}, we have
    \begin{equation}\label{ineq:E_Rademacher_empirical}
        \begin{aligned}
            \mathbb{E}_{\hat{\boldsymbol{\mu}}}\left\{\sup_{g\in\mathcal{G}}\left[ \bar{\mathcal{R}}_s(g)-\hat{\mathcal{R}}_s(g) \right]\right\}&\leq 2\,\mathbb{E}_{\hat{\boldsymbol{\mu}},\boldsymbol{\sigma}}\left\{\sup_{g\in\mathcal{G}}\frac{1}{M}\sum_{i=1}^M\sigma_i(\ell_{\gamma}\circ g)(\hat{\mu}_i)\right\}\\
            &\leq \frac{2}{\gamma}\,\mathbb{E}_{\hat{\boldsymbol{\mu}},\boldsymbol{\sigma}}\left\{\sup_{g\in\mathcal{G}}\frac{1}{M}\sum_{i=1}^M\sigma_ig(\hat{\mu}_i)\right\},
        \end{aligned}
    \end{equation}
    where
    \begin{equation}
        \begin{aligned}
            &\mathbb{E}_{\hat{\boldsymbol{\mu}},\boldsymbol{\sigma}}\left\{\sup_{g\in\mathcal{G}}\frac{1}{M}\sum_{i=1}^M\sigma_ig(\hat{\mu}_i)\right\}= \mathbb{E}_{\hat{\boldsymbol{\mu}},\boldsymbol{\sigma}}\left\{\sup_{\|c\|_{\mathcal{H}_k}\leq \sqrt{C}, 0\leq r^2\leq C} \frac{1}{M}\sum_{i=1}^M\sigma_i(r^2-\|\hat{\mu}_i-c\|_{\mathcal H_k}^2)\right\}\\
            \leq &~ 2\sqrt C\cdot\mathbb{E}_{\hat{\boldsymbol{\mu}},\boldsymbol{\sigma}}\left\{\left\|\frac{1}{M}\sum_{i=1}^M\sigma_i\hat{\mu}_{i}\right\|_{\mathcal{H}_k}\right\}+\frac{C}{\sqrt{M}}\leq \frac{3C}{\sqrt M}.\label{ineq:SMC_RC_upperbound2_empirical}\\
        \end{aligned}
    \end{equation}  
    Combining~\eqref{ineq:apply_Lemma_nu}, \eqref{ineq: first_term},~\eqref{ineq:E_Rademacher_empirical} and  \eqref{ineq:SMC_RC_upperbound2_empirical} with~\eqref{ineq:R_s_final_inequality} yields the claim of the proposition.

\subsection{Proof of Proposition~\ref{prop:asymptotic_behavior_U_empirical}}
Let \(\delta_M:=1/M\) and \(b_M:=\frac{1}{M}\sum_{i=1}^M\frac{1}{\sqrt{n_{i,M}}}+\frac{1}{\sqrt{M}}\), and recall that \(\gamma_M:=2\sqrt{b_M}\). Since \(b_M\to0\), we have \(\gamma_M\to0\), and hence \(\gamma_M\) lies in the admissible range for all sufficiently large \(M\). Applying Proposition~\ref{prop:empirical_generalization_error_bound} with \(\gamma=\gamma_M\) and \(\delta=\delta_M\), we obtain 
\[
    {\rm Pr}_{\Psi^M\otimes\mathbb{P}_1^{n_1}\otimes\cdots\otimes\mathbb{P}_{M}^{n_{M}}}\left\{{\rm Pr}_{\mathbb{P}'\sim\Psi}\left\{\mathbb{P}'\notin \widehat{\mathcal{U}}_{\nu,\gamma_M}^{\rm {SMC}}\right\}\leq\nu +r_M\right\}\geq 1-\delta_M,
\]
where
\[
    r_M:=\frac{4 C}{\gamma_M M}\sum_{i=1}^M\frac{1}{\sqrt{n_{i,M}}} +\frac{6C}{\gamma_M\sqrt M}+\sqrt{\frac{\log(1/\delta_M)}{2M}}.
\]
Next, we show that \(r_M\to 0\). Note that
\[
    0\leq r_M\leq \frac{6Cb_M}{\gamma_M}+\sqrt{\frac{\log M}{2M}}=3C\sqrt{b_M}+\sqrt{\frac{\log M}{2M}}\to 0.
\]

Given \(\epsilon>0\), there exists \(M_{\epsilon}\) such that \(r_M\leq \epsilon\) for every \(M\geq M_\epsilon\). On the event \({\rm Pr}_{\mathbb P'\sim\Psi}\{\mathbb P'\notin\widehat{\mathcal U}_{\nu,\gamma_M}^{\rm SMC}\}\leq \nu+r_M\), we therefore have
\[
    {\rm Pr}_{\mathbb P'\sim\Psi}\left\{\mathbb P'\in\widehat{\mathcal U}_{\nu,\gamma_M}^{\rm SMC}\right\}\geq 1-\nu-r_M\geq 1-\nu-\epsilon.
\]
It follows that for every \(M\geq M_\epsilon\),
\[
{\rm Pr}_{\Psi^M\otimes\mathbb P_1^{n_{1,M}}\otimes\cdots\otimes\mathbb P_M^{n_{M,M}}}\left\{\Pr_{\mathbb P'\sim\Psi}\left(\mathbb P'\in\widehat{\mathcal U}_{\nu,\gamma_M}^{\rm SMC}\right)\ge1-\nu-\epsilon\right\}\ge1-\frac1M.
\]
Letting \(M\to\infty\) yields the desired result.

\subsection{Proof of Proposition~\ref{prop:nonasymptotic_out_of_distribution_guarantee_DRO_finite_sample}}
\label{appendix:proof_prop_nonasymptotic_out_of_distribution_guarantee_DRO_finite_sample}
Proposition~\ref{prop:empirical_generalization_error_bound} has already established the nonasymptotic out-of-distribution generalization guarantee of \(\widehat{\mathcal{U}}_{\nu,\gamma}^{\rm SMC}\). If the event \(A:=\{\mathbb{P}'\sim\Psi\,|\,\mathbb{P}'\in\widehat{\mathcal{U}}_{\nu,\gamma}^{\rm SMC}\}\) occurs, then we have for any \(\mathbf{x}\in\mathcal{X}\), 
\begin{equation}
    \mathbb{E}_{\mathbb{P}'}\{f(\mathbf{x},\boldsymbol{\xi})\}\leq \sup_{\mathbb{P}\in\widehat{\mathcal{U}}_{\nu,\gamma}^{\rm SMC}}\mathbb{E}_{\mathbb{P}}\{f(\mathbf{x},\boldsymbol{\xi})\}.
\end{equation}
Therefore, 
\begin{equation}
    A\subseteq\left\{\mathbb{P}'\sim\Psi\,\middle|\,\forall\, \mathbf{x}\in\mathcal{X}: \mathbb{E}_{\mathbb{P}'}\{f(\mathbf{x},\boldsymbol{\xi})\}\leq \sup_{\mathbb{P}\in\widehat{\mathcal{U}}_{\nu,\gamma}^{\rm SMC}}\mathbb{E}_{\mathbb{P}}\{f(\mathbf{x},\boldsymbol{\xi})\}\right\}.
\end{equation}
Hence, with probability \(1-\delta\) over \((\mathbb{P}_1,\cdots, \mathbb{P}_M)\sim\Psi^{M}\) and \(\mathcal{D}_i\sim\mathbb{P}_i^{n_i}\) for \(i\in[M]\), we obtain
\[
    \begin{aligned}
    {\rm Pr}_{\mathbb{P}'\sim\Psi}\Biggl\{\forall\,\mathbf{x}\in\mathcal{X}: \mathbb{E}_{\mathbb{P}'}\{f(\mathbf{x},\boldsymbol{\xi})\}&\leq \sup_{\mathbb{P}\in\widehat{\mathcal{U}}_{\nu,\gamma}^{\rm SMC}}\mathbb{E}_{\mathbb{P}}\{f(\mathbf{x},\boldsymbol{\xi})\}\Biggr\}\\&\geq 1-\nu-\frac{4 C}{\gamma M}\sum_{i=1}^M\frac{1}{\sqrt{n_i}} -\frac{6C}{\gamma\sqrt M}-\sqrt{\frac{\log(1/\delta)}{2M}},
\end{aligned}
\]
which completes the proof.

\subsection{Proof of Proposition~\ref{prop:out_of_sample_guarantee_DRO_finite_sample}}
\label{appendix:proof_prop_out_of_sample_guarantee_DRO_finite_sample}
For each realization of the source distributions and the corresponding within-source samples, we define two events
\[
    A_M:=\left\{\mathbb{P}'\sim\Psi:\mathbb{P}'\in\widehat{\mathcal{U}}_{\nu,\gamma_M}^{\rm SMC}\right\},
\]
and 
\[
    B_M:=\left\{\mathbb{P}'\sim\Psi:\forall \mathbf x\in\mathcal X,\mathbb E_{\mathbb P'}\{f(\mathbf x,\boldsymbol{\xi})\}\leq \sup_{\mathbb P\in\widehat{\mathcal U}_{\nu,\gamma_M}^{\rm SMC}}\mathbb E_{\mathbb P}\{f(\mathbf x,\boldsymbol{\xi})\}\right\}.
\]
For any \(\mathbb P'\in A_M\), we have \(\mathbb{P}'\in\widehat{\mathcal{U}}_{\nu,\gamma_M}^{\rm SMC}\). Hence, for every \(\mathbf x\in\mathcal X\), \(\mathbb E_{\mathbb P'}\{f(\mathbf x,\boldsymbol{\xi})\}\leq \sup_{\mathbb P\in\widehat{\mathcal U}_{\nu,\gamma_M}^{\rm SMC}}\mathbb E_{\mathbb P}\{f(\mathbf x,\boldsymbol{\xi})\}\). Therefore, 
\[
    A_M\subseteq B_M.
\]
Consequently, conditionally on the observed source distributions and samples, we have
\[
    {\rm Pr}_{\mathbb P'\sim\Psi}\{B_M\}\geq {\rm Pr}_{\mathbb P'\sim\Psi}\{A_M\}.
\]
It follows that
\[
    \begin{aligned}
        \Bigl\{{\rm Pr}_{\mathbb P'\sim\Psi}&\bigl\{\mathbb{P}'\in\widehat{\mathcal{U}}_{\nu,\gamma_M}^{\rm SMC}\bigr\}\geq 1-\nu-\epsilon\Bigr\}\\
        &\subseteq\Biggl\{{\rm Pr}_{\mathbb P'\sim\Psi}\Bigl\{\forall \mathbf x\in\mathcal X:\mathbb E_{\mathbb P'}\{f(\mathbf x,\boldsymbol{\xi})\}\leq \sup_{\mathbb P\in\widehat{\mathcal U}_{\nu,\gamma_M}^{\rm SMC}}\mathbb E_{\mathbb P}\{f(\mathbf x,\boldsymbol{\xi})\}\Bigr\}\geq 1-\nu-\epsilon\Biggr\}.
    \end{aligned}
\]
Thus, 
\[
    \begin{aligned}
    {\rm Pr}_{R_M}\Biggl\{{\rm Pr}_{\mathbb P'\sim\Psi}\Bigl\{\forall \mathbf x\in\mathcal X:&\mathbb E_{\mathbb P'}\{f(\mathbf x,\boldsymbol{\xi})\}\leq \sup_{\mathbb P\in\widehat{\mathcal U}_{\nu,\gamma_M}^{\rm SMC}}\mathbb E_{\mathbb P}\{f(\mathbf x,\boldsymbol{\xi})\}\Bigr\}\geq 1-\nu-\epsilon\Biggr\}\\
    &\geq 
        {\rm Pr}_{R_M}\Bigl\{{\rm Pr}_{\mathbb P'\sim\Psi}\bigl\{\mathbb{P}'\in\widehat{\mathcal{U}}_{\nu,\gamma_M}^{\rm SMC}\bigr\}\geq 1-\nu-\epsilon\Bigr\}.\\
    \end{aligned}
\]
By Proposition~\ref{prop:asymptotic_behavior_U_empirical}, the right-hand side converges to one as \(M\to\infty\). Since the left-hand side is bounded above by one, the desired conclusion follows.


\subsection{Proof of Proposition~\ref{prop:dual_reformulation_SMC}}
\label{appendix:proof_prop_dual_reformulation_SMC}
    Let \(\hat{\mathcal{B}}_{\nu,\gamma} = \{\mu\in\mathcal{H}_k\mid\|\mu-\hat{\mu}_c^*\|_{\mathcal{H}_k}\leq \sqrt{(\hat{R}^*)^2+\gamma}\}\), where 
    \begin{equation}
    \hat{\mu}_c^*:=\sum_{i=1}^M\hat{\alpha}_i^*\hat{\mu}_{\mathbb{P}_i}, \quad (\hat{R}^*)^2=\max_{\alpha_i'<1/M\nu}\left(\hat{K}_{i'i'}-2\sum_{i=1}^M\hat{\alpha}^*_i\hat{K}_{ii'}+\sum_{i_1=1}^M\sum_{i_2=1}^M\hat{\alpha}^*_{i_1}\hat{\alpha}^*_{i_2}\hat{K}_{i_1i_2}\right).
    \end{equation} 
    We compute the corresponding support function as follows:  
    \begin{equation}
        \begin{aligned}
            h_{\hat{\mathcal{B}}_{\nu,\gamma}}^*(g) = \sup_{\mu\in \hat{\mathcal{B}}_{\nu,\gamma}}\langle g,\mu \rangle_{\mathcal{H}_k} 
             &= \langle g,\hat{\mu}_c \rangle_{\mathcal{H}_k} + \sup_{\|\mu-\hat{\mu}_c\|_{\mathcal{H}_k}\leq \hat{R}^*_\gamma}\langle g,\mu-\hat{\mu}_c \rangle\\
            &=\langle g,\hat{\mu}_c \rangle_{\mathcal{H}_k}+\hat{R}^*_\gamma\|g\|_{\mathcal{H}_k}\\
            &=\sum_{i\in{\rm SM}} \hat{\alpha}_i^*\mathbb{E}_{\boldsymbol{\xi}\sim\hat{\mathbb{P}}_i}\{g(\boldsymbol{\xi})\} + \hat{R}^*_\gamma\|g\|_{\mathcal{H}_k}\\
            &=\sum_{i\in{\rm SM}}\sum_{j=1}^{n_i}\frac{\hat{\alpha}_i^*}{n_i}g(\boldsymbol{\xi}_j^{(i)})+\hat{R}^*_\gamma\|g\|_{\mathcal{H}_k},
        \end{aligned}
    \end{equation}
    where the third equality follows from the Cauchy-Schwarz inequality. Therefore, by Lemma~\ref{lemma:generalized_duality}, problem~\eqref{eq:Kernel_SMM_DRO_Primary_finite} can be equivalently reformulated as problem~\eqref{eq:SMC_infinite_dual_formulation}.

\subsection{Proof of Theorem~\ref{thm:tightening_and_prior_bound}}
We first show that every feasible solution of problem~\eqref{eq:tightened_problem} is feasible for the original problem. Let \((\mathbf x, g,\beta)\) be feasible for problem~\eqref{eq:tightened_problem}, and fix an arbitrary \(\xi\in\Xi\). By the definition of \(r_N\), there exists an index \(l(\boldsymbol\xi)=\operatorname*{arg\,min}_{l\in[N]}\|\boldsymbol{\xi}-\boldsymbol{\zeta}_l\|\) such that \(\|\boldsymbol\xi-\boldsymbol\zeta_{l(\boldsymbol\xi)}\|\le r_N\). Since \(\omega_f\) majorizes the modulus of continuity of \(f\), \(f(\mathbf x, \boldsymbol{\xi})\leq f(\mathbf x,\boldsymbol{\zeta}_{l(\boldsymbol{\xi})})+\omega^f_{r_N}(\mathbf x)\). Moreover, by the reproducing property, \(|g(\boldsymbol{\xi})-g(\boldsymbol{\zeta}_{l(\boldsymbol{\xi})})|\leq d_k(\boldsymbol{\xi},\boldsymbol{\zeta}_{l(\boldsymbol{\xi})}) \|g\|_{\mathcal H_k}\leq \omega^k_{r_N}\|g\|_{\mathcal H_k}\). Feasibility for the constraint-tightened problem~\eqref{eq:tightened_problem} implies
\[
    f(\mathbf x, \boldsymbol{\zeta}_{l(\boldsymbol{\xi})})+\omega^f_{r_N}(\mathbf x)+\omega^k_{r_N}\|g\|_{\mathcal H_k}\leq g(\boldsymbol{\zeta}_{l(\boldsymbol{\xi})})+\beta.
\]
As a consequence,
\begin{equation} \label{ineq:bound_hold_f_g_tighten}
    f(\mathbf x, \boldsymbol{\xi})\leq f(\mathbf x, \boldsymbol{\zeta}_{l(\boldsymbol{\xi})})+\omega^f_{r_N}(\mathbf x)\leq g(\boldsymbol{\zeta}_{l(\boldsymbol{\xi})})+\beta-\omega^k_{r_N}\|g\|_{\mathcal H_k}\leq g(\boldsymbol{\xi})+\beta.
\end{equation}
Since~\eqref{ineq:bound_hold_f_g_tighten} holds for every \(\boldsymbol{\xi}\in\Xi\), \((\mathbf x,g,\beta)\) is feasible for the original problem~\eqref{eq:SMC_infinite_dual_formulation}. Therefore, the tightened feasible set is contained in the original feasible set, and hence \(z_{\rm D}\leq z_{\rm UB}\).

We next establish the optimality gap bound. Fix an optimal decision \(\mathbf x_{\rm D}^*\) of the problem~\eqref{eq:SMC_infinite_dual_formulation}. Since \(\hat\alpha^*_i/n_i\ge0\), \(\sum_{i\in{\rm SM}}\sum_{j=1}^{n_i}\hat\alpha^*_i/n_i=1\), for every \(g\in\mathcal H_k\),
\begin{align*}
\sup_{\boldsymbol\xi\in\Xi} \left(f(\mathbf x_{\rm D}^*,\boldsymbol\xi)-g(\boldsymbol\xi)\right)
&\ge \sum_{i\in{\rm SM}}\sum_{j=1}^{n_i} \frac{\hat\alpha^*_i}{n_i} \bigl( f(\mathbf x_{\rm D}^*,\boldsymbol\xi_j^{(i)}) -g(\boldsymbol\xi_j^{(i)}) \bigr).
\end{align*}
For this fixed decision, eliminating \(\beta\) gives the objective functional 
\[
   V(g)=\sum_{i\in {\rm SM}}\sum_{j=1}^{n_i}\frac{\hat{\alpha}_i^*}{n_i}g(\boldsymbol{\xi}_j^{(i)}) + \hat{R}^*_\gamma\|g\|_{\mathcal{H}_k} + \sup_{\boldsymbol{\xi}\in\Xi}(f(\mathbf x_{\rm D}^*,\boldsymbol{\xi})- g(\boldsymbol{\xi}))\geq \sum_{i\in {\rm SM}}\sum_{j=1}^{n_i}\frac{\hat{\alpha}_i^*}{n_i}f(\mathbf x_{\rm D}^*,\boldsymbol{\xi}_j^{(i)}) + \hat{R}^*_\gamma\|g\|_{\mathcal{H}_k}.
\]
Hence, since \(\hat{R}^*_\gamma>0\), the functional \(V\) is coercive on \(\mathcal H_k\). Under the standing compactness and continuity assumptions, \(V\) is finite-valued and proper. Furthermore, point evaluations are weakly continuous, the supremum term is weakly lower semicontinuous, and the RKHS norm is weakly lower semicontinuous. Thus, \(V\) is weakly lower semicontinuous. Since \(\mathcal H_k\) is a Hilbert space, \(V\) attains its minimum at some \(g^*\). Set \(\beta^*:=\sup_{\boldsymbol\xi\in\Xi}\bigl(f(\mathbf x_{\rm D}^*,\boldsymbol\xi)-g^*(\boldsymbol\xi)\bigr)\). Then \((\mathbf x_{\rm D}^*,g^*,\beta^*)\) is an optimal solution of problem~\eqref{eq:SMC_infinite_dual_formulation}, and \(V(g^*)=z_{\rm D}\).

The feasible comparison point \((\mathbf x_{\rm D}^*,0,\sup_{\boldsymbol\xi\in\Xi}f(\mathbf x_{\rm D}^*,\boldsymbol\xi))\) gives \(V(g^*)\leq V(0)=\sup_{\boldsymbol\xi\in\Xi}f(\mathbf x_{\rm D}^*,\boldsymbol\xi)\). Combining this inequality with the preceding lower bound yields 
\[ 
    \sum_{i\in{\rm SM}}\sum_{j=1}^{n_i}\frac{\hat\alpha_i^*}{n_i}f(\mathbf x_{\rm D}^*,\boldsymbol\xi_j^{(i)})+\hat R^*_\gamma\|g^*\|_{\mathcal H_k}\leq V(g^*)\leq\sup_{\boldsymbol\xi\in\Xi}f(\mathbf x_{\rm D}^*,\boldsymbol\xi). 
\] 
Therefore, 
\[
    \|g^*\|_{\mathcal H_k}\leq\frac{\sup_{\boldsymbol\xi\in\Xi}f(\mathbf x_{\rm D}^*,\boldsymbol\xi)-\sum_{i\in{\rm SM}}\sum_{j=1}^{n_i}\frac{\hat\alpha_i^*}{n_i}f(\mathbf x_{\rm D}^*,\boldsymbol\xi_j^{(i)})}{\hat R^*_\gamma} \leq\frac{\sup_{\boldsymbol\xi\in\Xi}f(\mathbf x_{\rm D}^*,\boldsymbol\xi)-\inf_{\boldsymbol\xi\in\Xi}f(\mathbf x_{\rm D}^*,\boldsymbol\xi)}{\hat R^*_\gamma}=\bar B_\gamma.
\]

Original feasibility gives \(f(\mathbf x_{\rm D}^*,\boldsymbol\zeta_l)\leq g^*(\boldsymbol\zeta_l)+\beta^*\) for every \(l\in[N]\). Therefore,
\[
    f(\mathbf x_{\rm D}^*, \boldsymbol\zeta_l) +\omega^f_{r_N}(\mathbf x_{\rm D}^*) +\omega^k_{r_N}\|g^*\|_{\mathcal H_k} \le g^*(\boldsymbol\zeta_l)+\beta^*+\omega^f_{r_N}(\mathbf x_{\rm D}^*) +\bar B_\gamma\omega^k_{r_N}.
\]
Thus, the triple \((\mathbf x_{\rm D}^*,g^*,\beta^*+\omega^f_{r_N}(\mathbf x_{\rm D}^*) +\bar B_\gamma\omega^k_{r_N})\) is feasible for problem~\eqref{eq:tightened_problem}. Its objective value gives
\[
    z_{\rm UB}\leq \sum_{i\in{\rm SM}}\sum_{j=1}^{n_i}\frac{\hat{\alpha}_i^*}{n_i}g^*(\boldsymbol{\xi}_j^{(i)}) + \hat{R}^*_\gamma\|g^*\|_{\mathcal{H}_k} + \beta^*+\omega^f_{r_N}(\mathbf x_{\rm D}^*) +\bar B_\gamma\omega^k_{r_N}=z_{\rm D}+\omega^f_{r_N}(\mathbf x_{\rm D}^*) +\bar B_\gamma\omega^k_{r_N}.
\]
which proves part (b).

Finally, consider any sequence of finite discretization sets \(\{\Xi_N\}\) satisfying \(r_N\to0\). Continuity of \(k\) on the compact set \(\Xi\times\Xi\) implies uniform continuity of its canonical feature map. Hence, \(\omega^k_{r_N}\to0\). By the assumed vanishing property of the selected majorant, \(\omega^f_{r_N}(\mathbf x_{\rm D}^*)\to0\). Moreover, \(\bar B_\gamma\) is finite and independent of \(N\). Combining the previous two parts, we obtain 
\[ 
    0\leq z_{\rm UB}-z_{\rm D}\leq\omega^f_{r_N}(\mathbf x_{\rm D}^*)+\bar B_\gamma\omega^k_{r_N}\rightarrow0. 
\] 
The conclusion in part~(c) follows from the squeeze theorem. The proof is complete.

\subsection{Proof of Proposition~\ref{prop:gap_bound}}
\label{appendix:proof_prop_gap_bound}
Let \((\hat{\mathbf x}^*, \hat g^*, \hat\beta^*, \hat t^*)\) and \((\tilde{\mathbf x}^*, \tilde g^*, \tilde\beta^*, \tilde t^*)\) be optimal solutions to problems~\eqref{eq:tightened_problem_wo_Nystrom} and~\eqref{eq:tightened_problem_Nystrom}, respectively. Keeping \(\hat{\mathbf x}^*\) fixed, let \(\bar g = {\rm Proj}_{\mathcal{H}_\mathcal{S}}\hat{g}^*\) and \(\bar t=\|\bar g\|_{\mathcal H_k}\leq \hat t^*\). Define \(\Delta:=\max_{1\leq l\leq N}|\hat{g}^*(\boldsymbol{\zeta}_l)-\bar g(\boldsymbol{\zeta}_l)|\) and a feasibility shift \(\tilde{\beta}=\hat{\beta}^*+\Delta\). For every \(l\in[N]\), feasibility of the original solution indicates
\[
    f(\hat{\mathbf x}^*, \boldsymbol{\zeta}_l)+\omega_{r_N}^f(\hat{\mathbf x}^*) +\omega_{r_N}^k \bar t\leq f(\hat{\mathbf x}^*, \boldsymbol{\zeta}_l)+\omega_{r_N}^f(\hat{\mathbf x}^*) +\omega_{r_N}^k\hat t^*\leq \hat g^*(\boldsymbol{\zeta}_l) + \hat\beta^*\leq \bar g(\boldsymbol{\zeta}_l)+\tilde \beta.
\]
Thus, \((\hat{\mathbf x}^*, \bar g, \tilde\beta, \bar t)\) is feasible for the restricted tightened problem~\eqref{eq:tightened_problem_Nystrom}. Besides, Assumption~\ref{ass:kernel} implies the boundedness of kernel diagonal terms, which in turn controls the feasibility shift:
\begin{equation}\label{ineq:beta_error_bound}
    \begin{aligned}
        \Delta&=\max_{1\leq l\leq N}|\hat{g}^*(\boldsymbol{\zeta}_l)-\bar g(\boldsymbol{\zeta}_l)|=\max_{1\leq l\leq N}\left| \langle \hat{g}^* -\bar g,k(\boldsymbol{\zeta}_l,\cdot) \rangle \right|\\
        &\leq \|\hat{g}^*-\bar g\|_{\mathcal{H}_k}\cdot\max_{1\leq l\leq N}\sqrt{k(\boldsymbol{\zeta}_l, \boldsymbol{\zeta}_l)}\leq \sqrt{C}\|\hat{g}^*-\bar g\|_{\mathcal{H}_k}.
    \end{aligned}
\end{equation}
Then, the optimal-value gap can be further bounded as:
\begin{equation}\label{eq:bound_gap_upperbound}
    \begin{aligned}
        z_{\rm UB+}-z_{\rm UB}&= \left\langle\tilde{g}^*-\hat{g}^*, \sum_{i\in{\rm SM}}\sum_{j=1}^{n_i} \frac{\hat\alpha^*_i}{n_i} \phi(\boldsymbol{\xi}_j^{(i)})\right\rangle_{\mathcal H_k}+\hat{R}^*_\gamma(\tilde t^*-\hat t^*)+(\tilde\beta^*-\hat\beta^*)\\
        &\leq \left\langle\bar g-\hat{g}^*, \sum_{i\in{\rm SM}}\sum_{j=1}^{n_i} \frac{\hat\alpha^*_i}{n_i} \phi(\boldsymbol{\xi}_j^{(i)})\right\rangle_{\mathcal H_k}+\hat{R}^*_\gamma(\bar t-\hat t^*)+\Delta\\
        &\leq \sqrt C\|\bar g-\hat g^*\|_{\mathcal H_k}+\Delta \leq 2\sqrt C\|\hat g^*-\bar g\|_{\mathcal H_k},
    \end{aligned}
\end{equation}
where the second inequality follows from the Cauchy--Schwarz inequality, the bound \(\|\hat\mu_c^*\|_{\mathcal H_k}\leq\sqrt C\), and \(\hat R^*_\gamma(\bar t-\hat t^*)\leq0\), while the last inequality uses~\eqref{ineq:beta_error_bound}.

Now, we turn to quantify \(\|\hat{g}^*-\bar g\|_{\mathcal{H}_k}\). By projecting, \(\bar g\) is the best approximation of \(\hat{g}^*\) in the restricted subspace \(\mathcal{H}_\mathcal{S} = {\rm span}(\{\phi(\boldsymbol\gamma_l)_{l=1}^{N_\mathcal S}\})\), i.e.,
\[
    \bar g = \arg\min_{g\in\mathcal{H}_\mathcal{S}}\|g-\hat{g}^*\|_{\mathcal{H}_k}.
\]
We aim to find \(\tilde{\boldsymbol{\vartheta}}\in\mathbb{R}^{N_{\mathcal S}}\) such that \(\bar g=\sum_{l=1}^{N_{\mathcal S}}\tilde{{\vartheta}}_l\phi(\boldsymbol{\gamma}_l)\), which satisfies the orthogonality condition: 
\begin{equation} 
    \label{eq:orthogonality}
    \left\langle \hat{g}^*-\sum_{l=1}^{N_{\mathcal S}}\tilde{{\vartheta}}_l\phi(\boldsymbol{\gamma}_l), \phi(\boldsymbol{\gamma}_s)\right\rangle_{\mathcal{H}_k}=0, \quad \forall s\in[N_{\mathcal S}].
\end{equation}
According to~\eqref{eq:empirical_representer_g_star}, the optimizer \(\hat{g}^*\) admits the finite expansion \(\hat{g}^*=\sum_{l=1}^{N_\Upsilon}\theta_l^*\phi(\boldsymbol{\upsilon}_l)\), with \(\boldsymbol{\theta}^*\in\mathbb{R}^{N_\Upsilon}\). By the reproducing property, for each \(s\in\{1,\ldots,N_{\mathcal S}\}\),
\begin{align}
    &\langle \hat{g}^*, \phi(\boldsymbol{\gamma}_{s})  \rangle_{\mathcal{H}_k} = \sum_{l=1}^{N_\Upsilon}\theta^*_lk(\boldsymbol{\upsilon}_l, \boldsymbol{\gamma}_{s})= (\mathbf{K}_{\mathcal{S}\Upsilon}\boldsymbol{\theta}^*)_{s}\label{eq:reproducing_ls}.
\end{align}
The orthogonality conditions~\eqref{eq:orthogonality} for the projection then read
\begin{equation}\label{eq:projection_normal_equation}
\mathbf{K}_{\mathcal{S}\mathcal{S}}\tilde{\boldsymbol{\vartheta}}=\mathbf{K}_{\mathcal{S}\Upsilon}\boldsymbol{\theta}^*.
\end{equation}
Note that if \(\mathbf a\in{\rm Null}(\mathbf{K}_{\mathcal{S}\mathcal{S}})\), then \(0=\mathbf a^\top \mathbf{K}_{\mathcal{S}\mathcal{S}}\mathbf a=\|\sum_{s=1}^{N_{\mathcal S}}a_s\phi(\boldsymbol{\gamma}_{s})\|_{\mathcal H_k}^2\), and hence \(\sum_{s=1}^{N_{\mathcal S}}a_s\phi(\boldsymbol{\gamma}_{s})=0\). It follows that \(\mathbf{K}_{\mathcal{S}\Upsilon}^\top\mathbf a=0\), so that \({\rm Null}(\mathbf{K}_{\mathcal{S}\mathcal{S}})\subseteq{\rm Null}(\mathbf{K}_{\mathcal{S}\Upsilon}^\top)\). Taking orthogonal complements gives \({\rm Range}(\mathbf{K}_{\mathcal{S}\Upsilon})\subseteq{\rm Range}(\mathbf{K}_{\mathcal{S}\mathcal{S}})\). Therefore, \(\mathbf{K}_{\mathcal{S}\Upsilon}\boldsymbol{\theta}^*\in{\rm Range}(\mathbf{K}_{\mathcal{S}\mathcal{S}})\), and~\eqref{eq:projection_normal_equation} is consistent. We may choose the solution
\begin{equation}\label{eq:project_coef}
    \boldsymbol{\tilde{\vartheta}} = \mathbf{K}_{\mathcal{S}\mathcal{S}}^{\dagger}\mathbf{K}_{\mathcal{S}\Upsilon}\boldsymbol{\theta}^*,
\end{equation}
which has the minimum Euclidean norm. If \(\mathbf{K}_{\mathcal{S}\mathcal{S}}\) is singular, the coefficient vector need not be unique. Nevertheless, any two solutions differ by an element of \(\operatorname{Null}(\mathbf{K}_{\mathcal{S}\mathcal{S}})\) and therefore represent the
same element of \(\mathcal{H}_{\mathcal{S}}\). Then, we have
\begin{equation} \label{eq:project_error}
    \begin{aligned}
        \|\hat{g}^*-\bar g\|_{\mathcal{H}_k}^2 &= \left\langle\sum_{l=1}^{N_\Upsilon}\theta_l^*\phi(\boldsymbol{\upsilon}_l)-\sum_{s=1}^{N_{\mathcal S}}\tilde{\vartheta}_s\phi(\boldsymbol{\gamma}_s), \sum_{l=1}^{N_\Upsilon}\theta_l^*\phi(\boldsymbol{\upsilon}_l)-\sum_{s=1}^{N_{\mathcal S}}\tilde{\vartheta}_s\phi(\boldsymbol{\gamma}_s)\right\rangle_{\mathcal{H}_k}\\
        &=\boldsymbol{\theta}^{*\top}\mathbf{K}_{\Upsilon\Upsilon}\boldsymbol{\theta}^{*}-2\tilde{\boldsymbol{\vartheta}}^\top\mathbf{K}_{\mathcal{S}\Upsilon}\boldsymbol{\theta}^*+\tilde{\boldsymbol{\vartheta}}^{\top}\mathbf{K}_{\mathcal{S}\mathcal{S}}\tilde{\boldsymbol{\vartheta}}.
    \end{aligned}
\end{equation}
Plugging \eqref{eq:project_coef} into \eqref{eq:project_error}, we obtain
\[
    \|\hat{g}^*-{\rm Proj}_{\mathcal{H}_\mathcal{S}}\hat{g}^*\|_{\mathcal{H}_k}^2= \boldsymbol{\theta}^{*\top}\left(\mathbf{K}_{\Upsilon\Upsilon}-\mathbf{K}_{\mathcal{S}\Upsilon}^\top\mathbf{K}^\dagger_{\mathcal{S}\mathcal{S}}\mathbf{K}_{\mathcal{S}\Upsilon}\right)\boldsymbol{\theta}^*=\boldsymbol{\theta}^{*\top}\mathbf R_{\mathcal S}\boldsymbol{\theta}^{*}.
\]
In particular, \(\mathbf R_{\mathcal S}=\mathbf{K}_{\Upsilon\Upsilon}-\mathbf{K}_{\mathcal{S}\Upsilon}^\top\mathbf{K}^\dagger_{\mathcal{S}\mathcal{S}}\mathbf{K}_{\mathcal{S}\Upsilon}\succeq 0\) because, for every \(\boldsymbol\theta\in\mathbb R^{N_\Upsilon}\), 
\[
    \boldsymbol{\theta}^\top\mathbf R_{\mathcal S}\boldsymbol\theta=\min_{\boldsymbol\vartheta\in\mathbb R^{N_{\mathcal S}}}\left\|\sum_{l=1}^{N_\Upsilon}\theta_l\phi(\boldsymbol\upsilon_l)-\sum_{s=1}^{N_{\mathcal S}}\vartheta_s\phi(\boldsymbol\gamma_s)\right\|_{\mathcal H_k}^2\geq 0.
\]
Equivalently, \(\mathbf R_{\mathcal S}\) is the generalized Schur complement of \(\mathbf{K}_{\mathcal{S}\mathcal{S}}\). Given any symmetric PSD matrix \(\mathbf{K}\), we have \(\boldsymbol{\theta}^\top\mathbf{K}\boldsymbol{\theta}\leq \|\mathbf{K}\|_2\|\boldsymbol{\theta}\|_2^2\). Therefore, 
\[
    \|\hat{g}^*-{\rm Proj}_{\mathcal{H}_\mathcal{S}}\hat{g}^*\|_{\mathcal{H}_k}\leq \|\boldsymbol{\theta}^*\|_2\sqrt{\|\mathbf{K}_{\Upsilon\Upsilon}-\mathbf{K}_{\mathcal{S}\Upsilon}^\top\mathbf{K}^\dagger_{\mathcal{S}\mathcal{S}}\mathbf{K}_{\mathcal{S}\Upsilon}\|_2}=\|\boldsymbol{\theta}^*\|_2\sqrt{\|\mathbf{R}_{\mathcal S}\|_2}.
\]
This concludes the proof.

\section{Experimental Details}
\subsection{Experiments on \emph{Coffee Shop Sales} dataset}
\label{app:coffee_experiment_details}
We use the \emph{Coffee Shop Sales} dataset and consider five brewed-coffee products as five related demand domains. Specifically, the five coffee products considered are Organic brewed coffee with Brazilian (Product 1), Gourmet brewed coffee with Colombian Medium Roast (Product 2), Gourmet brewed coffee with Ethiopia (Product 3), Premium brewed coffee with Jamaican Coffee River (Product 4), and Drip coffee with Our Old Time Diner Blend (Product 5). The following two tables summarize the economic and operational parameters and the settings used for each method, respectively.

\begin{table}[!ht]
\centering
\caption{Economic and operational parameters for the coffee product
launch experiment}
\label{tab:appendix_coffee_parameters}
\renewcommand{\arraystretch}{1.12}
\begin{tabular}{lll}
\toprule
Parameter & Description & Value \\
\midrule
\(d\) & Number of cup sizes & \(3\) \\ 
\(\mathbf{p}\) & Selling prices & \((2.20,\,3.00,\,3.50)\) \\
\(\mathbf{c}\) & Ordering costs, \(\mathbf{c}=0.50\mathbf{p}\) & \((1.10,\,1.50,\,1.75)\) \\
\(\mathbf{h}\) & Holding costs, \(\mathbf{h}=0.15\mathbf{p}\) & \((0.33,\,0.45,\,0.525)\) \\
\(\mathbf{b}\) & Backorder costs, \(\mathbf{b}=0.25\mathbf{p}\) & \((0.55,\,0.75,\,0.875)\) \\
\(\boldsymbol{\kappa}\) & Inventory resource costs & \((1.00,\,1.30,\,1.60)\) \\
\(B_{\mathrm{launch}}\) & Initial procurement budget & \(67\) \\
\(Q_{\max}\) & Total inventory capacity & \(55\) \\
\bottomrule
\end{tabular}
\end{table}

\begin{table}[!ht]
\centering
\caption{Method-specific hyperparameters for the coffee product launch experiment}
\label{tab:real_coffee_hyperparameters}
\renewcommand{\arraystretch}{1.12}
\begin{tabular}{clc}
\toprule
Method & Hyperparameter & Value \\
\midrule
SAA & -- & -- \\
DRO-W & Wasserstein radius \(\epsilon\) & \(0.04\) \\
SAA-WB & -- & -- \\
DRO-K & RKHS-ball radius \(\epsilon\) & \(0.1\) \\
DRO-Conv & -- & -- \\
\midrule
\multirow{5}{*}{RooD-SO}
& Kernel & Gaussian \\
& Kernel bandwidth & Median heuristic \\
& \(\nu\) & \(0.7\) \\
& \(\gamma\) & \(0.001\) \\
& Number of supplementary points \(N\) & \(4000\) \\
& Discretization scheme & Regular grid \\
& Nystr\"om residual tolerance \(\varepsilon\) & 0.1\% \\
\bottomrule
\end{tabular}
\end{table}

\subsection{Portfolio optimization on real-world ETF data}
\label{app:portfolio_optimization_details}
In this subsection, we introduce the details of the ETF selection, model parameters, and method settings used in the portfolio optimization experiment.
\begin{table}[!ht]
\centering
\caption{Three-Sector ETF Combinations with Regional Labels}
\label{tab:sector_etfs}
\resizebox{\linewidth}{!}{%
\begin{tabular}{ccccccccccc}
\toprule
& & Set 1 & Set 2 & Set 3 & Set 4 & Set 5 & Set 6 & Set 7 & Set 8 & Set 9 \\
\midrule
\multirow{3}{*}{Sector}&Utilities & XLU  & VPU & IDU & FUTY & PSCU & FXU & JXI  & WUTI  & XDWU \\
&\makecell{Information\\Technology} & XLK & VGT & IYW & FTEC & PSCT & FXL & IXN  & WTEC  & XDWT \\
&Healthcare & XLV  & VHT & IYH & FHLC & PSCH & FXH & IXJ  & WHEA  & XDWH \\
\midrule
\multirow{2}{*}{Region}&Exposure & U.S. & U.S. & U.S.  & U.S. & U.S. & U.S. & Global & Global  & Global \\
&Listing & U.S.  & U.S.   & U.S.   & U.S.  &  U.S.  & U.S.   & U.S.    & London  & London \\
\bottomrule
\end{tabular}}
\end{table}

Table~\ref{tab:sector_etfs} lists nine ETF sets, each comprising one ETF from each of three sectors: utilities, information technology, and healthcare. Sets 1--6 report U.S. sector proxies constructed from different fund families and methodologies: large-cap sector ETFs (Set 1), broad-market coverage (Sets 2--4), small-cap sector ETFs (Set 5), and strategy-based AlphaDEX sector ETFs (Set 6). Sets 7--9 provide global-exposure sector proxies, including U.S.-listed iShares ETFs (Set 7) and UCITS ETFs listed in London (Sets 8--9). ``Exposure Region'' refers to the geographic exposure of underlying holdings, while ``Listing Region'' indicates the exchange where the ETF is traded.

\begin{table}[!ht]
\centering
\caption{Model parameters for the real-ETF portfolio problem}
\label{tab:real_etf_model_parameters}
\renewcommand{\arraystretch}{1.12}
\begin{tabular}{clc}
\toprule
Parameter & Description & Value \\
\midrule
\(d\) & Number of assets in each portfolio & \(3\) \\
\(\delta_{1}\) & Weight assigned to the CVaR-type risk term & \(10\) \\
\(\delta_{2}\) & Tail probability of the CVaR-type risk term & \(0.2\) \\
\bottomrule
\end{tabular}
\end{table}

The parameters of the portfolio optimization model are summarized in Table~\ref{tab:real_etf_model_parameters}, including the number of assets and the weight and tail probability of the CVaR-type risk term. The following table reports the settings used for each method, including the Wasserstein radius for DRO-W and the kernel, regularization, discretization, and approximation settings for RooD-SO.

\begin{table}[!ht]
\centering
\caption{Method-specific hyperparameters for the real-ETF portfolio experiment}
\label{tab:real_etf_hyperparameters}
\renewcommand{\arraystretch}{1.12}
\begin{tabular}{clc}
\toprule
Method & Hyperparameter & Value \\
\midrule
SAA & -- & -- \\
DRO-W & Wasserstein radius \(\epsilon\) & \(0.007\) \\
SAA-WB & -- & -- \\
DRO-K & RKHS-ball radius \(\epsilon\) & \(0.05\) \\
DRO-Conv & -- & -- \\
\midrule
\multirow{5}{*}{RooD-SO}
& Kernel & Gaussian \\
& Kernel bandwidth & Median heuristic \\
& \(\nu\) & \(0.05\) \\
& \(\gamma\) & \(0.001\) \\
& Number of supplementary points \(N\) & \(1000\) \\
& Discretization scheme & Latin hypercube sampling \\
& Nystr\"om residual tolerance \(\varepsilon\) & 2\% \\
\bottomrule
\end{tabular}
\end{table}

\end{appendices}

\end{document}